\documentclass[10pt,reqno]{amsproc}

\usepackage{amssymb}
\usepackage[mathscr]{eucal}
\usepackage{enumitem}

\usepackage[all]{xy}
\ifx\pdfoutput\undefined \xyoption{dvips}\fi

\newtheorem{thm}{Theorem}
\newtheorem{lemma}[thm]{Lemma}
\newtheorem{cor}[thm]{Corollary}

\theoremstyle{definition}
\newtheorem{defn}[thm]{Definition}
\newtheorem{notation}[thm]{Notation}
\newtheorem{obs}[thm]{Observation}
\newtheorem{recall}[thm]{Recall}
\newtheorem{eg}[thm]{Example}

\def\parhead[#1]{\vspace{1ex}%
\noindent{\bf\boldmath #1.}}

\newcommand{\mlaux}[3]{\setbox0=\hbox{$\mathsurround=0pt #2{#3}$}%
  \dimen0=\dp0\advance\dimen0 by \ht0\lower#1\dimen0\box0}
\newcommand{\mlower}[2]{\mathpalette{\mlaux{#1}}{#2}}

\newcommand{\makekdp}[2]{\setbox0=\hbox{$\mathsurround=0pt #1{#2}$} \dp0=0pt \box0}
\newcommand{\killdepth}{\relax\mathpalette\makekdp}
\newcommand{\category}[1]{\underline{\killdepth{\text{#1}}}}

\newcommand{\makellapm}[2]{\hbox to 0pt{\hss$\mathsurround=0pt #1{#2}$}}

\newcommand{\makerlapm}[2]{\hbox to 0pt{$\mathsurround=0pt #1{#2}$\hss}}

\newcommand{\bicat}[1]{\underline{\mathcal{#1}}}

\newcommand{\icat}{\mathbb}
\newcommand{\lcat}{\mathcal}

\newcommand{\defeq}{\mathrel{\mlower{0.15}{\stackrel{\text{def}}=}}}

\newcommand{\ifundef}[1]{\expandafter\ifx\csname#1\endcsname\relax}

\ifundef{Top}
  \DeclareSymbolFont{dsymbolsC}{U}{txsyc}{m}{n}
  \SetSymbolFont{dsymbolsC}{bold}{U}{txsyc}{bx}{n}
  \DeclareFontSubstitution{U}{txsyc}{m}{n}

  \def\re@DeclareMathSymbol#1#2#3#4{%
      \let#1=\undefined
      \DeclareMathSymbol{#1}{#2}{#3}{#4}}

  \re@DeclareMathSymbol{\Top}{\mathord}{dsymbolsC}{120}
  \re@DeclareMathSymbol{\Bot}{\mathord}{dsymbolsC}{121}
\fi

\newcommand{\partinj}{\Bot}
\newcommand{\partproj}{\Top}

\let\newmop=\DeclareMathOperator

\newmop{\im}{im}
\newmop{\coim}{coim}

\newmop{\cell}{cell}
\newmop{\cof}{cof}
\newmop{\fib}{fib}

\newmop{\triv}{th}
\newmop{\sst}{sp}
\newmop{\sk}{sk}
\newmop{\cosk}{ck}
\newmop{\dec}{dec}
\newmop{\lax}{lax}

\newcommand{\cdec}{\dec^c}
\newcommand{\clax}{\lax^c}
\newcommand{\chom}{\hom^c}

\newcommand{\Hom}{\mathop{\mathcal{H}\text{om}}}

\newmop{\cls}{cls}
\newcommand{\ccls}{\cls^c}

\newmop{\depth}{dp}
\newmop{\id}{id}

\newmop{\colim}{colim}

\newmop{\obj}{obj}
\newmop{\arr}{arr}
\newmop{\sq}{sq}

\newmop{\dis}{dis}

\newmop{\cod}{cod}
\newmop{\dom}{dom}

\newmop{\rev}{rev}

\newcommand{\dual}{^\circ}
\newcommand{\oth}{^{\mathord{\text{th}}}}

\newcommand{\op}{^{\mathord{\text{\rm op}}}}

\newcommand{\boundary}{\partial}

\newcommand{\anytens}{\mathbin{\odot}}
\newcommand{\pretens}{\mathbin{\boxtimes}}
\newcommand{\spretens}{\mathbin{\boxdot}}
\newcommand{\join}{\mathbin{\oplus}}
\newcommand{\gray}{\mathbin{\circledast}}
\newcommand{\laxgray}{\mathbin{\otimes}}

\newcommand{\canytens}{\anytens_c}
\newcommand{\cpretens}{\pretens_c}
\newcommand{\cjoin}{\oplus_c}
\newcommand{\cgray}{\gray_c}
\newcommand{\claxgray}{\laxgray_c}

\newcommand{\sgray}{_{\gray}}

\newcommand{\tens}{\mathbin{\bullet}}

\newcommand{\oriental}{\mathop{\mathcal{O}}}

\newcommand{\indec}{}
\def\indec<#1,#2>{\langle#1,#2\rangle}

\newcommand{\sint}{}
\newcommand{\hint}{}
\def\sint(#1,#2){{(#1,#2)}}
\def\hint(#1,#2]{{(#1,#2]}}

\let\newfunctor=\newmop

\newfunctor{\Dec}{Dec}

\newcommand{\inerve}{\nerv_\infc}

\newfunctor{\pathecat}{P}
\newfunctor{\recons}{R}

\newfunctor{\forget}{U}
\newfunctor{\refl}{L}
\newfunctor{\dgm}{D}
\newfunctor{\ladj}{F}
\newfunctor{\gfunc}{G}
\newfunctor{\incl}{I}
\newfunctor{\nerv}{N}

\newfunctor{\efunc}{E}
\newfunctor{\kfunc}{K}
\newfunctor{\lfunc}{L}

\newcommand{\adjoint}{\dashv}

\newcommand{\yoneda}[1]{\ulcorner{#1}\urcorner}

\newcommand{\joinidt}{\{\ast\}}

\newdir{ >}{{}*!/-10pt/@{>}}
\newdir{u(}{{}*!/-5pt/@^{(}}
\newdir{d(}{{}*!/-5pt/@_{(}}
\newdir{|>}{%
  !/4.5pt/@{|}*:(1,-.2)@^{>}*:(1,+.2)@_{>}*!/-3pt/@{ }}

\def\arrow#1:#2->#3.{{#1\colon #2}\xy (0,0)*+{} %
  \ar (7,0)*+{}\endxy{#3}}

\def\epi#1:#2->#3.{{#1\colon #2}\xy (0,0)*+{} %
  \ar@{->>} (7,0)*+{}\endxy{#3}}

\def\overepi#1:#2->#3.{{#2}\xy (0,0)*+{} %
  \ar@{->>}^-{\smash{#1}} (7,0)*+{}\endxy{#3}}

\def\cover#1:#2->#3.{{#1\colon #2}\xy (0,0)*+{} %
  \ar@{-|>} (7,0)*+{}\endxy{#3}}

\def\overinc#1:#2->#3.{{#2}\xy (0,0)*+{} %
  \ar@{u(->}^-{\smash{#1}} (9,0)*+{}\endxy{#3}}

\def\overarr#1:#2->#3.{{#2}\xy (0,0)*+{} %
  \ar^-{\smash{#1}} (7,0)*+{}\endxy{#3}}

\def\spanarr#1:#2->#3.{{#1\colon #2}\xy (0,0)*+{} %
  \ar|-{\object@{|}} (9,0)*+{}\endxy{#3}}

\def\inc#1:#2->#3.{{#1\colon #2}\xy (0,0)*+{} %
  \ar@{u(->} (9,0)*+{}\endxy{#3}}

\def\nattrans#1:#2->#3.{{#1\colon #2}\xy (0,0)*+{} %
  \ar^-{.} (7,0)*+{}\endxy{#3}}

\def\iso#1:#2->#3.{{#1\colon #2}\xy (0,0)*+{} %
  \ar^-{\smash{\simeq}} (7,0)*+{}\endxy{#3}}

\def\equiv#1:#2->#3.{{#1\colon #2}\xy (0,0)*+{} %
  \ar^-{\smash{\sim}} (7,0)*+{}\endxy{#3}}

\def\anoniso#1->#2.{{#1}\xy (0,0)*+{} %
  \ar^-{\simeq} (7,0)*+{}\endxy{#2}}

\def\twocell#1:#2->#3.{{#1\colon #2}\xy (0,0)*+{} %
  \ar@{=>} (7,0)*+{}\endxy{#3}}

\def\isotwocell#1:#2->#3.{{#1\colon #2}\xy (0,0)*+{} %
  \ar@{=>}^-{\sim} (7,0)*+{}\endxy{#3}}

\newcommand{\Wcs}{\category{Wcs}}

\newcommand{\tDelta}{t\Delta}
\newcommand{\Simp}{\category{Simp}}
\newcommand{\Set}{\category{Set}}

\newcommand{\Cat}{\category{Cat}}

\newcommand{\Strat}{\category{Strat}}

\newcommand{\aSimp}{\Simp_{\mathord{+}}}

\newcommand{\infc}{\ensuremath\omega}

\newcommand{\aDelta}{\Delta_{\mathord{+}}}
\newcommand{\aStrat}{\Strat_{\mathord{+}}}

\newcommand{\CFib}{\category{CFib}}
\newcommand{\TFib}{\category{TFib}}
\newcommand{\WEqv}{\category{WEqv}}

\newcommand{\face}{\delta}
\newcommand{\vertex}{\varepsilon}
\newcommand{\degen}{\sigma}

\newcommand{\cube}[1]{{[\![#1]\!]}}

\newcommand{\usimp}{\tilde}

\newcommand{\pocorner}{\hbox to 10pt{{\vrule height10pt depth0pt width0.5pt}%
    \vbox to 10pt{{\hrule height0.5pt width9.5pt depth0pt}\vfill}}}
\newcommand{\poexcursion}{\save[]-<12pt,-12pt>*{\pocorner}\restore}
\newcommand{\pbcorner}{\vbox to 0pt{\kern 5pt\hbox to 0pt{\kern 5pt%
      \vbox{{\hrule height0.5pt width9.5pt depth0pt}}%
      {\vrule height10pt depth0pt width0.5pt}\hss}\vss}}
\newcommand{\pbexcursion}{\save[]-<-2pt,2pt>*{\pbcorner}\restore}

\def\funcdisplay #1:#2->#3. {{%
    \xymatrix@R=0em@C=14em{%
      {#2}\ar[r]^-{#1} & {#3}}}}

\def\adjdisplay#1->#2,#3-|#4.{{%
    \xymatrix@R=0em@C=7em{%
      {#1} \ar@/_0.65pc/[rr]_-{#4} & {\bot} &
      {#2}\ar@/_0.65pc/[ll]_-{#3}}}}

\def\adjinline#1->#2,#3-|#4.{\arrow {#3}\dashv{#4}:#1->#2.}

\newcommand{\pent}[2]{\xybox{\xymatrix@!=#2{
    & {\c}\ar[rr]^{\cd} && {\d}\ar[rdd]^{\de} & \\
    &&&& \\
    {\b}\ar[uur]^{\bc} &&&& {\e} \\
    && {\a}\ar[ull]^{\ab}\ar[urr]_{\ae} && {\ifcase #1
      \ar"4,3";"1,2"|*+{\labelstyle \ac}="one"
      \ar"4,3";"1,4"|*+{\labelstyle \ad}="two"
      \ar@{}"3,1";"one"|(0.55){\objectstyle \abc}
      \ar@{}@<1.5em>"one";"two"|{\objectstyle \acd}
      \ar@{}"two";"3,5"|(0.45){\objectstyle \ade} \or
      \ar"3,1";"1,4"|*+{\labelstyle \bd}="one"
      \ar"3,1";"3,5"|*+{\labelstyle \be}="two"
      \ar@{}"1,2";"one"|(0.55){\objectstyle \bcd}
      \ar@{}@<1.5em>"one";"two"|{\objectstyle \bde}
      \ar@{}"two";"4,3"|(0.45){\objectstyle \abe} \or
      \ar"4,3";"1,2"|*+{\labelstyle \ac}="two"
      \ar"1,2";"3,5"|*+{\labelstyle \ce}="one"
      \ar@{}"1,4";"one"|(0.55){\objectstyle \cde}
      \ar@{}@<1.5em>"one";"two"|{\objectstyle \ace}
      \ar@{}"two";"3,1"|(0.45){\objectstyle \abc} \or
      \ar"3,1";"1,4"|*+{\labelstyle \bd}="two"
      \ar"4,3";"1,4"|*+{\labelstyle \ad}="one"
      \ar@{}"3,5";"one"|(0.55){\objectstyle \ade}
      \ar@{}@<1.5em>"one";"two"|{\objectstyle \abd}
      \ar@{}"two";"1,2"|(0.45){\objectstyle \bcd} \or
      \ar"3,1";"3,5"|*+{\labelstyle \be}="one"
      \ar"1,2";"3,5"|*+{\labelstyle \ce}="two"
      \ar@{}"4,3";"one"|(0.55){\objectstyle \abe}
      \ar@{}@<1.5em>"one";"two"|{\objectstyle \bce}
      \ar@{}"two";"1,4"|(0.45){\objectstyle \cde} \else\fi} }}}
\newcommand{\pentofpent}[1]{
  \def\baselen{#1}
  \begin{xy}
    0;<\baselen,0mm>:
    *{\xybox{
        \POS(0,-4)*[o]{\pent 0{\baselen}}="zero"
        \POS(16,40)*[o]{\pent 3{\baselen}}="three"
        \POS(72,40)*[o]{\pent 1{\baselen}}="one"
        \POS(88,-4)*[o]{\pent 4{\baselen}}="four"
        \POS(44,-36)*[o]{\pent 2{\baselen}}="two"
        \ar@<1ex>"zero";"three"^-{\objectstyle\abcd}
        \ar@<1ex>"three";"one"^-{\objectstyle\abde}
        \ar@<1ex>"one";"four"^-{\objectstyle\bcde}
        \ar@<-1ex>"zero";"two"_-{\objectstyle\acde}
        \ar@<-1ex>"two";"four"_-{\objectstyle\abce}
        \ar(44,-5);(44,+15)^{\objectstyle\abcde}
     }}
  \end{xy}
}

\title[Weak Complicial Sets I]
      { Weak Complicial Sets \\
        A Simplicial Weak \infc-Category Theory\\
        Part I: Basic Homotopy Theory
}

\author[Verity]{Dominic Verity}
\dedicatory{To Ross Street on the occasion of his $60\oth$ birthday.}
\address{
  Centre of Australian Category Theory \\
  Macquarie University \\
  NSW 2109 \\
  Australia
}
\email{domv@ics.mq.edu.au}

\thanks{The author would like to thank Fitzwilliam College Cambridge
  whose visiting fellowship programme supported the completion of this
  work.}

\date{April 16, 2006}


\begin{document}

\begin{abstract}
  This paper develops the foundations of a simplicial theory of weak
  \infc-categories, which builds upon the insights originally
  expounded by Ross Street in his 1987 paper on oriented
  simplices. The resulting theory of {\em weak complicial sets\/}
  provides a common generalisation of the theories of (strict)
  \infc-categories, Kan complexes and Joyal's {\em
    quasi-categories}. We generalise a number of results due to the
  current author with regard to complicial sets and strict
  \infc-categories to provide an armoury of well behaved technical
  devices, such as joins and Gray tensor products, which will be used
  to study these the weak \infc-category theory of these structures in
  a series of companion papers. In particular, we establish their
  basic homotopy theory by constructing a Quillen model structure on
  the category of {\em stratified simplicial sets\/} whose fibrant
  objects are the weak complicial sets. As a simple corollary of this
  work we provide an independent construction of Joyal's model
  structure on simplicial sets for which the fibrant objects are the
  quasi-categories.
\end{abstract}

\maketitle
\tableofcontents


\section{Overview and History}
\label{intro.sec}

The theory of {\em complicial sets\/} dates to the mid-1970s and to
the work of the mathematical physicist John
Roberts~\cite{Roberts:1978:Complicial}. His original interest in this
topic grew from his conviction that (strict) \infc-categories
constituted the appropriate algebraic structures within which to value
non-abelian cohomology theories~\cite{Roberts:1977:Complicial}. This
led him to define complicial sets to be simplicial sets equipped with
a distinguished set of {\em neutral\/} or {\em thin\/} simplices,
called {\em stratified simplicial sets}, and satisfying a certain kind
of {\em unique thin horn filler\/} condition. He conjectured that it
should be possible to generalise the classical categorical nerve
construction to provide a functor from the category of strict
\infc-categories to the category of complicial sets and this would
provide an equivalence between these two categories.

The first step in realising his vision was made by Ross Street in his
paper on orientals~\cite{Street:1987:Oriental}, which provided the
first fully rigorous description of such a nerve construction and
re-formulated Roberts' vision into a specific conjecture. More
recently, the original program outlined in these papers was completed
by the current author~\cite{Verity:2006:Complicial} who provided the
first proof of the full Street-Roberts conjecture. That work
demonstrates that it is indeed the case that Street's nerve
construction provides the equivalence that Roberts proposed.

This result itself provides a new and powerful approach to studying
strict \infc-categories themselves. For instance,
in~\cite{Verity:2006:Complicial} we show how to construct a
combinatorially simple tensor product of stratified sets whose
reflection into the category of strict \infc-categories is the {\em
  lax Gray tensor product}. Calculations involving this latter
structure are known to be complicated by the fact that it is usually
presented as a colimit of strict \infc-categories freely generated by
geometric products of {\em globs}. On the other hand, if we are
willing to work in the world of stratified simplicial sets then we may
instead describe the corresponding tensor directly by equipping the
product of underlying simplicial sets with a suitably defined set of
thin simplices. In contrast this latter structure is eminently well
suited to direct combinatorial calculation.

However, it was not this kind of application to strict \infc-category
theory which originally encouraged the current author's interest in
the Street-Roberts conjecture. Instead it was piqued by parenthetic
comments in Street's paper on oriented
simplices
~\cite{Street:1987:Oriental} to the effect that we might use
it as a foundation upon which to develop a useful generalisation of
the theory of {\em bicategories\/} to higher dimensions. At that time
no truly workable theory of such {\em weak \infc-categories\/} had
been constructed, although a growing group of researchers were
becoming aware of the role that such structures might play in
algebraic topology, theoretical physics, computer science and higher
category theory itself. In brief, Street's idea was that we might
obtain such a theory by again working with stratified simplicial sets
but this time weakening the axioms that characterised complicial sets
by only insisting that horns of the kind identified in that theory
should have {\em some\/}, not necessarily unique, thin filler.

The subsequent 20 years has been a fertile one in weak \infc-category
theory and we might now identify in the literature three or four quite
distinct approaches to defining such structures, each of which splits
into a plethora of definitional sub-varieties. In this time Street's
original side remark has remained largely under investigated, indeed
for a number of years the current author has avoided writing up his
own ideas along these lines for fear of simply launching yet another
weakened higher category definition on the world. Spurred on, however,
by Street's 2003 work on weak
\infc-categories~\cite{Street:2003:WomCats}, which reformulated and
refined his original insight and introduced the term {\em weak
  complicial set\/} and Joyal's work on {\em
  quasi-categories\/}~\cite{Joyal:2002:QuasiCategories}, I am now of
the view that a thorough explication of this theory is well overdue.

So why might we be interested in studying weak \infc-category theories
based upon simplicial rather than globular geometries? From a
philosophical perspective, Street himself sums the case up best in the
following passage from~\cite{Street:2003:WomCats}:
\begin{quotation}
  \em Simplicial sets are lovely objects about which algebraic
  topologists know a lot. If something is described as a simplicial
  set, it is ready to be absorbed into topology. Or, in other words,
  no matter which definition of weak \infc-category eventually becomes
  dominant, it will be valuable to know its simplicial nerve.
\end{quotation}
In short, any weak \infc-category theory worth its salt should come
equipped with a simplicial nerve functor describing its place in
algebraic topology. Furthermore, it is reasonable to expect that this
would, at the very least, map each weak \infc-category to a weak
complicial set. It follows, therefore, that any study of weak
complicial sets themselves will remain valuable regardless of which
particular formulation of the weak \infc-category notion might become
dominant in the future.

More pragmatically, the answer to this question is really one of
utility. As we shall see here the theory of weak complicial sets is
one which immediately generalises the most widely accepted $0$-trivial
and $1$-trivial weak \infc-categorical structures (Kan sets and
quasi-categories respectively) and at the same time encompasses the
theory of strict \infc-categories. Furthermore, we shall also
demonstrate here that it supports a plethora of well behaved
\infc-categorical constructions, such as joins
(section~\ref{joins.sec}) and Gray tensor products
(section~\ref{gray.tensor.sec}), and admits a rich homotopy theory
(section~\ref{quillen.model.sec}). In a companion
paper~\cite{Verity:2005:GrayNerves}, we derive a nerve construction
for categories enriched, in the classical sense
of~\cite{Kelly:1982:ECT}, in weak complicial sets with respect to a
Gray tensor product (called {\em complicial Gray-categories}) which
faithfully represents such structures as weak complicial sets. In
particular, this demonstrates that the totality of all weak complicial
sets and their homomorphisms, strong transformations, modifications
and so forth is itself representable as a very richly structured
(large) weak complicial set.

The actual category theory of these structures will be explored in
another companion paper~\cite{Verity:2006:IntQCat}, wherein we
represent weak complicial sets as certain kinds of {\em complicially
  enriched quasi-categories}. This provides us with a natural context
in which to generalise traditional category theory to a kind of
homotopy coherent quasi-category theory within the Quillen model
category of weak complicial sets itself. This approach allows us to
translate all of the basic constructions of $n$-category theory into
the weak complicial context and at the same time to establish for it
homotopical versions of the theories of discrete fibrations, Yoneda's
lemma, adjunctions, limits and colimits and so forth.

While all of this speaks to the expressiveness of weak complicial set
theory, we must also convince ourselves that it provides a strong
enough framework within which to establish certain natural coherence
theorems. While work in this direction is still at a relatively early
stage, studies to date indicate that it is likely that a direct
analogue of the well-known coherence theorems for bicategories and
tricategories holds in this context as well. To be precise, there are
strong reasons to suspect that every weak complicial set satisfying
certain very mild conditions on its thin $1$-simplices (related to our
work here in section~\ref{equiv.sec}) is equivalent to the nerve of
some complicial Gray-category.

Herein, however, we restrict ourselves to the modest task of
establishing the foundational homotopy theory of weak complicial sets
upon which all of our later work in this area will be
based. Section~\ref{wcs.intro.sec} introduces these structures and
establishes the associated theory of anodyne extensions between
stratified simplicial sets. Section~\ref{joins.sec} studies the join
operation on stratified sets, introduces the corresponding
d{\'e}calage construction and demonstrates that these are
appropriately well behaved with respect to weak compliciality. This
work is then applied in section~\ref{equiv.sec} to show that we can
usefully replace the condition which stipulates that weak complicial
sets must have thin fillers for outer complicial horns with one which
simply states that all thin $1$-simplices are actually equivalences in
some suitable sense. 

Section~\ref{gray.tensor.sec}, which is combinatorially the most
involved of this work, re-introduces the (lax) Gray tensor products
of~\cite{Verity:2006:Complicial} and studies their properties with
regard weak complicial sets and anodyne extensions. In particular,
this allows us to show that we may construct weak complicial sets of
homomorphisms, (lax) transformations, (lax) modifications and so forth
between any pair of weak complicial sets and thereby enrich the
category of these structures over itself in three distinct but related
ways. Subsection~\ref{str.comp.char} provides a new characterisation
of strict complicial sets as those weak ones which are well-tempered,
in the sense that for these thinness is a sufficient property for the
detection of degenerate simplices.

Finally section~\ref{quillen.model.sec} draws together these various
threads by constructing a Quillen model structure on the category of
stratified simplicial sets whose cofibrations are the inclusions and
whose fibrant objects are precisely the weak complicial
sets. Furthermore, we show that this is a monoidal model category with
respect to the Gray tensor products studied in
section~\ref{gray.tensor.sec}. Finally, we round out our presentation
by localising our model structure and transporting it to the category
of simplicial sets itself, in order to provide an independent
construction of a model category structure on that latter category
whose fibrant objects are Joyal's
quasi-categories~\cite{Joyal:2002:QuasiCategories}.



\section{Introducing Weak Complicial Sets}
\label{wcs.intro.sec}

Here we recall the standard notation of the theory of simplicial sets,
introduce their {\em stratified\/} generalisations and establish the
basic machinery required to define and study {\em weak complicial
  sets}.

\subsection{Stratified Simplicial Sets}

\begin{notation}[simplicial operators]
  As usual we let $\aDelta$ denote the (skeletal) category of finite
  ordinals and order preserving maps between them and use the notation
  $\Delta$ to denote its full subcategory of non-zero
  ordinals. Following tradition we let $[n]$ denote the ordinal $n+1$
  as an object of $\aDelta$ and refer to arrows of $\aDelta$ as {\em
    simplicial operators}.  We will generally use lower case Greek
  letters $\arrow\alpha,\beta,\gamma...:[m]->[n].$ to denote
  simplicial operators and let $\im(\alpha)$ denote the subset
  $\{i\in[n]\mid \exists j\in[m]\mathrel{.}
  \alpha(j)=i\}\subseteq[n]$ known as the {\em image\/} of the
  operator $\alpha$. We will also use the following standard notation
  and nomenclature throughout:
  \begin{itemize}
  \item The injective maps in $\aDelta$ are referred to as {\em face
      operators}. For each $j\in[n]$ we use the
    $\arrow\face^n_j:[n-1]->[n].$ to denote the {\em elementary\/}
    face operator distinguished by the fact that its image does
    not contain the integer $j$.
  \item The surjective maps in $\aDelta$ are referred to as {\em
      degeneracy operators}. For each $j\in[n]$ we use
    $\arrow\degen^n_j:[n+1]->[n].$ to denote the {\em elementary\/}
    degeneracy operator determined by the property that two integers
    in its domain map to the integer $j$ in its codomain.
  \item For each $i\in[n]$ the operator $\arrow\vertex^n_i:[0]->[n].$
    given by $\vertex^n_i(0)=i$ is called the {\em $i\oth$
      vertex operator of $[n]$}. 
  \item We also use the notations $\arrow\eta^n:[n]->[0].$ and
    $\arrow\iota^n:[-1]->[n].$ to denote the unique such simplicial
    operators.
  \end{itemize}
  Unless doing so would introduce an ambiguity, we will tend to reduce
  notational clutter by dropping the superscripts of these elementary
  operators.
\end{notation}

\begin{notation}[simplicial sets]
  \label{stand.simp.def}
  The category $\Simp$ of {\em simplicial sets\/} and {\em simplicial
    maps\/} between them is simply the functor category
  $[\Delta\op,\Set]$, where $\Set$ denotes the (large) category of all
  (small) sets and functions between them.  If $\arrow
  X:\Delta\op->\Set.$ is a simplicial set then we will often simplify
  our notation by using $X_n$ for the object $X([n])\in\lcat{M}$
  and $X_\alpha$ for the function $\arrow X(\alpha):
  X([m])->X([n]).$. We also adopt the standard latin notations
  $d^n_i$, $s^n_i$ and $v^n_i$ for the actions of the elementary
  simplicial operators $\face^n_i$, $\degen^n_i$ and $\vertex^n_i$
  respectively.

  In practice, it is often easier to think of a simplicial set as a
  single set endowed with a partially defined right action of the
  simplicial operators.  To be more precise, this description presents
  a simplicial set as a triple consisting of a single set $X$, a
  dimension function $\arrow\dim:X->\mathbb{N}.$, and a partial right
  action $x\cdot\alpha$ of $\alpha\in\arr(\Delta)$ on $x\in X$ which
  is defined whenever the dimension of $x\in X$ equals that of the
  codomain of $\alpha$. Under this presentation, simplicial maps
  become functions of underlying sets which preserve both dimension
  and action.  We say that $X$ is a {\em simplicial subset\/} of a
  simplicial set $Y$, denoted $X\subseteq_s Y$, if $X$ is a subset of
  $Y$ which is closed in there under the action of simplicial
  operators and thus inherits a simplicial set structure from it.
  We adopt the following traditional denotations of a few fundamental
  simplicial sets:
  \begin{itemize}
  \item The {\em standard $n$-simplex\/} $\Delta[n]$ which is the
    representable simplicial set on $[n]$, whose $r$-simplices are
    operators $\arrow\alpha:[r]->[n].\in\Delta$ acted upon by right
    composition.
  \item The {\em boundary of the $n$-simplex\/}
    $\boundary\Delta[n]$ which is the simplicial subset of $\Delta[n]$
    of those simplices $\arrow\alpha:[r]->[n].$ which are {\bf not}
    surjective. Notice that the boundary of the $0$-simplex
    $\boundary\Delta[0]$ is simply the empty stratified set $\emptyset$.
  \item The {\em $(n-1)$-dimensional $k$-horn\/} $\Lambda^k[n]$ which is the
    simplicial subset of $\Delta[n]$ consisting of those simplices
    $\arrow\alpha:[r]->[n].$ for which there is some $i\in[n]$ which
    is neither in the image of $\alpha$ nor equal to $k$ (that is for
    which $[n]\neq \im(\alpha)\cup\{k\}$). In other words, this is the
    smallest simplicial subset of $\Delta[n]$ containing the
    set of $(n-1)$-faces $\{ \degen^n_i:i\in [n]\setminus\{k\} \}$. 
  \end{itemize}

  We say that a simplex $x$ of a simplicial set $X$ is {\em
    degenerate\/} iff there is some non-identity degeneracy operator
  $\alpha$ and a simplex $x'\in X$ such that $x=x'\cdot\alpha$. More
  specifically we say that $x$ is {\em degenerate at $k$\/} if
  $x=x'\cdot\degen_k$ for some simplex $x'\in X$, in which case we
  would have $x'=x\cdot\face_k=x\cdot\face_{k+1}$. The {\em
    Eilenberg-Zilber lemma\/} tells us that every simplex $y\in X$ may
  be represented uniquely as $y=x\cdot\beta$ where $\beta$ is a
  degeneracy operator and $x$ is a {\bf non-}degenerate simplex.

  Finally, recall that {\em Yoneda's lemma\/} for simplicial sets
  tells us that there exists a natural bijection between the
  $n$-simplices of a simplicial set $X$ and simplicial maps
  $\overarr:\Delta[n]->X.$. This identifies $x\in X_n$ with the
  simplicial map $\yoneda{x}$ that carries the simplex
  $\alpha\in\Delta[n]$, which is simply a simplicial operator with
  codomain $[n]$, to the simplex $\yoneda{x}(\alpha)\defeq
  x\cdot\alpha$ in $X$.
\end{notation}

\begin{notation}\label{simp.of.1}
  We introduce the following notations to denote the simplices of the
  standard simplex $\Delta[1]$:
  \begin{itemize}
  \item $\arrow 0^r:[r]->[1].$ is the operator which maps each
    $i\in[r]$ to $0\in[1]$.
  \item $\arrow 1^r:[r]->[1].$ is the operator which maps each
    $i\in[r]$ to $1\in[1]$.
  \item $\arrow\rho^r_i:[r]->[1].$ ($1\leq i\leq r$) is the operator
    defined by
    \begin{equation*}
      \rho^r_i(j) =
      \begin{cases}
        0 & \text{if $j< i$,} \\
        1 & \text{if $j\geq i$.}
      \end{cases}
    \end{equation*}
  \end{itemize}
  As above, we shall adopt the convention of omitting the superscripts
  on these operators unless doing so would introduce an ambiguity.
  Later on it will become convenient to index the $r$-simplices of
  $\Delta[1]$ using the {\em doubly pointed\/} set
  $\cube{r}\defeq\{-,+,1,2,...,r\}$, by letting $\rho^r_{-}=0^r$,
  $\rho^r_{+}=1^r$ and defining $\rho^r_i$ as above for an arbitrary
  integer (non-point) in $\cube{r}$.
\end{notation}

\begin{obs}[nerves of categories]\label{cat.nerves}
  We shall also assume that the reader is familiar with the classical
  {\em nerve construction\/} which functorially associates a
  simplicial set $\nerv(\icat{C})$ to each category $\icat{C}$.  This is
  formed by regarding the ordered sets $[n]$ to be categories in the
  usual way and applying {\em Kan's
    construction\/}~\cite{Kan:1958:Adjoint} to the inclusion of
  $\Delta$ as a full subcategory into $\Cat$ (the category of small
  categories), to obtain an adjoint pair:
  \begin{equation*}
    \adjdisplay \Cat->\Simp, \ladj-| \nerv.
  \end{equation*}
  In other words, the $n$-simplices of $\nerv(\icat{C})$ are functors
  $\arrow f:[n]->\icat{C}.$ upon which simplicial operators act by
  pre-composition.
\end{obs}

\begin{defn}[stratified simplicial sets]\label{strat.sets.def}
  A {\em stratification\/} on a simplicial set $X$ is a
  subset\footnote{Note that $tX$ is merely a subset of $X$, {\bf not}
    a simplicial subset, in general it will not be closed in $X$ under
    the action of simplicial operators.}  $tX$ of its simplices
  satisfying the conditions that
  \begin{itemize}
  \item no $0$-simplex of $X$ is in $tX$, and
  \item all of the degenerate simplices of $X$ are in $tX$.
  \end{itemize}

  A {\em stratified set\/} is a pair $(X,tX)$ consisting of a
  simplicial set $X$ and a chosen stratification $tX$ the elements of
  which we call {\em thin\/} simplices. In practice, we
  will elect to notionally confuse stratified sets with their
  underlying simplicial sets $X,Y,Z,...$ and uniformly adopt the
  notation $tX,tY,tZ,...$ for corresponding sets of thin
  simplices. Then, where disambiguation is required, we use the notation
  $\usimp{X},\usimp{Y},\usimp{Z},...$ to denote the
  underlying simplicial sets of these stratified sets.

  A {\em stratified map\/} $\arrow f:X->Y.$ is simply a simplicial map
  of underlying simplicial sets which {\em preserves thinness\/} in
  the sense that for all $x\in tX$ we have $f(x)\in tY$. Identities
  and composites of stratified maps are clearly stratified maps, from
  which it follows that we have a category $\Strat$ of stratified sets
  and maps.
\end{defn}

\begin{defn}[stratified subsets, inverse and direct images]
  Suppose that $U$ and $X$ are stratified sets, then we say that $U$
  is a {\em stratified subset\/} of $X$, denoted $U\subseteq_s X$, if
  $\usimp{U}$ is a simplicial subset of $\usimp{X}$ and its
  stratification $tU$ is a subset of $tX$. If $\arrow f:X->Y.$ is a
  stratified map then the:
  \begin{itemize}
  \item {\em direct image\/} of the stratified subset $U\subseteq_s X$
    along $f$ is the stratified subset $f(U)\subseteq_s Y$ with
    underlying simplicial set $\{f(x)\mid x\in U\}$ and in which $y\in
    f(U)$ is thin iff there is some $x\in tU$ with $f(x)=y$.
  \item {\em inverse image\/} of the stratified subset $V\subseteq_s
    Y$ along $f$ is the stratified subset $f^{-1}(V)\subseteq_s X$
    with underlying simplicial set $\{x\in X\mid f(x)\in V\}$ and in
    which $x\in f^{-1}(V)$ is thin iff $f(x)$ is thin in $V$.
  \end{itemize}
\end{defn}

\begin{obs}[inclusions of stratified sets]\label{strat.inc.obs}
  We call the monomorphisms in $\Strat$ {\em stratified inclusions\/}
  and these are customarily denoted by arrows with hooked domains
  $\inc i:X->Y.$. A stratified subset $X$ of $Y$ clearly gives rise to
  a corresponding stratified inclusion which we denote by
  $\overinc\subseteq_s:X->Y.$. Indeed, wherever necessary we may
  always replace an arbitrary stratified inclusion by an isomorphic
  subset inclusion.

  The forgetful functor from $\Strat$ to $\Set$ which carries a
  stratified set to its set of simplices preserves colimits and
  reflects monomorphisms. It follows that the class of stratified
  inclusions is closed in $\Strat$ under pushout, transfinite
  composition and retraction since this is the case for the class of
  injective functions in $\Set$. Furthermore the class of all
  stratified inclusions is the cellular completion of the set of
  boundary and thin simplex inclusions:
  \begin{equation*}
    \{\overinc\subseteq_r:\boundary\Delta[n]->\Delta[n].\mid
    n=0,1,...\}\cup\{\overinc\subseteq_e:\Delta[n]->\Delta[n]_t.\mid 
    n=1,2,...\}
  \end{equation*}
\end{obs}

\begin{defn}[stratified subsets, regularity and entirety]
  \label{reg.ent.def}
  We say that a stratified subset $X$ of $Y$ is:
  \begin{itemize}
  \item {\em regular}, denoted $X\subseteq_r Y$, if
    $tX=\usimp{X}\cap tY$, and
  \item {\em entire}, denoted \ $X\subseteq_e Y$, if
    $\usimp{X}=\usimp{Y}$.
  \end{itemize}
  The terms {\em regular subset\/} and {\em entire subset\/} will
  always be taken to denote stratified subsets which possess the
  appropriate property. If $W$ is a subset of simplices of the
  stratified set $X$ then the stratified (resp.\ regular or entire)
  subset of $X$ {\em generated\/} by $W$ is defined to be the
  smallest such stratified subset of $X$ which contains $W$.

  Extending these definitions to all stratified maps, we say that
  $\arrow f:X->Y.\in\Strat$ is regular if it {\em reflects thin
    simplices}, meaning that whenever $f(x)$ is thin in $Y$ it follows
  that $x$ is thin in $X$, and entire if it is surjective on
  simplices.

  A stratified map $\arrow f:X->Y.$ admits two well-behaved canonical
  factorisations:
  \begin{itemize}
  \item {\em regular image factorisation\/} $\overarr
    f_e:X->\im_r(f).\overinc\subseteq_r:->Y.$ in which the stratified map
    $f_e$ is entire and $\im_r(f)$, the {\em regular image\/} of $f$,
    is the regular subset of $Y$ whose set of simplices is
    $\{y\in Y\mid \exists x\in X\mathrel{.} f(x)=y\}$.
  \item {\em entire coimage factorisation\/} $\overinc
    \subseteq_e:X->\coim_e(f).\overarr f_r:->Y.$ in which the stratified
    map $f_r$ is regular and $\coim_e(f)$, the {\em entire coimage\/}
    of $f$, is the entire superset of $X$ whose thin
    simplices are those $x\in X$ for which $f(x)$ is thin in $Y$.
  \end{itemize}
\end{defn}

\begin{notation}[complicial simplices and horns]
  \label{comp.simp.def}
  The functor from $\Strat$ to $\Simp$ which forgets stratifications
  has both a left and a right adjoint, which assign to a simplicial
  set its minimal and maximal stratification respectively. We will
  implicitly promote any simplicial set $X\in\Simp$ to a stratified
  set using the (left adjoint) minimal stratification, under which its
  sets of degenerate and thin simplices coincide, and thereby regard
  $\Simp$ as a full subcategory of $\Strat$. In particular, the
  representable simplicial sets $\Delta[n]\in\Simp$ provide us with
  geometrical models for the standard simplices in $\Strat$.

  A few other stratified sets will take on particular importance in
  our deliberations later on:
  \begin{itemize}
  \item The {\em standard thin $n$-simplex\/} $\Delta[n]_t$
    constructed from $\Delta[n]$ by making thin its unique
    non-degenerate $n$-simplex
    $\arrow\id_{[n]}:[n]->[n].\in\Delta[n]$.
  \item The {\em $k$-complicial $n$-simplex \/} $\Delta^k[n]$
    constructed from $\Delta[n]$ by making thin all those simplices
    $\arrow\alpha:[r]->[n].$ whose image contains the set of integers
    $\{k-1,k,k+1\}\cap[n]$. Non-degenerate simplices satisfying this
    latter condition are said to be {\em $k$-admissible}.
  \item The {\em $(n-1)$-dimensional $k$-complicial horn\/}
    $\Lambda^k[n]$ which is the regular subset of
    $\Delta^k[n]$ of those simplices $\arrow\alpha:[r]->[n].$ for
    which the set $[n]\setminus(\im(\alpha)\cup\{k\})$ is
    non-empty. In other words, this is the regular subset
    of $\Delta^k[n]$ generated by its set $\{\face_i\mid
    i\in[n]\setminus\{k\}\}$ of all $(n-1)$-faces except $\face_k$.
  \item The stratified set $\Delta^k[n]''$ and its regular subset
    $\Lambda^k[n]'$ which are obtained from $\Delta^k[n]$ and
    $\Lambda^k[n]$ (respectively) by making all 
    $(n-1)$-simplices thin. 
  \item The union $\Delta^k[n]'\defeq\Delta^k[n]\cup
    \Lambda^k[n]'\subseteq_e\Delta^k[n]''$ which may be
    constructed from $\Delta^k[n]$ by making thin the
    $(n-1)$-simplices $\face^n_{k-1}$ and $\face^n_{k+1}$.
  \end{itemize}
  While the stratifications of these complicial simplices may seem a
  little less than intuitive, they are however fundamental to much of
  the theory that follows. Motivation for these choices is provided by
  the various works of Roberts~\cite{Roberts:1977:Complicial} 
  and~\cite{Roberts:1978:Complicial}, Street \cite{Street:1987:Oriental}
  and~\cite{Street:1988:Fillers} and Verity \cite{Verity:2006:Complicial}.
\end{notation}

\begin{obs}[$k$-admissibility recast]\label{kadmis}
  It is sometime useful to recast the definition of $k$-admissibility
  slightly. To that end, it is easily shown that a non-degenerate
  $r$-simplex $\alpha\in\Delta[n]$ is $k$-admissible if and only if
  there exists some $l\in[r]$ such that $\alpha(i)=k+i-l$ for each
  $i\in[r]\cap\{l-1,l,l+1\}$.
\end{obs}

\begin{obs}[$\Strat$ as a LFP quasi-topos]\label{LFP.quasitopos}
  The full subcategory $\tDelta$ of standard simplices and standard
  thin simplices is dense in $\Strat$ (cf.\ chapter 5 of
  Kelly~\cite{Kelly:1982:ECT}), thereby providing us with a reflective
  full embedding of $\Strat$ into the presheaf category
  $[\tDelta\op,\Set]$. More explicitly, $\tDelta$ may be obtained from
  $\Delta$ by appending extra objects $[n]_t$ for $n=1,2,...$ and
  extra operators $\arrow\varsigma^n_k:[n+1]_t->[n].$ and
  $\arrow\varphi^n:[n]->[n]_t.$ satisfying the relations
  $\varsigma^n_k\circ\varphi^{n+1}=\degen^n_k$. A presheaf
  $F\in[\tDelta\op,\Set]$ is isomorphic to some stratified set if and
  only if it maps each operator $\arrow\varphi^n:[n]->[n]_t.$ to a
  monomorphism in $\Set$. It follows that the category $\Strat$ is
  locally finitely presentable, since it is equivalent to the category
  of models for a finite limit sketch on $\tDelta$.

  The utility of this observation is immediately clear, for instance
  it tells us that $\Strat$ has limits which are calculated pointwise,
  colimits which are constructed in $[\tDelta\op,\Set]$ and then
  reflected into $\Strat$ and that its finitely presented objects are
  those stratified sets with only a finite number of non-degenerate
  simplices. Furthermore, as observed by Street
  in~\cite{Street:2003:WomCats}, the left adjoint to the inclusion
  $\overinc:\Strat->[\tDelta\op,\Set].$ preserves pullbacks of pairs
  of morphisms into (images of) stratified sets from which it follows
  that $\Strat$ is a quasi-topos. In other words, for each stratified
  set $X$ the slice category $\Strat/X$ is cartesian closed and
  $\Strat$ has a classifier for regular subobjects.
\end{obs}

\begin{notation}[skeleta and superstructures]\label{sup.defn}
  We say that a stratified set is {\em $n$-skeletal\/} if all of its
  simplices of dimension greater than $n\in\mathbb{N}$ are
  degenerate. The {\em $n$-skeleton\/} $\sk_n(X)$ of a stratified set
  $X$ is its regular subset consisting of those of its simplices whose
  faces of dimension greater than $n$ are all degenerate. This
  construction provides us with an endo-functor of $\Strat$ whose
  range is the full subcategory of $n$-skeletal stratified sets and
  which has a right adjoint $\cosk_n$ called the {\em
    $n$-coskeleton\/} functor.

  Playing the same game with thinness, we say that a stratified set is
  {\em $n$-trivial\/} if all of its simplices of dimension greater
  than $n$ are thin. The {\em $n$-trivialisation\/} $\triv_n(X)$ of a
  stratified set $X$ is its entire superset constructed by making thin
  all of its simplices of dimension greater than $n$. Again this
  construction provides us with an endo-functor of $\Strat$ whose
  range is the full-subcategory of $n$-trivial stratified sets and
  which has a right adjoint $\sst_n$ called the {\em
    $n$-superstructure\/} functor. The $n$-superstructure $\sst_n(X)$
  may be realised as the regular subset of $X$ of those simplices
  whose faces of dimension greater than $n$ are all thin. 
\end{notation}

\subsection{Weak Complicial Sets}

Now we are ready to embark on defining and studying weak complicial
sets: 

\begin{notation}[lifting problems and properties]
  A commutative square in some category $\lcat{C}$
  \begin{equation*}
    \xymatrix@=1.8em{
      {U}\ar[r]^{u}\ar[d]_{i} & {E}\ar[d]^{p} \\
      {V}\ar[r]_{v}\ar@{..>}[ur] & {A} }
  \end{equation*}
  is called a {\em lifting problem\/} from $i$ to $p$ and it is said
  to {\em have a solution\/} if there exists some diagonal
  map (dotted in the diagram) which makes both triangles commute.
  When such a solution exists we say that $i$ has the {\em left
    lifting property (LLP)\/} with respect to $p$ or that $p$ had the
  {\em right lifting property (RLP)\/} with respect to $i$.

  We say that an object $C\in\lcat{C}$ has the RLP with respect to
  the morphism $\arrow i:U->V.$ iff the unique map $\arrow:C->1.$ into
  the terminal object of $\lcat{C}$ enjoys that property. In such a
  case, a lifting problem amounts to a morphism $\arrow u:U->C.$ and
  a solution to this is simply a morphism $\arrow\bar{u}:V->C.$ for
  which $\bar{u}\circ i=u$.
\end{notation}  

\begin{defn}[elementary anodyne extensions and weak complicial sets]
  \label{wcs.defn}
  The set of {\em elementary anodyne extensions\/} in $\Strat$
  consists of two families of subset inclusions:
  \begin{itemize}
  \item $\overinc\subseteq_r:\Lambda^k[n]->\Delta^k[n].$ for
    $n=1,2,...$ and $k\in[n]$, these are called {\em complicial horn
      extensions}, and
  \item $\overinc\subseteq_e:\Delta^k[n]'->\Delta^k[n]''.$ for
    $n=2,3,...$ and $k\in[n]$, these are called {\em complicial
      thinness extensions}.
  \end{itemize}
  We classify these elementary anodyne extensions into two
  sub-classes, the {\em inner\/} ones for which the index $k$
  satisfies $0<k<n$ and the remaining {\em left\/} and {\em right
    outer\/} ones for which $k=0$ or $k=n$ respectively. 
  Now we say that a stratified set $A$ is a
  \begin{itemize}
  \item {\em weak inner complicial set\/} if it has the RLP with
    respect to all inner elementary anodyne extensions.
  \item {\em weak left (resp.\ right) complicial set\/} if it is a
    weak inner complicial set which also has the RLP with respect to
    all left (resp.\ right) outer elementary anodyne extensions.
  \item {\em weak complicial set\/} if it has the RLP with respect to
    all elementary anodyne extensions.
  \end{itemize}
  Informally we might simply say that a weak complicial set has {\em
    fillers\/} for all complicial horns.
\end{defn}

\begin{eg}[Kan complexes and Joyal's quasi-categories]\label{kan.joyal.eg}
  The theory of weak complicial sets generalises and subsumes those of
  Kan complexes and Joyal's quasi-categories. In particular, 
  if $X$ is a simplicial set then it is:
  \begin{itemize}
  \item a Kan complex iff $\triv_0(X)$ is a weak complicial set,
    and
  \item a quasi-category iff it admits some $1$-trivial stratification 
    which makes it into a weak complicial set.
  \end{itemize}
  The first of these observations is trivial, the second is a direct
  consequence of theorem 1.3 in Joyal's paper on
  quasi-categories~\cite{Joyal:2002:QuasiCategories}. We return to
  this example in section~\ref{equiv.sec}, where we generalise and
  reprove Joyal's result in the current context.
\end{eg}

\begin{eg}[complicial sets]\label{comp.sets}
  Definitions 121 and 154 of \cite{Verity:2006:Complicial} tell us
  that any {\em complicial set\/} satisfies a {\em unique\/} horn
  filler condition with respect to elementary inner anodyne
  extensions.  Furthermore, lemma 163 of {\em loc.\ cit.\/}
  demonstrates that any left (resp. right) outer complicial
  $n$-simplex in a complicial set is degenerate at $0$ (resp.\
  $n-1$). From this fact it is easily demonstrated that in a
  complicial set any outer complicial horn may be (uniquely) filled by
  a degenerate simplex. It follows that any complicial set is actually
  a weak complicial set. A converse to this result, providing an
  alternative characterisation of complicial sets amongst the weak
  complicial sets, may be found in theorem~\ref{wcs+welltemp=>scs}.

  We will sometimes say that the complicial sets of
  \cite{Verity:2006:Complicial} are {\em strict\/} in order to
  differentiate them more clearly from the far more general weak
  complicial sets of this paper.
\end{eg}

\begin{eg}[stratifying \infc-categorical nerves]
  \label{estrat.nerve}
  The combinatorial calculations of Street~\cite{Street:1988:Fillers}
  demonstrate that the nerve $\inerve(\icat{C})$ of any (strict)
  \infc-category $\icat{C}$ may be made into a (strict) complicial set
  by endowing it with the {\em Roberts stratification\/} in which the
  commutative simplices are thin. However, the same calculations may
  be pushed a little further to show that $\inerve(\icat{C})$ is also
  made into a (generally non-strict) weak complicial set, denoted
  $\inerve^e(\icat{C})$, by endowing it with the stratification under
  which a $n$-simplex is thin if it maps the unique non-trivial
  $n$-cell of the $n\oth$ oriental $\oriental_n$ to an {\em
    \infc-categorical $n$-equivalence\/} in $\mathbb{C}$. The precise
  formulation and proof of this fact, which we shall not require
  further here, is a matter of routine (strict) \infc-category theory,
  which we leave as an exercise to the reader.
\end{eg}

\begin{eg}[nerves of complicial Gray-categories] 
  As discussed later, in section~\ref{gray.tensor.sec}, the cartesian
  product of stratified sets plays the role of the {\em Gray tensor
    product\/} in the theory of weak complicial sets. Consequently, it
  is natural to define a {\em complicial Gray-category\/} to be a
  category enriched over the cartesian category of weak complicial
  sets. In the companion paper~\cite{Verity:2005:GrayNerves} we
  generalise the {\em homotopy coherent nerve\/} construction of
  Cordier and Porter~\cite{Cordier:1986:HtyCoh} to provide a nerve
  functor which faithfully represents such complicial Gray-categories
  as weak complicial sets. Later in this work we show that the
  category of weak complicial sets is itself a complicial
  Gray-category to which we may apply this nerve construction and
  thereby represent the universe of all (small) weak complicial sets
  canonically as a (large) weak complicial set.
\end{eg}

\begin{notation}[fibrations and cofibrations]
  If $I$ is a set of morphisms in some category $\lcat{C}$ then we
  adopt the following standard notations:
  \begin{itemize}
  \item $\cell(I)$ denotes the {\em cellular completion\/} of $I$,
    that is the closure of the class of pushouts of elements of $I$
    under transfinite composition, whose elements are called {\em
      relative $I$-cell complexes},
  \item $\cof(I)$  denotes the closure of  $\cell(I)$ under retraction,
    whose elements are called {\em $I$-cofibrations}, and
  \item $\fib(I)$ denotes the class of maps which have the RLP with
    respect to $I$, whose elements are called {\em $I$-fibrations}. We
    say that an object $A$ is {\em $I$-fibrant\/} if the unique morphism
    $\arrow:A->1.$ to the terminal object is an $I$-fibration.
  \end{itemize}
  These definitions ensure that each $I$-fibration $\arrow p:A->B.$
  has the RLP with respect to any $I$-cofibration
  $\arrow i:U->V.$.

  We will assume from hereon that the reader is familiar the basic
  properties of classes of fibrations and cofibrations in a form that
  usually accompanies modern presentations of categorical homotopy
  theory. If this is not the case then any one of the commonly cited
  introductions to the basic theory of Quillen model categories should
  provide the suitable background. Certainly a familiarity with Dwyer
  and Spalinski's excellent survey article~\cite{Dwyer:1995:ModelCats}
  would suffice for our purposes here.
\end{notation}

\begin{defn}[anodyne extensions and complicial fibrations]
  We say that a stratified inclusion $\inc e:U->V.\in\Strat$ is an
  {\em (inner) anodyne extension\/} if it is in the cellular
  completion of the set of elementary (inner) anodyne extensions.
  Correspondingly, we say that a stratified map $\arrow p:E->A.$ is a
  {\em (inner) complicial fibration\/} if it is a fibration with
  respect to the set of elementary (inner) anodyne extensions.
  
  We also sometimes say that $e$ is an {\em right\/} (resp.\ {\em
    left}) anodyne extension if it is in the cellular completion of
  the union of the sets of inner and right outer (resp.\ left outer)
  anodyne extensions. Members of the corresponding class of
  fibrations are known as {\em right (resp.\ left)
    complicial fibrations}.

  Of course, we may rephrase definition~\ref{wcs.defn} in these terms
  by saying that $A$ is a weak (inner, left, right) complicial set
  iff the unique map $\arrow p:A->1.$ into the terminal stratified set
  is an (inner, left, right) complicial fibration.
\end{defn}

\begin{defn}[thinness extensions]\label{thinness.ext}
  We say that a stratified map $\overinc e:U->V.$ is a {\em thinness
    extension\/} if it is both an anodyne extension and an entire
  inclusion. By definition all elementary thinness extensions and any
  (transfinite) composite of pushouts of such things are also thinness
  extensions.  

  In general it is clearly that solutions of lifting problems whose
  domains are entire maps are unique. Furthermore, this uniqueness
  property immediately implies that if a stratified map $A$ has the
  RLP with respect to some entire map then any stratified map $\arrow
  p:A->B.$ also has that property. Consequently, it follows that any
  stratified map whose domain is a weak complicial set has the RLP
  with respect to any thinness extension.
\end{defn}

\begin{recall}[glueing squares]\label{glueing}
  A {\em glueing square\/} is a commutative square in some category
  which is both a pushout and a pullback. When constructing anodyne
  extensions we will often construct the pushouts we need as glueing
  squares, using the simple observation that if $\inc i:U->X.$ is a
  stratified inclusion and $V\subseteq_s X$ is a stratified subset of
  its codomain then the first of the following squares
  \begin{displaymath}
    \xymatrix@=2em{
      {f^{-1}(V)}\pbexcursion
      \ar@{u(->}[r]^-{i}
      \ar@{u(->}[d]_{\subseteq_s} &
      {V}\ar@{u(->}[d]^{\subseteq_s} \\
      {X}\ar@{u(->}[r]_-{i} & {f(X)\cup V}\poexcursion}
    \mkern40mu
    \xymatrix@=2em{
      {U\cap V}\pbexcursion\ar@{u(->}[r]^-{\subseteq_s}
      \ar@{u(->}[d]_{\subseteq_s} &
      {V}\ar@{u(->}[d]^{\subseteq_s} \\
      {U}\ar@{u(->}[r]_-{\subseteq_s} & {U\cup V}\poexcursion}
  \end{displaymath}
  of inclusions is a glueing square in $\Strat$. When the inclusion
  $i$ is actually a subset inclusion $\overinc\subseteq_s:U->X.$, this
  may be re-drawn to give the glueing square to its right.

  For instance, to prove that the regular inclusion $\overinc
  \subseteq_r:\Lambda^k[n]'->\Delta^k[n]''.$ is an anodyne
  extension we start with the diagram:
  \begin{displaymath}
    \xymatrix@R=2em{
      {\Lambda^k[n]}\ar@{^(->}[r]^{\subseteq_r}
      \ar@{^(->}[d]_{\subseteq_e} &
      {\Delta^k[n]}\ar@{^(->}[d]^{\subseteq_e} & \\
      {\Lambda^k[n]'}\ar@{^(->}[r]_{\subseteq_r} &
      {\Delta^k[n]'}\poexcursion\ar@{^(->}[r]_{\subseteq_e} &
      {\Delta^k[n]''}}
  \end{displaymath}
  Applying the observation of the last paragraph, we show that the
  square here is a glueing square in $\Strat$ since
  $\Lambda^k[n]=\Lambda^k[n]'\cap\Delta^k[n]$ and
  $\Delta^k[n]'=\Lambda^k[n]'\cup\Delta^k[n]$. Its upper horizontal is
  a complicial horn extension so it follows that its lower
  horizontal is an anodyne extension, which we compose with
  the elementary thinness extension to its right to obtain the
  desired presentation of $\overinc \subseteq_r:\Lambda^k[n]'
  ->\Delta^k[n]''.$ as an anodyne extension. We call this inclusion an
  {\em (inner) thin horn extension}.
\end{recall}

\begin{lemma}[superstructures of weak complicial sets]\label{super.wcs}
  For each $n\in\mathbb{N}$ the $n$-trivialisation functor $\triv_n$
  of notation~\ref{sup.defn} preserves (inner) anodyne extensions. It
  follows that its right adjoint, the superstructure functor $\sst_n$,
  preserves (inner) complicial fibrations and (inner) weak complicial
  sets.
\end{lemma}




\begin{proof} (essentially lemma~150 and lemma~171
  of~\cite{Verity:2006:Complicial}) Since $\triv_n$ is a left adjoint
  it preserves all colimits and so it is enough to check that it maps
  each {\em elementary\/} (inner) anodyne extension to an (inner)
  anodyne extension. Considering cases:

  \vspace{1ex}\noindent{\boldmath $n\geq m-1$} Observe that each of
  the stratified sets $\Delta^k[m]$, $\Lambda^k[m]$, $\Delta^k[m]'$
  and $\Delta^k[m]''$ is $(m-1)$-trivial, so if $n\geq m-1$ then they
  are also $n$-trivial. It follows that the endo-functor $\triv_n$
  maps each of these sets, and thus each of the elementary anodyne
  extensions $\overinc\subseteq_r:\Lambda^k[m]->\Delta^k[m].$ and
  $\overinc\subseteq_e:\Delta^k[m]'->\Delta^k[m]''.$, to itself.
  
  \vspace{1ex}\noindent{\boldmath $n<m-1$} Then we know that
  $\Delta^k[m]$, $\Delta^k[m]'$ and $\Delta^k[m]''$ only differ in as
  much as they have different sets of thin simplices at dimension
  $m-1$ and consequently, since $n<m-1$, we know that
  $\triv_n(\Delta^k[m])= \triv_n(\Delta^k[m]')=
  \triv_n(\Delta^k[m]'')$.  It follows that the functor $\triv_n$ maps
  the elementary thinness extension $\overinc\subseteq_e:
  \Delta^k[m]'->\Delta^k[m]''.$ to the identity on
  $\triv_n(\Delta^k[m])$.

  Finally, observe that when $n<m-1$ we know that $\Delta^k[m]''$ is
  an entire subset of $\triv_n(\Delta^k[m])$ and that its union with
  the regular subset
  $\triv_n(\Lambda^k[m])\subseteq_r\triv_n(\Delta^k[m])$ is equal to
  $\triv_n(\Delta^k[m])$ itself. Furthermore, the intersection
  $\Delta^k[m]''\cap\triv_n(\Lambda^k[m])$ is equal to the regular
  subset $\Lambda^k[m]'\subseteq_r\Delta^k[m]'')$ and it follows that
  \begin{displaymath}
    \xymatrix@=2em{
      {\Lambda^k[m]'}\ar@{u(->}[r]^{\subseteq_r}
      \ar@{u(->}[d]_{\subseteq_s} &
      {\Delta^k[m]''}\ar@{u(->}[d]^{\subseteq_s} \\
      {\triv_n(\Lambda^k[m])}\ar@{u(->}[r]_{\subseteq_r} &
      {\triv_n(\Delta^k[m])}\poexcursion}
  \end{displaymath}
  is a glueing square $\Strat$. We demonstrated that the upper
  horizontal map in this square is an anodyne extension in
  recollection~\ref{glueing}, so it follows that its pushout the lower
  horizontal is also an anodyne extension as required.

  The second sentence of the statement follows directly from the first
  under the adjunction $\triv_n\adjoint\sst_n$.
\end{proof}

\begin{obs}[alternating duals of weak complicial sets]\label{alt.dual}
  The canonical idempotent endo-functor
  $\arrow(-)\dual:\aDelta->\aDelta.$ which acts as the identity on
  objects and maps an operator $\arrow\alpha:[n]->[m].$ to the
  operator defined by $\alpha\dual(i)=m-\alpha(n-i)$ may be extended
  to a idempotent endo-functor on the category of stratified sets
  $\Strat$ called the {\em alternating dual}. This carries a
  stratified set $X$ to $X\dual$ which has the same sets of simplices
  and thin simplices as $X$ but has a dual action $\ast$ under which a
  simplicial operator $\alpha$ acts on a simplex $x$ according to the
  formula $x\ast\alpha=x\cdot\alpha\dual$.

  The action of $(-)\dual$ on operators provides us with a canonical
  isomorphism between the standard $n$-simplex $\Delta[n]$ and its
  dual $\Delta[n]\dual$ and it is clear that this underlies an
  isomorphism between the complicial simplex $\Delta^{n-k}[n]$ and the
  dual $\Delta^k[n]\dual$. Consequently, we see that on taking duals
  of the elementary anodyne extensions
  $\overinc\subseteq_r:\Lambda^k[n]->\Delta^k[n].$ and
  $\overinc\subseteq_e:\Delta^k[n]'->\Delta^k[n]''.$ we obtain
  inclusions which are isomorphic to
  $\overinc\subseteq_r:\Lambda^{n-k}[n]->\Delta^{n-k}[n].$ and
  $\overinc\subseteq_e:\Delta^{n-k}[n]'->\Delta^{n-k}[n]''.$
  respectively.

  As an idempotent functor the alternating dual is its own
  (left and right) adjoint, so in particular it preserves the colimits
  of $\Strat$ and it follows that we may extend the result of the last
  paragraph to demonstrate the preservation of all (inner) anodyne
  extensions. Furthermore, applying this adjointness property directly
  to the right lifting properties that define weak (inner)complicial
  sets and (inner) complicial fibrations, we see that these are also
  preserved by taking alternating duals.
\end{obs}

\subsection{Generalised Horns}

Before moving on we prove a simple technical lemma, of some use later
on, which shows that the inclusions associated with certain kinds of
generalised horns are (inner) anodyne extensions.  It should be noted
that this is not the most general result of this kind possible and more
powerful results of this nature may be found in
Verity~\cite{Verity:2006:Complicial} or Ehlers and
Porter~\cite{Porter:2005:SpecialPasting}. However, our simpler result
below is exactly what we will require in the sequel.

\begin{defn}\label{gen.horn}
  Suppose that $\vec{k}=\{k_1,k_2,...,k_t\}\subset[n]$ is a non-empty
  family of integers with $k_i+1<k_{i+1}$ for each $i=1,2,...,t-1$,
  then we say that an entire superset $N$ of the standard simplex
  $\Delta[n]$ is a {\em $\vec{k}$-complicial $n$-simplex\/} if it
  satisfies the following conditions for each $k_i\in\vec{k}$ and each
  $k_i$-admissible simplex $\alpha\in\Delta[n]$:
  \begin{enumerate}[label=(\alph*)]
  \item\label{gen.horn.a} $\alpha$ is thin in $N$, and
  \item\label{gen.horn.b} if $l\in[r]$ is the (unique)
    integer such that $\alpha(l)=k_i$ and $\alpha\circ\face_l$ is thin
    in $N$ then so is $\alpha\circ\face_j$ for each
    $j\in\{l-1,l+1\}\cap[r]$.
  \end{enumerate}
  Notice that if $N$ and $N'$ are two $\vec{k}$-complicial
  $n$-simplices then their intersection $N\cap N'$ is also a
  $\vec{k}$-complicial $n$-simplex. It follows that there is a
  minimal stratification which makes $\Delta[n]$ into such a
  $\vec{k}$-complicial $n$-simplex which we call $\Delta^{\vec{k}}[n]$.

  The {\em $(n-1)$-dimensional $\vec{k}$-complicial horn\/}
  $\Lambda^{\vec{k}}N$ is simply the regular subset of $N$ of those
  simplices $\arrow\alpha:[r]->[n].$ for which there is some $i\in[n]$
  which is neither in the image of $\alpha$ nor in the set $\vec{k}$
  (that is for which $[n]\neq\im(\alpha)\cup\vec{k}$).  In other
  words, this is the regular subset of $N$ generated by the set of
  $(n-1)$-simplices $\{\face_i\in\Delta[n]\mid i\in [n]\setminus
  \vec{k}\}$.  We say that $\Lambda^{\vec{k}}N$ is an {\em inner
    generalised horn\/} if $0<k_1$ and $k_t<n$.
\end{defn}

\begin{lemma}[generalised horn lemma]\label{gen.horn.lem}
  If $N$ is a (inner) $\vec{k}$-complicial $n$-simplex then the
  associated horn inclusion
  $\overinc\subseteq_r:\Lambda^{\vec{k}}N->N.$ is an (inner)
  anodyne extension.
\end{lemma} 

\begin{proof}
  Our proof is by induction on the length of $\vec{k}$. For the base
  case, if $\vec{k}=\{k\}$ then our generalised horns are no more than
  ordinary $k$-complicial horns with extra thin simplices. To be
  precise, there are two possibilities for the $(n-1)$-simplex
  $\face_k\in\Delta[n]$:
 
  \begin{enumerate}[label=\textbf{case (\roman*)},
                    fullwidth, leftmargin=1em, itemsep=1ex]
  \item It is {\em not\/} thin in $N$, in which case
    $\Delta^k[n]\subseteq_e N$ (by condition~(\ref{gen.horn.a}) of
    definition~\ref{gen.horn}), $\Delta^k[n]\cap
    \Lambda^{\vec{k}}N=\Lambda^k[n]$ and $\Delta^k[n]\cup
    \Lambda^{\vec{k}}N=N$ (since $\face_k$ is not thin in $N$) so we
    get a glueing square which displays
    $\overinc\subseteq_r:\Lambda^{\vec{k}}N->N.$ as a pushout of the
    complicial horn $\overinc\subseteq_r: \Lambda^k[n]->\Delta^k[n].$.
  \item It {\em is\/} thin in $N$, in which case
    condition~(\ref{gen.horn.b}) of definition~\ref{gen.horn} applied
    to the $n$-simplex $\arrow\id_{[n]}:[n]->[n].$ ensures that we
    actually have $\Delta^k[n]''\subseteq_e N$ and consequently that
    $\Delta^k[n]''\cap \Lambda^{\vec{k}}N=\Lambda^k[n]'$ and
    $\Delta^k[n]''\cup \Lambda^{\vec{k}}N=N$ so we get a glueing
    square which displays $\overinc\subseteq_r:\Lambda^{\vec{k}}N->N.$
    as a pushout of the (inner) thin horn extension
    $\overinc\subseteq_r: \Lambda^k[n]'->\Delta^k[n]''.$ of
    observation~\ref{glueing}.
  \end{enumerate}
  In either case it follows that the inclusion of the statement is
  an (inner) anodyne extension (as a pushout of such).

  To establish the inductive case, suppose that the result holds for
  the vector $\vec{k}=\{k_1,k_2,...,k_t\}$ and consider the extended
  vector $\vec{k}'= \vec{k}\cup\{k\}$ where $k_t+1<k\in[n]$. Suppose
  also that $N$ satisfies the conditions of definition~\ref{gen.horn}
  with respect to $\vec{k}'$. The $(n-1)$-simplex $\face_k\in N$
  corresponds to a stratified inclusion
  $\inc\yoneda{\face_k}:\Delta[n-1]->N.$ (by Yoneda's lemma) which we
  nay factor, as in definition~\ref{reg.ent.def}, to obtain an entire
  superset $M$ of $\Delta[n-1]$ and a regular inclusion
  $\overinc:M->N.$. Explicitly, the simplex $\alpha\in\Delta[n-1]$ is
  thin in $M$ iff $\face_k\circ\alpha$ is thin in $N$, from which
  description it is a matter of routine verification, using the fact
  that $k_t+1<k$, to check that $M$ also satisfies the conditions of
  definition~\ref{gen.horn} with respect to $\vec{k}$. By
  construction, the union of the image of $\overinc:M->N.$ and the
  horn $\Lambda^{\vec{k}'}N \subseteq_r N$ is the more complete
  generalised horn $\Lambda^{\vec{k}}N\subseteq_r N$ and, furthermore,
  the subset $\Lambda^{\vec{k}}M\subseteq_r M$ is easily seen to be
  the inverse image of the regular subset
  $\Lambda^{\vec{k}'}N\subseteq_r N$ along that inclusion.  So we
  obtain a commutative diagram
  \begin{equation*}
    \xymatrix@R=2em@C=3em{
      {\Lambda^{\vec{k}}M}\ar@{u(->}[r]^-{\subseteq_r}
      \ar@{u(->}[d] & {M}\ar@{u(->}[d] & \\
      {\Lambda^{\vec{k}'}N}\ar@{u(->}[r]_-{\subseteq_r} &
      {\Lambda^{\vec{k}}N}\poexcursion
      \ar@{u(->}[r]_-{\subseteq_r} & {N}}
  \end{equation*}
  in which the left hand square is a glueing square and the upper
  horizontal and right hand lower horizontal maps are both (inner)
  anodyne extensions by the induction hypothesis. It follows that the
  lower left horizontal is also an (inner) anodyne extension, since it
  is a pushout of such, and thus that its composite
  $\overinc\subseteq_r:\Lambda^{\vec{k}'}N->N.$ with the inclusion to
  its right is also an (inner) anodyne extension as required.
\end{proof}

\begin{cor}\label{gen.horn.cor}
  If $N$ is an (inner) $\vec{k}$-complicial $n$-simplex then the
  entire inclusion
  $\overinc\subseteq_e:\Lambda^{\vec{k}}N\cup\Delta^{\vec{k}}[n]->N.$
  is an (inner) anodyne extension.
\end{cor}

\begin{proof}
  A routine reprise of the method used in the proof of the last lemma,
  replacing pushouts of horn extensions by pushouts of related
  thinness extensions wherever necessary. We leave the details to the
  reader.
\end{proof}



\section{Joins of Stratified Sets}
\label{joins.sec}

Here we generalise the ever useful simplicial join operation (see for
instance~\cite{Porter:1998:Joins}
or~\cite{Joyal:2002:QuasiCategories}) to stratified sets and prove
that it gives rise to d\'{e}calage constructions which are well
behaved with respect to weak compliciality.

\subsection{Augmented Simplicial Sets and the Join Construction}

\begin{recall}[ordinal sum]
  The ordinal sum functor
  $\arrow\join:\aDelta\times\aDelta->\aDelta.$ is defined on objects
  by letting $[n]\join[m]=[n+m+1]$ and defining the sum of two
  operators $\arrow \alpha:[n]->[n'].$ and $\arrow\beta:[m]->[m'].$ by
  \begin{equation*}
    \alpha\join\beta(k) =
      \begin{cases}
        \alpha(k) & \text{if $k\leq n$,} \\
        \beta(k-n-1)+n'+1 \mkern10mu & \text{otherwise}
      \end{cases}
  \end{equation*}
  This functor makes $\aDelta$ into a strict monoidal category, whose
  unit is the empty ordinal $[-1]$.
\end{recall}

Simplicial joins are constructed by extending $\join$ to the category
of {\em augmented\/} simplicial sets
$\aSimp=[\aDelta\op,\Set]$. Correspondingly, stratified joins are
defined most naturally on augmented stratified sets.

\begin{defn}[augmented stratified sets]
  An {\em augmented stratified set\/} $X$ consists of an augmented
  simplicial set equipped with a stratification $tX\subseteq X$
  satisfying the single condition that no $(-1)$-dimensional simplices
  should be members of the subset $tX$. In other words, this is no
  more than a stratified set $X$ together with a chosen {\em
    augmentation}, that being a delegated subset of thin $0$-simplices
  $tX_0\subseteq X_0$, a set of $(-1)$-simplices $X_{-1}$
  and a function $\arrow d_0:X_0->X_{-1}.$ satisfying the simplicial
  identity $\arrow d_0\circ d_0 = d_0\circ d_1:X_1->X_{-1}.$.

  All of the basic definitions and results of the theory of stratified
  sets carry over to this context, in particular an (augmented)
  stratified map between such structures is simply an (augmented)
  simplicial map which preserves thinness.  We let $\aStrat$ denote
  the category of augmented stratified sets and their stratified maps.
\end{defn}

\begin{obs}
  The canonical functor $\arrow:\aStrat->\Strat.$ which forgets
  augmentations has both a left and a right adjoint, providing us with
  two ``opposed'' augmentations of any stratified set $X$:
  \begin{itemize}
  \item The {\em canonical augmentation\/} (left adjoint) with
    $tX_0=\emptyset$ and $X_{-1}\defeq\pi_0(X)$ constructed by
    coequalising the pair of maps $\arrow d_0, d_1:X_1->X_0.$.
  \item The {\em trivial augmentation\/} (right adjoint) with
    with $tX_0=X_0$ and $X_{-1}=\{*\}$ the one point set.
  \end{itemize}
  We make no particular choice of default augmentation, preferring
  instead to specify an appropriate augmentation on a case by case
  basis. 
\end{obs}

\begin{defn}[joins of augmented stratified sets]
  Day's convolution construction \cite{Day:1970:ClosedFunc} allows us
  to extend the monoidal structure $\join$ on $\aDelta$ to a monoidal
  closed structure on $\aSimp$. The tensor product of this structure,
  also denoted by $\join$, is called the {\em simplicial join\/} and
  the corresponding closures $\dec_l(X,Z)$ and $\dec_r(Y,Z)$ are
  called the {\em left and right d\'ecalage\/} constructions
  respectively. In line with traditional usage we will sometimes use
  the notations $\dec_l(Z)$ and $\dec_r(Z)$ to denote the d\'ecalages
  $\dec_l(\Delta[0],Z)$ and $\dec_r(\Delta[0],Z)$ respectively.

  If $X$ and $Y$ are two (augmented) simplicial sets then Day's
  convolution formula tells us that their join is given by the coend
  formula:
  \begin{equation*}
    (X\join Y)_r = \int^{[n],[m]\in\aDelta} X_n\times Y_m \times
    \aDelta([r],[n]\join [m])
  \end{equation*}
  A routine calculation demonstrates that an $r$-simplex of this join
  corresponds to a pair $\langle x,y \rangle$ with $x\in X_n$ and
  $y\in Y_m$ for some pair of integers $n,m\geq -1$ with $[n]\join
  [m]=[r]$. Under this representation if $\arrow\beta:[n']->[n].$ and
  $\arrow\gamma:[m']->[m].$ are simplicial operators then we have
  $\langle x,y\rangle\cdot(\beta\join\gamma)=\langle
  x\cdot\beta,y\cdot\gamma\rangle$ and this identity completely
  determines the action of $\aDelta$ on $X\join Y$ since any operator
  $\arrow \alpha:[r']->[r].$ with $[r]=[n]\join[m]$ may be decomposed
  as $\alpha=\beta\join\gamma$ for a unique pair of such operators.

  We now extend this to (augmented) stratified sets $X,Y\in\aStrat$ by
  applying $\join$ to their underlying (augmented) simplicial sets
  and letting $\langle x,y\rangle$ be thin in $X\join Y$ if and only
  if $x$ is thin in $X$ or $y$ is thin in $Y$. It is clear that this
  provides a monoidal structure on $\aStrat$ and it is a routine
  matter to check that each of the endo-functors $X\join{-}$ and
  ${-}\join Y$ preserve the colimits of $\aStrat$ simply by
  observing that by definition they do so on underlying (augmented)
  simplicial sets and checking that the resulting comparison
  isomorphisms reflect thinness appropriately. It follows therefore
  that these functors have right adjoints, which we again denote by
  $\dec_l(X,{*})$ and $\dec_r(Y,{*})$ respectively. 
\end{defn}

\begin{obs}[joins, d{\'e}calage and augmentation]
  We must, of course, augment all stratified sets before applying the
  join or d{\'e}calage constructions to them. We shall adopt different
  implicit augmentation conventions for each of these operations:
  \begin{enumerate}[fullwidth, leftmargin=1em, itemsep=1ex]
  \item[{\bf Joins}] we apply canonical augmentation in either
    variable. This ensures that joins preserve colimits of stratified
    sets independently in each variable and that the join of two
    stratified sets is again canonically augmented.
  \item[{\bf D{\'e}calages}] we apply canonical augmentation in the
    first (contravariant) variable and trivial augmentation in the
    second (covariant) variable. This ensures that d{\'e}calages carry
    colimits of stratified sets in the contravariant variable and
    limits of stratified sets in the covariant variable to limits in
    $\aStrat$.
  \end{enumerate}
\end{obs}

\begin{obs}[augmented standard simplices and their boundaries]
  \label{aug.simplices}
  In the context of augmented simplicial sets the notation $\Delta[n]$
  will stand for the representable on the object $[n]$ {\em as an
    object of\/} $\aDelta$ and $\boundary\Delta[n]$ will stand for its
  subset of non-surjective operators. These are all trivial
  augmentations of the corresponding un-augmented structures, and in
  most cases they coincides with the corresponding canonical
  augmentation. Indeed, the only exceptions to this rule are the sets
  $\Delta[-1]$, $\boundary\Delta[0]$ and $\boundary\Delta[1]$.
\end{obs}

\begin{obs}[joins and alternating duals]
  Joins of (augmented) stratified sets are well behaved with respect
  to the alternating dual of observation~\ref{alt.dual}.  To be
  precise, observe that if $\alpha$ and $\beta$ are a pair of
  simplicial operators then we have
  $(\alpha\join\beta)\dual=\beta\dual\join\alpha\dual$ from which it
  follows immediately that the ``swap'' function, which carries a pair
  $\langle x,y\rangle$ to the reversed pair $\langle y,x\rangle$,
  provides us with a stratified isomorphism between $(X\join Y)\dual$
  and $Y\dual\join X\dual$ which is natural in $X$ and $Y$. By
  adjointness, these isomorphisms provide us with canonical natural
  isomorphisms $\dec_l(X,Z)\dual\cong\dec_r(X\dual,Z\dual)$ and 
  $\dec_r(Y,Z)\dual\cong\dec_l(Y\dual,Z\dual)$.

  It follows that in the sequel it will be enough to consider left
  joins $X\join{-}$ and the corresponding left d\'ecalage closures
  $\dec_l(X,{*})$, since the properties of right joins and d\'ecalage
  follow on applying alternating duals and the isomorphisms of the
  last paragraph.
\end{obs}

\subsection{D{\'e}calage and Weak Compliciality}

\begin{obs}[corner joins]\label{corner.join}
  Applying the construction of recollection~\ref{corner.tensor} to the
  join of augmented stratified sets we obtain the {\em corner join and
    corner d\'ecalage\/} constructions which we denote by $\cjoin$,
  $\cdec_l$ and $\cdec_l$ respectively. Generally we
  are interested in taking the corner join of two (augmented)
  stratified subset inclusions $\overinc\subseteq_s:U->V.$ and
  $\overinc\subseteq_s:X->Y.$. In which case we know that $U\join Y$
  and $V\join X$ are stratified subsets of $V\join Y$ with $(U\join
  Y)\cap(V\join X)=U\join X$ and it follows, by
  observation~\ref{glueing}, that we have a glueing square
  \begin{equation*}
    \xymatrix@R=2em@C=3em{
      {U\join X}\ar@{u(->}[r]^{\subseteq_s}
      \ar@{u(->}[d]_{\subseteq_s} &
      {V\join X}\ar@{u(->}[d]^{\subseteq_s} & \\
      {U\join Y}\ar@{u(->}[r]_-{\subseteq_s} &
      {(U\join Y)\cup(V\join X)}\poexcursion
      \ar@{u(->}[r]_-{\subseteq_s} & {V\join Y}
    }
  \end{equation*}
  which demonstrates that the inclusion to the right of its lower
  right hand corner is (isomorphic to) the corner join of our
  inclusions. One useful observation that follows from this is that if
  the inclusion of $U$ into $V$ is actually entire then $V\join Y$ and
  $U\join Y$ have the same underlying (augmented) simplicial sets from
  which it follows that the corner join depicted above is also an
  entire inclusion.
\end{obs}

\begin{obs}[joins and anodyne extensions]\label{join.anodyne}
  By construction the join operation on $\aStrat$ extends ordinal sum
  on $\aDelta$ so it follows that the ordinal sum of operators
  provides a canonical isomorphism $\Delta[n]\join
  \Delta[m]\cong\Delta[n+m+1]$ in $\aStrat$ which maps each pair
  $\langle \alpha,\beta \rangle\in\Delta[n] \join\Delta[m]$ to
  $\alpha\join\beta\in\Delta[n+m+1]$. Furthermore if $0\leq k <n$ then
  a simplex $\alpha$ is thin in the complicial simplex $\Delta^k[n]$
  if and only if $\alpha\join\beta$ is thin in $\Delta^k[n+m+1]$ for
  each simplex $\beta$ in $\Delta[m]$, so it follows that the
  isomorphism of the last sentence extends to a stratified isomorphism
  $\Delta^k[n]\join\Delta[m]\cong \Delta^k[n+m+1]$.

  Now observe that if $\alpha\in\Delta[n]$ and $\beta\in\Delta[m]$
  then $\alpha\join\beta\in\Delta^k[n+m+1]$ is in the complicial horn
  $\Lambda^k[n+m+1]$ if and only if $\alpha$ is in $\Lambda^k[n]$ or
  $\beta$ is in $\boundary\Delta[m]$. So the isomorphism of the last
  paragraph restricts to provide an isomorphism between the inclusion
  $\overinc\subseteq_s:(\Lambda^k[n]\join\Delta[m])
  \cup(\Delta^k[n]\join\boundary\Delta[m])->\Delta^k[n]\join
  \Delta[m].$, which is simply the corner join of the inclusions
  $\overinc\subseteq_r: \Lambda^k[n]->\Delta^k[n].$ and
  $\overinc\subseteq_r:\boundary \Delta[m]->\Delta[m].$ by
  observation~\ref{corner.join}, and the complicial horn
  $\overinc\subseteq_r:\Lambda^k[n+m+1]-> \Delta^k[n+m+1].$.

  This argument does not quite apply when $k=n$ for then it is not
  the case that $\Delta^n[n]\join\Delta[m]$ is isomorphic to
  $\Delta^n[n+m+1]$ since they then have slightly different
  stratifications. However, it may be adapted to show that in this
  case the corner join of the last paragraph can be presented as the
  lower horizontal map in a pushout
  \begin{equation}\label{right.inner.po}
    \xymatrix@R=1.5em@C=3em{
      {\Lambda^n[n+m+1]}\ar@{u(->}[d]
      \ar@{u(->}[r]^{\subseteq_r} &
      {\Delta^n[n+m+1]}\ar@{u(->}[d] \\
      {(\Lambda^n[n]\join\Delta[m])\cup(\Delta^n[n]\join\boundary
        \Delta[m])}\ar@{u(->}[r]^<>(0.5){\subseteq_r} &
      {\Delta^n[n]\join\Delta[m]}\poexcursion
    }
  \end{equation}
  and is thus an anodyne extension.

  We may apply a similar argument to the corner join of 
  $\overinc\subseteq_r: \Lambda^k[n]->\Delta^k[n].$ and
  $\overinc\subseteq_e:\Delta[n]->\Delta[n]_t.$, which is an entire
  subset inclusion since the second of these inclusions is itself
  entire (observation~\ref{corner.join}). Now, if a non-degenerate
  simplex $\langle \alpha,\beta\rangle$ is thin in the codomain
  $\Delta^k[n]\join \Delta[m]_t$ of this corner join then by
  definition either $\alpha$ is thin in $\Delta^k[n]$, in which case
  $\langle\alpha,\beta\rangle$ is thin in $\Delta^k[n]\join\Delta[m]$,
  or $\beta$ is thin in $\Delta[m]_t$, in which case it is also thin
  in $\Lambda^k[n]\join\Delta[m]_t$ unless $\alpha$ fails to be a
  simplex of the horn $\Lambda^k[n]$ altogether. In that latter case
  either $\alpha=\id_{[n]}$ or $\alpha=\face_k$, and again the first
  of these would make $\langle\alpha,\beta\rangle$ thin in
  $\Delta^k[n]\join\Delta[m]$. Summarising this argument, we see that
  if the simplex $\langle\alpha,\beta\rangle$ is thin in the codomain
  of the corner join under consideration and non-thin in its domain
  then it can only be the simplex
  $\langle\face_k,\id_{[m]}\rangle\in\Delta[n]\join\Delta[m]$, which
  corresponds to the simplex $\face_k\in\Delta[n+m+1]$ under the
  canonical ordinal sum isomorphism
  $\Delta[n]\join\Delta[m]\cong\Delta[n+m+1]$. It follows,
  immediately, that we may present this corner join as the lower
  horizontal of a pushout square
  \begin{displaymath}
    \xymatrix@R=1.5em@C=3em{
      {\Delta^k[n+m+1]'}\ar@{u(->}[d]
      \ar@{u(->}[r]^{\subseteq_r} &
      {\Delta^k[n+m+1]''}\ar@{u(->}[d] \\
      {(\Lambda^k[n]\join\Delta[m]_t)\cup(\Delta^k[n]\join
        \Delta[m])}\ar@{u(->}[r]^<>(0.5){\subseteq_r} &
      {\Delta^k[n]\join\Delta[m]_t}\poexcursion
    }
  \end{displaymath}
  and we may infer that it is thus an inner anodyne extension.
  Finally, an analogous analysis of the thinness extension
  $\overinc\subseteq_e:\Delta^k[n]'-> \Delta^k[n]''.$ shows that its
  corner joins with the boundary and thin simplex inclusions can again
  be obtained as pushouts of the thinness extension
  $\overinc\subseteq_e: \Delta^k[n+m+1]'->\Delta^k[n+m+1]''.$, which
  again demonstrates that they are both anodyne extensions.
\end{obs}

Summarising these observations we get the following lemma and its
obvious dual involving left handed d\'ecalage:

\begin{lemma}[weak complicial sets and d\'{e}calage]\label{wcs.decalage}
  If $\inc e:U->V.$ is an anodyne extension and $\inc i:X->Y.$ is any
  inclusion (monomorphism) of augmented stratified sets then their
  corner join $e\cjoin i$ is an anodyne extension. It follows that
  if $\arrow p:E->B.$ is a complicial fibration then so is the right
  corner d\'ecalage $\cdec_r(i,p)$.

  Thus, if $Y$ is any augmented stratified set then the
  endo-functor ${-}\join Y$ preserves all anodyne extensions and if
  $A$ is a weak complicial set then so is $\dec_r(Y,A)$.
\end{lemma}

\begin{proof}
  To prove the first part, observe that the class of all inclusions of
  augmented stratified sets is the cellular completion of the set of boundary
  and thin simplex inclusions:
  \begin{equation}\label{inc.gen}
    \{\overinc\subseteq_r:\boundary\Delta[n]->\Delta[n].\mid
    n=-1,0,1,...\}\cup\{\overinc\subseteq_e:\Delta[n]->\Delta[n]_t.\mid 
    n=0,1,2,...\}
  \end{equation}
  The calculations of the last observation demonstrated that the
  corner join of any of the inclusions in this set with an elementary
  anodyne extension is again an anodyne extension, so we may apply
  lemma~\ref{corner.cofibration} to extend this result to all
  inclusions and anodyne extensions as required. The second sentence
  of the statement now follows by applying
  observation~\ref{corner.fibration}.

  Finally, observe that if $\overinc i:\emptyset->Y.$ is the unique
  inclusion from the empty augmented stratified set into $Y$ then the
  corner join $e\cjoin i$ is clearly isomorphic to $\arrow e\join
  Y:U\join Y->V\join Y.$, since joins preserve the initial object
  $\emptyset$ which implies that the domain
  $(V\join\emptyset)\vee_{U\join\emptyset} (U\join Y)$ of our
  corner join is isomorphic to $U\join Y$. Similarly the corner
  d\'ecalage $\cdec_r(i,p)$ is isomorphic to the map
  $\arrow\dec_r(Y,p): \dec_r(Y,E)->\dec_r(Y,B).$. It follows that the
  result of the last paragraph specialises to establish the final
  sentence of the statement.
\end{proof}

\begin{obs}[inner anodyne extensions and joins]\label{inner.joins}
  Notice that  observation~\ref{join.anodyne} actually demonstrates
  that if we corner join an elementary {\bf inner} anodyne extension with a
  boundary or thin simplex inclusion then the resulting map is in fact
  also an inner anodyne extension. This immediately implies that
  lemma~\ref{wcs.decalage} has a direct analogue in which anodyne
  extensions, complicial fibrations and weak complicial sets are
  replaced by their inner counterparts.

  However, consulting observation~\ref{join.anodyne} again in greater
  detail we see that a little more is true. In particular, observe
  that in the pushout of display~(\ref{right.inner.po}) the upper
  horizontal map is actually an inner horn extension whenever $m\geq
  0$ and so it follows that its lower horizontal, the corner join of
  the right outer horn $\overinc\subseteq_r:\Lambda^n[n]
  ->\Delta^n[n].$ and the boundary inclusion $\overinc\subseteq_r:
  \boundary\Delta[m]->\Delta[m].$, is also an inner anodyne
  extension. A similar comment holds for the other three cases of
  observation~\ref{join.anodyne} in which an elementary right outer
  anodyne extension is corner joined with a boundary or thin simplex
  inclusion.  It follows therefore, by applying
  lemma~\ref{corner.cofibration} again, that the corner join
  $e\cjoin i$ of a {\bf right} anodyne extension $e$ and an
  inclusion $i$ in the cellular completion of the set obtained by removing
  $\overinc\subseteq_r:\boundary\Delta[-1]->\Delta[-1].$ from the set
  in display~\eqref{inc.gen} is actually an {\bf inner} anodyne
  extension. Furthermore, a simple argument demonstrates that a map $i$
  in this latter cellular completion if and only if it is an inclusion of
  augmented stratified sets which acts isomorphically on sets of
  $(-1)$-simplices.
\end{obs}
 


\section{Equivalences in Weak Complicial Sets}
\label{equiv.sec}

In this section we provide an alternative characterisation of weak
complicial sets, which replaces outer complicial horn fillers with an
{\em equivalence\/} condition on thin 1-simplices. This theory
directly generalises the analysis of quasi-isomorphisms given by Joyal
in his paper on
quasi-categories~\cite{Joyal:2002:QuasiCategories}. Herein the
material of this section primarily serves to simplify subsequent work,
by freeing us from directly analysing certain special cases involving
outer horns.

\subsection{Equivalences in Simplicial Sets}

\begin{defn}[the generic simplicial equivalence]
  \label{fs.equiv.defn}
  Let $\mathbb{I}$ denote the chaotic category on two objects
  $\{-,+\}$, which is generally referred to as the {\em generic
    isomorphism}, and let $E\in\Simp$ denote its nerve (cf.\
  observation~\ref{cat.nerves}).  In other words, $E$ is the
  simplicial set whose $m$-simplices are, not necessarily order
  preserving, functions $\arrow e:[m]->\{-,+\}.$ upon which simplicial
  operators act by pre-composition.  We can think of an $m$-simplex of
  $E$ as a sequence $e_0e_1...e_m$ of the symbols $-$ and $+$ of
  length $n+1$ upon which a simplicial operator
  $\arrow\alpha:[n]->[m].$ acts by re-indexing
  $(e_0e_1...e_m)\cdot\alpha=e_{\alpha(0)}e_{\alpha(1)}...e_{\alpha(n)}$.
  For reasons that will become apparent, we call $E$ the {\em generic
    simplicial equivalence\/} and we say that a 1-simplex $v$ of a
  simplicial set $X$ is a {\em (simplicial) equivalence\/} if there
  exists some simplicial map $\arrow f:E->X.$ with $f(-+)=v$.
\end{defn}

In what follows, we sometimes use the symbols $p$ and $q$ to represent
elements of $\{-,+\}$ and use the notation $\neg$ to denote the
function which swaps $+$ and $-$.

\begin{obs}[decomposing $E$]\label{equiv.decomp}
  An $n$-simplex $e\in E$ is degenerate iff there is some $i\in[n-1]$
  for which $e_i=e_{i+1}$. It follows that $E$ has exactly 2
  non-degenerate $n$-simplices, these being the two alternating
  sequences of length $n+1$ starting from $-$ and $+$ respectively for
  which we reserve the notation $e^-_n=-+-+...$ and $e^+_n=+-+-...$.

  Let $E^p_n$ denote the simplicial subset of $E$ generated by the
  simplex $e^p_n$ and observe that the obvious identities
  $e^p_{n+1}\cdot\face_0=e^{\neg p}_n$ and $e^p_{n+1}\cdot\face_{n+1}=
  e^p_n$ imply that we have $E^q_n\subseteq_s E^p_{n+1}$ for each
  $n\in\mathbb{N}$ and $p,q\in\{-,+\}$ and that $E$ itself is the
  union of the increasing chain $E^p_1\subseteq_s E^p_2\subseteq_s
  \cdots \subseteq_s E^p_n\subseteq_s\cdots$.  Furthermore, they also
  imply that the only two non-degenerate simplices of $E^p_{n+1}$
  which are not in $E^p_n$ (resp.\ $E^{\neg p}_n$) are $e^p_{n+1}$
  itself and its face $e^{\neg p}_n=e^p_{n+1}\cdot\face_0$
  (resp. $e^p_n=e^p_{n+1}\cdot \face_{n+1}$). It follows that we have
  canonical pushout squares
  \begin{equation}\label{equiv.decomp.po}
    \xymatrix@=2em{
      {\Lambda^0[n+1]}\ar@{u(->}[r]^<>(0.5){\subseteq_s}\ar[d]&
      {\Delta[n+1]}\ar[d]^{\yoneda{e^p_{n+1}}} \\
      {E^p_n}\ar@{u(->}[r]_{\subseteq_s} & {E^p_{n+1}}\poexcursion }\mkern50mu
    \xymatrix@=2em{
      {\Lambda^{n+1}[n+1]}\ar@{u(->}[r]^<>(0.5){\subseteq_s}\ar[d] &
      {\Delta[n+1]}\ar[d]^{\yoneda{e^p_{n+1}}} \\
      {E^{\neg p}_n}\ar@{u(->}[r]_{\subseteq_s} & {E^p_{n+1}}\poexcursion }
  \end{equation}
  in $\Simp$, in which $\arrow\yoneda{e^p_{n+1}}:\Delta[n+1]->
  E^p_{n+1}.$ is the simplicial map that corresponds to the
  $(n+1)$-simplex $e^p_{n+1}\in E^p_{n+1}$ via Yoneda's lemma. We will
  also use the notation $E_n$ to denote the union of the subsets
  $E^-_n$ and $E^+_n$ in $E$.
\end{obs}

\begin{obs}[equivalences in weak complicial sets]
  \label{equiv.decomp.strat}
  From hereon we will adopt the (slightly nonstandard) convention that
  the simplicial sets $E^-_n$, $E^+_n$ and $E_n$ are all stratified
  with the $0$-trivialised stratification, in which a simplex is thin
  iff its dimension is greater than $0$. To recover the default minimal
  stratification on these sets we apply the underlying simplicial set
  notation $\usimp{E}^-_n$, $\usimp{E}^+_n$ and $\usimp{E}_n$
  and appeal to the default stratification rule.

  Lifting the left hand pushout of display~\eqref{equiv.decomp.po} to
  $\Strat$ we find that the inclusion
  $\overinc\subseteq_r:E^p_n->E^p_{n+1}.$ is a pushout of the left
  horn extension $\overinc
  \subseteq_r:\Lambda^0[n+1]'->\Delta^0[n+1]''.$ and is thus itself a
  left anodyne extension.  Arguing dually we see that the inclusion
  $\overinc\subseteq_r: E^{\neg p}_n-> E^p_{n+1}.$ is a right anodyne
  extension.  Taking composites of these it follows that each
  inclusion $\overinc \subseteq_r:E^p_m->E^p_n.$ is a left anodyne
  extension and that $\overinc\subseteq_r:E^q_m-> E^p_n.$ is a right
  anodyne extension if $(n-m)$ is even and $p=q$ or if $(n-m)$ is odd
  and $p=\neg q$.

  In particular, the inclusion $\overinc \subseteq_r:E^-_1->E.$ may be
  constructed as a countable composite of the inclusions $E^-_1
  \subseteq_rE^-_2\subseteq_rE^-_3\subseteq_r \cdots$ so it follows
  from the last paragraph that it is a left anodyne extension. Indeed,
  we may also construct it as a countable composite of the
  (alternating) sequence $E^-_1
  \subseteq_rE^+_2\subseteq_rE^-_3\subseteq_r E^+_4\subseteq_r\cdots$
  which demonstrates that it is also a right anodyne extension.

  Now observe that the stratified set $E^-_1$ is simply isomorphic to
  the standard thin $1$-simplex $\Delta[1]_t$, so it follows that any
  thin $1$-simplex $v$ of a weak complicial set $A$ gives rise to a
  unique stratified map $\arrow\yoneda{v}:E^-_1->A.$ with
  $\yoneda{v}(e^-_1)=v$. Since $A$ is a weak complicial set this may
  be lifted along the anodyne extension $\overinc
  \subseteq_r:E^-_1->E.$ to give a stratified map which demonstrates
  that $v$ is a equivalence in the underlying simplicial set of $A$.
\end{obs}

\begin{obs}[some symmetries of $E$]\label{iso.symm}
  In the sequel we will have use for a couple of canonical
  isomorphisms defined upon $E$:
  \begin{enumerate}[label=\textbf{symmetry (\roman*)},
                    fullwidth, leftmargin=1em, itemsep=1ex]
  \item The function $\arrow \neg:E->E.$ which applies the parity swapping
    function $\neg$ pointwise to the symbols comprising each simplex
    of $E$ and is clearly an idempotent map of simplicial
    sets. Furthermore, this restricts to provide an isomorphism between
    $E^p_n$ and $E^{\neg p}_n$ for each $n\in\mathbb{N}$ and
    $p\in\{-,+\}$.
  \item The function ``$\rev$'' which reverses the order of the symbols in
    each simplex of $E$ and is clearly the underlying function of a
    mutually inverse pair of simplicial isomorphisms
    $\arrow\rev:E->E\dual.$ and $\arrow\rev:E\dual->E.$. Furthermore
    these restrict to provide an isomorphism between $(E^p_n)\dual$ and
    $E^p_n$ if $n$ is even and $E^{\neg p}_n$ if $n$ is odd.
  \end{enumerate}
\end{obs}

\subsection{Equivalences and Inner Compliciality}

Conversely, the following sequence of observations and lemmas
demonstrate that an equivalence property on thin 1-simplices is enough
to ensure that a weak inner complicial set has outer horn
fillers. This result may be considered to be a complicial
generalisation of Joyal's analysis of {\em special horn fillers\/} in
quasi-categories~\cite{Joyal:2002:QuasiCategories}.

\begin{obs}\label{inner.equiv.obs}
  Our primary goal over the next few lemmas will be to show that a
  weak inner complicial set $A$ is actually a weak complicial set iff
  it has the RLP with respect to the inclusion
  $\overinc\subseteq_r:E^-_1 ->E^-_3.$.

  To that end we start by observing that the inclusion
  $\overinc\subseteq_r:E^-_0-> E^-_1.$ is isomorphic
  to $\inc\Delta(\vertex_0):\Delta[0]-> \Delta[1]_t.$ and arguing as
  in observation~\ref{join.anodyne} to show that its corner join with
  the inclusion $\overinc\subseteq_r:\boundary\Delta[m]-> \Delta[m].$
  is isomorphic to the left outer horn
  $\overinc\subseteq_r:\Lambda^0[m+2]->\Delta^0[m+2].$.  This leads us
  to considering the following increasing sequence of stratified
  subsets of $E^-_3\join\Delta[m]$ which starts with the
  domain of this corner join
  \begin{align*}
    X_0 &\defeq (E^-_0\join \Delta[m])\cup (E^-_1
    \join\boundary\Delta[m]) \\
    X_1 &\defeq  (E^-_0\join \Delta[m])\cup (E^-_1
    \join\boundary\Delta[m])\cup (E^-_3\join\joinidt) \\
    X_2 &\defeq (E^-_0\join \Delta[m])\cup (E^-_3
    \join\boundary\Delta[m]) \\
    X_3 &\defeq (E^-_2\join \Delta[m])\cup (E^-_3
    \join\boundary\Delta[m])
  \end{align*}
  and ends with a stratified set containing its codomain
  $E^-_1\join\Delta[m]$.  Here $\joinidt$ in the definition of $X_1$
  represents the stratified subset of $\Delta[m]$ containing only its
  unique $-1$ dimensional simplex $\arrow\iota^m:[-1]->[m].$. The
  following observations follow directly from these definitions:

  \begin{enumerate}[label=\textbf{observation (\roman*)},
                    fullwidth, leftmargin=1em, itemsep=1ex]
  \item We have $X_0\cup (E^-_3\join\joinidt)=X_1$
    and $X_0\cap (E^-_3\join\joinidt)=E^-_1\join\joinidt$
    so we obtain a glueing square which presents the inclusion
    $\overinc\subseteq_r: X_0->X_1.$ as a pushout of the inclusion
    $\overinc\subseteq_r: E^-_1\join\joinidt->E^-_3
    \join\joinidt.$. Furthermore this, in turn, is isomorphic to the
    inclusion $\overinc\subseteq_r:E^-_1-> E^-_3.$
    of the statement, since $\joinidt$ is isomorphic to $\Delta[-1]$ the
    identity for $\join$.

  \item We have the equalities $X_1\cup (E^-_3\join\boundary
    \Delta[m])= X_2$ and $X_1\cap (E^-_3\join\boundary \Delta[m]) =
    (E^-_1 \join\boundary\Delta[m])\cup (E^-_3 \join\joinidt)$ thus
    ensuring that we have a glueing square which presents the
    inclusion $\overinc\subseteq_r:X_1->X_2.$ as a pushout of the
    inclusion $\overinc\subseteq_r: (E^-_1 \join
    \boundary\Delta[m])\cup (E^-_3\join\joinidt)->E^-_3
    \join\boundary\Delta[m].$ which, in turn, is the corner join of
    $\overinc\subseteq_r:E^-_1->E^-_3.$ and
    $\overinc\subseteq_r:\joinidt->\boundary\Delta[m].$ (by
    observation~\ref{corner.join}). Now the first of these is a right
    anodyne extension, as demonstrated in
    observation~\ref{equiv.decomp.strat}, so we may apply
    observation~\ref{inner.joins} to show that their corner join is an
    inner anodyne extension.

  \item We have the equalities $X_2\cup (E^-_2\join\Delta[m]) = X_3$
    and $X_2\cap (E^-_2\join\Delta[m]) =
    (E^-_0\join\Delta[m])\cup(E^-_2\join\boundary \Delta[m])$ thus
    ensuring that we have a glueing square which presents the
    inclusion $\overinc\subseteq_r:X_2->X_3.$ as a pushout of the
    inclusion $\overinc\subseteq_r: (E^-_0 \join \Delta[m])\cup
    (E^-_2\join\boundary\Delta[m])-> E^-_2\join\Delta[m].$ which, in
    turn, is the corner join of $\overinc\subseteq_r:E^-_0->E^-_2.$ and
    $\overinc\subseteq_r:\boundary\Delta[m]->\Delta[m].$ (by
    observation~\ref{corner.join}). Now the first of these is a right
    anodyne extension, as demonstrated in
    observation~\ref{equiv.decomp.strat}, so we may apply
    observation~\ref{inner.joins} to show that their corner join is an
    inner anodyne extension.
  \end{enumerate}

  From these it follows immediately that the inclusion $\overinc
  \subseteq_r:X_0->X_3.$ is in the cellular completion of the set of
  inclusions obtained by adding $\overinc\subseteq_r:E^-_1
  ->E^-_3.$ to the set of elementary inner anodyne
  extensions. 
\end{obs}

\begin{obs}\label{inner.equiv.obs.2}
  A similar sequence of observations hold for complicial thinness
  extensions, however this time we consider the corner join of
  $\overinc\subseteq_r:E^-_0->E^-_1.$ and
  $\overinc\subseteq_e:\Delta[m]->\Delta[m]_t.$ and argue along the
  lines presented in the latter part of observation~\ref{join.anodyne}
  to show that it is isomorphic to the left outer thinness extension
  $\overinc\subseteq_e:\Delta^0[m+2]'->\Delta^0[m+2]''.$. This again
  leads us to considering an increasing sequence of stratified subsets
  of $E^-_3\join\Delta[m]_t$ which starts with the domain of our
  corner join
  \begin{align*}
    Y_0 &\defeq (E^-_0\join \Delta[m]_t)\cup (E^-_1
    \join\Delta[m]) \\
    Y_1 &\defeq  (E^-_0\join \Delta[m]_t)\cup (E^-_1
    \join\Delta[m])\cup (E^-_3\join\joinidt) \\
    Y_2 &\defeq (E^-_0\join \Delta[m]_t)\cup (E^-_3
    \join\Delta[m]) \\
    Y_3 &\defeq (E^-_2\join \Delta[m]_t)\cup (E^-_3
    \join\Delta[m])
  \end{align*}
  and ends with a stratified set containing its codomain
  $E^-_1\join\Delta[m]_t$. Arguing exactly as we did in the
  subclauses of the last observation, we see that the first inclusion
  in this sequence is a pushout of $\overinc\subseteq_r:E^-_1
  ->E^-_3.$ and that its last two are both inner anodyne
  extensions.
\end{obs}

\begin{lemma}[lifting of equivalences is enough]
  \label{inner.equiv.lem}
  Suppose that $B$ is a weak complicial set and $\arrow p:A->B.$ is an
  inner complicial fibration which has the RLP with
  respect to the inclusions $\overinc\subseteq_r:E^-_1 ->E^-_3.$ and
  $\overinc\subseteq_r:E^-_0->E^-_1.$ then $p$ is a complicial
  fibration. Consequently, $A$ is a weak complicial set iff it is a
  weak inner complicial set which has the RLP with
  respect to $\overinc\subseteq_r:E^-_1 ->E^-_3.$.
\end{lemma}

\begin{proof}
  First observe that we may apply (both of) the symmetry isomorphisms
  of observation~\ref{iso.symm} to show that the dual inclusion
  $\overinc\subseteq_r:(E^-_1)\dual->(E^-_3)\dual.$ is actually
  isomorphic to $\overinc\subseteq_r:E^-_1->E^-_3.$ itself. On the
  other hand the inclusion $\overinc\subseteq_r:E^-_0->E^-_1.$ is not
  self dual, instead its dual is isomorphic to
  $\overinc\subseteq_r:E^-_0->E^+_1.$. However, it is still the case
  that the conditions of the statement are enough to ensure that
  $\arrow p:A->B.$ also has the RLP with respect to
  this latter inclusion. 

  To see that this is the case, simply observe that $E^+_1$ is also a
  stratified subset of $E^-_3$ from which it follows that we may solve
  a lifting problem $(u,v)$ in the following diagram
  \begin{equation*}
    \xymatrix@R=0.5em@C=2em{
      {E^-_0}\ar[rr]^-{u}\ar@{u(->}[dd]_{\subseteq_r} &&
      {A}\ar[dd]^{p} \\
      & {E^-_3}\ar@{..>}[dr]^{v'}\ar@{..>}[ur]_{w} & \\
      {E^+_1}\ar[rr]_-{v}\ar@{u(->}[ru]^-{\subseteq_r}_{i}
      && {B}}
  \end{equation*}
  in two steps. First use the weak compliciality of $B$ to extend the
  stratified map $v$ along the anodyne extension
  $\overinc\subseteq_r:E^+_1->E^-_3.$ (cf.\
  observation~\ref{equiv.decomp.strat}) to obtain the dotted map
  $v'$. Then observe that the inclusion
  $\overinc\subseteq_r:E^-_0->E^-_3.$ may be decomposed as a composite
  of the inclusions identified in the statement, from which it follows
  that $p$ also has the RLP with respect this
  latter inclusion. This allows us to solve the new lifting problem
  $(u,v')$ and obtain the stratified map $w$, which we compose with
  the inclusion $i$ to finally construct the desired solution to the
  original lifting problem.

  Applying this result and the fact that the class of inner complicial
  fibrations is closed under alternating duals, we have demonstrated
  that $\arrow p:A->B.$ satisfies the conditions given in the
  statement if and only if its alternating dual $\arrow
  p\dual:A\dual->B\dual.$ satisfies them. Consequently, to demonstrate
  that $\arrow p:A->B.$ is a complicial fibration it is enough to show
  that it has the RLP with respect to all left
  outer horns and thinness extensions, because then we may demonstrate
  the corresponding right handed result for $p$ by appealing to the
  already established left handed one for $\arrow
  p\dual:A\dual->B\dual.$.

  Now to prove that $p$ is a left complicial fibration, first observe
  that the inclusions $\overinc \subseteq_r:E^-_0->E^-_1.$ and
  $\overinc\subseteq_r:\Lambda^0[1]->\Delta^1[1].$ are isomorphic and
  so the statement already assumes left outer horn lifting at
  dimension $1$. At dimension $2$ and above we apply
  observation~\ref{inner.equiv.obs} to replace the left outer horn
  $\overinc\subseteq_r:\Lambda^0[m+2]->\Delta^0[m+2].$ ($m=0,1,2,...$)
  by the isomorphic inclusion
  $\overinc\subseteq_r:X_0->E^-_1\join\Delta[m].$ and then seek to solve the
  lifting problems $(u,v)$ of the form in the diagram
  \begin{equation}\label{lift.prob.y}
    \xymatrix@R=4em@C=12em{
      {X_0}\ar[r]^-{u}\ar@{u(->}[d]_{\subseteq_r} &
      {A}\ar[d]^{p} \\
      {E^-_1\join\Delta[m]}\ar[r]_-{v} 
      \save[]+<4em,3.25em>*+{X_3}="one"\restore
      \save[]+<9.6em,1.8em>*+{E^-_3\join\Delta[m]}="two"\restore
      & {B}
      \ar@{u(->}^-{\subseteq_r}_-{i} "2,1";"one"
      \ar@{u(->}_-{j} "one";"two"
      \ar@{..>}^-{v'} "two";"2,2"
      \ar@{..>}_-{w} "one";"1,2"
    }
  \end{equation}
  in a couple of steps. First we use the weak compliciality of $B$ to
  construct the map $v'$ by extending the stratified map $v$ along the
  inclusion $\overinc\subseteq_r:E^-_1\join\Delta[m]->
  E^-_3\join\Delta[m].$, which is an anodyne extension as it may be
  constructed by applying the anodyne extension preserving right join
  functor $-\join\Delta[m]$ (cf.\ observation~\ref{wcs.decalage}) to
  the anodyne extension $\overinc\subseteq_r:E^-_1->E^-_3.$ (cf.\
  observation~\ref{equiv.decomp.strat}). Now we again consult
  observation~\ref{inner.equiv.obs} to see that the inclusion
  $\overinc\subseteq_r:X_0->X_3.$ is in the cellular closure of the
  set consisting of the elementary inner anodyne extensions and the
  inclusion $\overinc\subseteq_r:E^-_1->E^-_3.$ of the statement, so
  in particular it follows that the assumed injectivity properties of
  $p$ imply that is has the RLP with respect to
  this inclusion. Using this fact we may solve the new lifting problem
  $(u,v'\circ j)$ and obtain the stratified map $w$, which we compose
  with the inclusion $i$ to finally construct the desired solution to
  the original lifting problem.

  An identical argument which, this time, uses the results described
  in observation~\ref{inner.equiv.obs.2} demonstrates that $p$ also
  has the RLP with respect to each elementary
  thinness extension $\overinc\subseteq_e:\Delta^0[m+2]'->
  \Delta[m+2]''.$ ($m=0,1,2,...$). This completes our proof that $p$
  is a left complicial fibration and finally establishes the first
  sentence of the statement, by applying our comments on duality
  above. The second sentence of the statement follows from the first
  simply by observing that $1$ is a weak complicial set and thus
  that $A$ satisfies the conditions of the latter iff the unique map
  $\arrow p:A->1.$ satisfies the conditions of the former.
\end{proof}

\begin{obs}
  Notice that in the last lemma we did not need to explicitly assume
  that $A$ had the RLP with respect to the
  inclusion $\overinc\subseteq_r:E^-_0->E^-_1.$ in order for it to be
  a weak complicial set. Indeed, it is the case that all stratified
  sets have this property, since the unique stratified map from
  $E^-_1$ to $E^-_0\cong\Delta[0]$ is (trivially) left inverse to the
  $\overinc\subseteq_r:E^-_0->E^-_1.$ and may thus be composed with
  any lifting problem $\overarr:E^-_0->X.$ to construct its solution.
\end{obs}

\begin{cor}[lifting of left or right outer horn fillers is enough]
  \label{left+inner=outer.cor}
  Suppose that $B$ is a weak complicial set and $\arrow p:A->B.$ is a
  left complicial fibration then $p$ is a complicial
  fibration. Consequently, $A$ is a weak complicial set iff it is a
  weak left complicial set. Applying these results to the alternating
  duals $\arrow p\dual: A\dual->B\dual.$ and $A\dual$ we also obtain
  the corresponding results for right compliciality.
\end{cor}

\begin{proof}
  By observation~\ref{equiv.decomp.strat} the inclusions
  $\overinc\subseteq_r:E^-_1->E^-_3.$ and
  $\overinc\subseteq_r:E^-_0->E^-_1.$ are both left anodyne
  extensions, so the assumption that $\arrow p:A->B.$ is a left
  complicial fibration implies that it has the RLP
  with respect to those inclusions and thus satisfies the conditions
  of the last lemma.
\end{proof}

\begin{cor}\label{degen.obs}
  If $B$ is a weak complicial set and $\arrow p:A->B.$ is an inner
  complicial fibration then any lifting problem
  \begin{equation}\label{lift.prob.x}
    \xymatrix@R=2em@C=3em{
      {\Lambda^0[m+2]}\ar[r]^-{u}\ar@{u(->}[d]_{\subseteq_r} &
      {A}\ar[d]^p \\
      {\Delta^0[m+2]}\ar[r]_-{v} & {B}
    }
  \end{equation}
  ($m\geq 0$) has a solution {\em so long as\/} $u$ maps the
  $1$-simplex with vertices $0$ and $1$ to a degenerate $1$-simplex in
  $A$. 
\end{cor}

\begin{proof}
  Simply a minor modification of that part of the proof of
  lemma~\ref{inner.equiv.lem} surrounding display~(\ref{lift.prob.y}),
  the details of which we leave to the reader.
\end{proof}

\subsection{Equivalence Inverses}
\label{equiv.inv.subsect}

In this subsection we refine lemma~\ref{inner.equiv.lem} one step
further.

\begin{obs}[illustrating low dimensional calculations]
  In some of what follows, it will be useful to illustrate certain low
  dimensional calculations in our stratified sets. To do so we resort
  to drawing simplices as {\em oriental diagrams}, which were
  introduced by Street in~\cite{Street:1987:Oriental}. It should be
  noted, however, that for us this is simply a convenient way of
  drawing simplices on the $2$-dimensional page rather than a way of
  describing them as free (strict) \infc-categories.

  For example, doing so immediately illuminates the meaning of the the
  RLP with respect to the inclusion
  $\overinc\subseteq_r:E^-_1->E^-_3.$ which was so central to the work
  of the last subsection. Diagrammatically it states that for each
  thin $1$-simplex $v\in tA_1$ there exists some $3$-simplex
  $t=\hat{v}(e^-_3)\in A$ which may be pictured as:
  \begin{equation}\label{equiv.3cocycle}
      \begin{xy}
      \POS(0,0) *[o]{\xybox{\xymatrix@R=2.5em@C=1.5em{
          & {+}\ar[rr]^w && {-}\ar[rd]^v & \\
          {-}\ar[ur]^v\ar[urrr]|(0.4){=}^{}="two" \ar[rrrr]_v^{}="one"
          &&&& {+} \ar@{}"one";"1,4"|(0.4){=}
          \ar@{}"1,2";"two"|{\objectstyle\stackrel{a}{\simeq}} }}}="a"
      \POS(55,0) *[o]{\xybox{\xymatrix@R=2.5em@C=1.5em{ &
          {+}\ar[rr]^w\ar[drrr]|(0.6){=}^{}="one"
          && {-}\ar[rd]^v & \\
          {-}\ar[ur]^v="two"\ar[rrrr]_v^{}="two" &&&& {+}
          \ar@{}"two";"1,2"|(0.4){=}
          \ar@{}"one";"1,4"|{\objectstyle\stackrel{b}{\simeq}} }}}="b"
     \ar@{->}"a";"b"^t_{\simeq}
    \end{xy}
  \end{equation}
  Here we adopt the diagrammatic convention of labelling degenerate
  simplices using the equality symbol $=$ and thin simplices with the
  equivalence symbol $\simeq$. When drawn in this form the intention
  of our definition immediately becomes plain, viz {\em ``$v$ has
    equivalence inverse $w$ and the associated thin $2$-simplices
    $\id_x\simeq w\circ v$ and $\id_y\simeq v\circ w$ have been chosen
    to satisfy a certain $3$-cocycle condition''}.
\end{obs}

\begin{defn}\label{equiv.coher.obs}
  Let $E'_2$ be the stratified set $E^-_2\vee_{E^+_1} E^+_2$, that is
  to say let it be constructed by forming the pushout:
  \begin{equation}\label{lr.equiv.po}
    \xymatrix@=2em{
      {E^+_1}\ar@{u(->}[r]^{\subseteq_r}
      \ar@{u(->}[d]_{\subseteq_r} & {E^-_2}\ar@{u(->}[d]^{i_0} \\
      {E^+_2}\ar@{u(->}[r]_-{i_1} & {E'_2}\poexcursion
    }
  \end{equation}
  Of course, each of the inclusions $i_0$ and $i_1$ is an anodyne
  extension since they are, by definition, pushouts of regular subset
  inclusions which we showed to be anodyne extensions in
  observation~\ref{equiv.decomp.strat}. We will also use $\inc
  i:E^-_1->E'_2.$ to denote the anodyne extension obtained by
  composing $\inc i_0:E^-_2->E'_2.$ of the last sentence and the
  anodyne extension $\overinc\subseteq_r:E^-_1->E^-_2.$. of
  observation~\ref{equiv.decomp.strat}.

  More explicitly, we may represent a stratified map $\arrow
  f:E'_2->A.$ diagrammatically as a pair of thin $2$-simplices
  \begin{equation}\label{lr.equiv.inv}
    \xymatrix@R=2.4em@C=1.2em{
      & {+}\ar[dr]^-w & \\
      {-} \ar[rr]_{=}^{}="one"
      \ar[ur]^-v && {-} 
      \ar@{}"1,2";"one"|<>(0.6){\objectstyle\stackrel{a}{\simeq}}
    }\mkern40mu
    \xymatrix@R=2.4em@C=1.2em{
      & {-}\ar[dr]^-{v'} & \\
      {+} \ar[rr]_{=}^{}="one"
      \ar[ur]^-{w} && {+} 
      \ar@{}"1,2";"one"|<>(0.6){\objectstyle\stackrel{a'}{\simeq}}
    }
  \end{equation}
  in $A$. In other words, this amounts to a $1$-simplex $v$ in $A$
  with a right equivalence inverse $w$ which itself, in turn, has a
  right equivalence inverse $v'$. 

  Notice that the regular subset $E_2\subseteq_r E$ is {\bf not}
  isomorphic to $E'_2$, a fact which follows as soon as we observe
  that a stratified map $\arrow f:E_2->A.$ simply amounts to a pair of
  $1$-simplices $v$ and $w$ which are mutual equivalence inverses.
  It is clear, therefore, that we may construct $E_2$ from $E'_2$ by
  taking a quotient which identifies the $1$-simplices 
  labelled with $v$ and $v'$ in diagram~(\ref{lr.equiv.inv}).
\end{defn}


\begin{lemma}\label{equiv.coher.lem}
  Suppose $B$ is a weak complicial set and that the inner complicial
  fibration $\arrow p:A->B.$ has the RLP with respect to the inclusion
  $\inc i:E^-_1->E'_2.$ then it has the RLP with respect to the
  inclusion $\overinc\subseteq_r:E^-_1->E^-_3.$.
\end{lemma}

\begin{proof}
  Let $C$ denote the stratified set shown in
  figure~\ref{equiv.4.simp}, which is constructed from the
  $0$-trivialised $4$-simplex $\triv_0(\Delta[4])$ by quotienting to
  make degenerate those simplices designated with an $=$ symbol. The
  $1$-simplices labelled $v$ and $w$ and the $2$-simplices labelled
  $a$, $a'$ and $b$ correspond to the simplices of $E'_2$ and $E^-_3$
  labelled in displays~(\ref{lr.equiv.inv}) and~(\ref{equiv.3cocycle})
  respectively, thus allowing us to identify these sets with regular
  subsets of $C$. We've labelled the vertices here with the integers
  used to label the vertices of the original $4$-simplex from which
  $C$ was derived, although of course the quotienting involved in its
  construction means that $0$,$2$ and $3$ actually denote the same
  vertex in there (called $-$) whereas $1$ and $4$ both denote a
  second vertex (called $+$). The remaining simplices have been given
  alphabetic labels in order to discuss them in the arguments that
  follow and to aid the reader in identifying them uniquely in the
  various parts of the diagram in which they are drawn.
 
  \begin{figure}
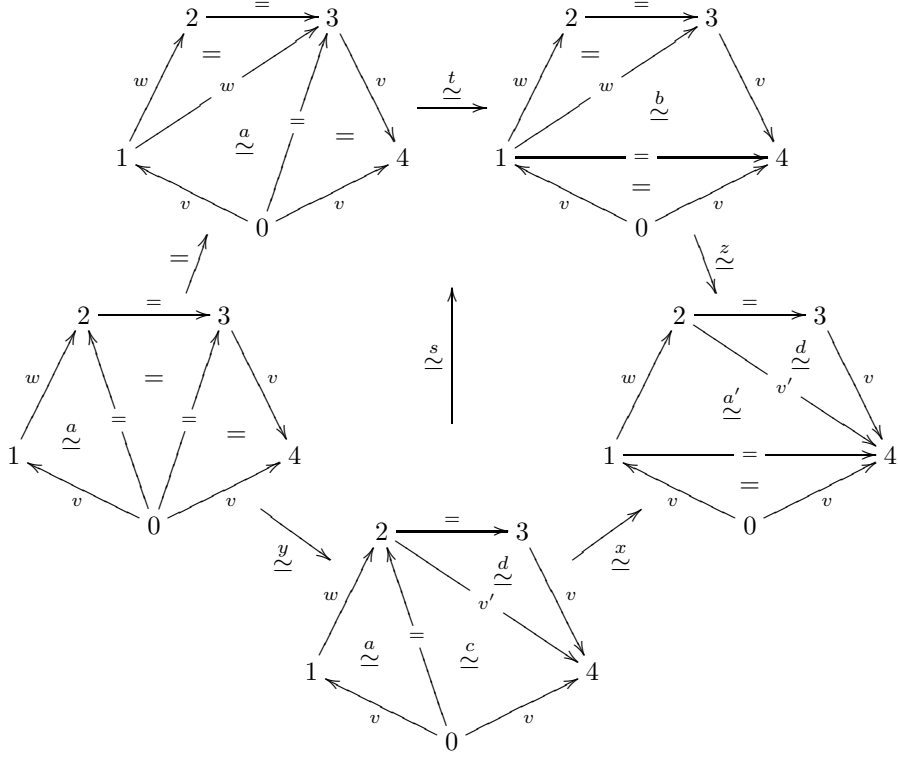

    \begin{equation*}
      \def\a{0}\def\b{1}\def\c{2}\def\d{3}\def\e{4}
      \def\ab{v}\def\bc{w}\def\cd{=}\def\de{v}\def\ae{v}
      \def\ac{=}\def\ad{=}\def\ae{v}\def\bd{w}\def\be{=}
      \def\ce{v'}
      \def\abc{\stackrel{a}\simeq}\def\abd{\abc}\def\abe{=}
      \def\acd{=}\def\ace{\stackrel{c}\simeq}\def\ade{=}
      \def\bcd{=}\def\bce{\stackrel{a'}\simeq}
      \def\bde{\stackrel{b}\simeq}
      \def\cde{\stackrel{d}\simeq}
      \def\abcd{=}\def\abce{\stackrel{x}\simeq}
      \def\abde{\stackrel{t}\simeq}\def\acde{\stackrel{y}\simeq}
      \def\bcde{\stackrel{z}\simeq}
      \def\abcde{\stackrel{s}\simeq}
      \pentofpent{0.9mm}
    \end{equation*}
    \caption{A $4$-simplex}\label{equiv.4.simp}
  \end{figure}

  We start by showing that the inclusion $\overinc\subseteq_r:E'_2->C.$
  enjoys the LLP with respect to any inner complicial fibration
  $\arrow p:A->B.$ whose codomain is a weak complicial set. To do so
  define an increasing sequence of regular subsets of $C$ by
  \begin{equation*}
    U_1 \defeq E'_2 \cup \{ x \}^* \mkern30mu 
    U_2 \defeq U_1 \cup \{ y \}^* \mkern30mu 
    U_3 \defeq U_2 \cup \{ z \}^*
  \end{equation*}
  where the notation $\{{-}\}^*$ denotes the regular subset generated
  by the given set of simplices. We will show that the inclusion of
  each of these in the next may be constructed as a pushout of a thin
  horn extension which has the LLP with respect to $p$ as follows:
  \begin{itemize}[itemsep=1ex, labelindent=2ex, leftmargin=*]
  \item The regular subset $E'_2$ includes those simplices labelled
    $v$, $v'$, $w$, $a$ and $a'$ in figure~\ref{equiv.4.simp}, so in
    particular it includes the data for a $1$-complicial horn on the
    vertices labelled $0,1,2,4$ and it is clear therefore that we may
    construct the inclusion $\overinc\subseteq_r:E'_2->U_1.$ as a
    pushout of the thin inner horn
    $\overinc\subseteq_r:\Lambda^1[3]'->\Delta^1[3]''.$ along the
    evident stratified map which carries the vertices of its domain to
    those labelled $0,1,2,4$ in $E'_2$.
  \item The regular subset $U_1$ contains the $3$-simplex $x$ and thus
    also contains its $2$-face $c$ so we see that this set includes
    the data for a $0$-complicial horn on the vertices labelled
    $0,2,3,4$ and it is thus clear that we may construct
    $\overinc\subseteq_r:U_1->U_2.$ as a pushout of the thin horn
    $\overinc\subseteq_r:\Lambda^0[3]'->\Delta^0[3]''.$ along the
    evident stratified map which carries the vertices of its domain to
    those labelled $0,2,3,4$ in $U_1$. While this is an outer horn,
    the $1$-simplex with vertices labelled $0$ and $2$ is degenerate
    in $C$ and so corollary~\ref{degen.obs} applies in this
    case to show that the inclusion $\overinc\subseteq_r:U_1->U_2.$
    does have the LLP with respect to $p$ as required.
  \item The regular subset $U_2$ contains the $3$-simplex $y$ and its
    $2$-face $d$ so we see that this set includes the data for a
    $1$-complicial horn on the vertices labelled $1,2,3,4$ and it is
    thus clear that we may construct $\overinc\subseteq_r:U_2->U_3.$
    as a pushout of the thin inner horn
    $\overinc\subseteq_r:\Lambda^1[3]'->\Delta^1[3]''.$ along the
    evident stratified map which carries the vertices of its domain to
    those labelled $1,2,3,4$ in $U_2$.
  \item Finally the regular subset $U_3$ contains all of the
    $3$-simplices labelled $x$, $y$ and $z$ and it contains the
    $3$-simplex with vertices labelled $0,2,3,4$, since that is a
    degenerate $3$-simplex with 2-face $a$ which is in $E'_2$, so we
    see that this set includes the data for a $4$-dimensional
    $2$-complicial horn and it is thus clear that we may construct
    $\overinc\subseteq_r:U_3->C.$ as a pushout of the thin inner horn
    $\overinc\subseteq_r:\Lambda^2[4]'->\Delta^2[4]''.$ along the
    evident stratified map which carries the vertices of its domain to
    those labelled $0,1,2,3,4$ in $U_3$.
  \end{itemize}
  So each one of these inclusions enjoys the LLP with respect to $p$
  and it follows therefore that their composite does. 

  Now observe that we may construct a stratified map $\arrow
  r:C->E^-_3.$ which maps the vertex variously labelled $0$, $2$, and
  $3$ in figure~\ref{equiv.4.simp} to $-$ and the one labelled $1$ and
  $4$ to $+$. This is a retraction in the sense that if we pre-compose
  it with the inclusion $\overinc\subseteq_r:E^-_3->C.$ we obtain the
  identity on $E^-_3$. So suppose that $\arrow p:A->B.$ is an inner
  complicial fibration, $B$ is a weak complicial set and that $p$ has
  the RLP with respect to the inclusion $\overinc i:E^-_1->E'_2.$
  and consider the lifting problem depicted as the outer square in:
  \begin{equation*}
    \xymatrix@R=3em@C=10em{
      {E^-_1}\ar[r]^f\ar@{u(->}[d]_{\subseteq_r} 
      \save []+<4em,-1.2em>*+{E'_2}="two"\restore &
      {A}\ar[d]^p \\
      {E^-_3}\ar[r]_g 
      \save []+<6em,1.4em>*+{C}="one" \restore & {B}
      \ar"one";"2,2"^{g\circ r}\ar@{u(->}"2,1";"one"
      \ar@{u(->}"two";"one"\ar@{d(->}"1,1";"two"_{i}
      \ar@{..>}"two";"1,2"_(0.35){h'}\ar@{..>}"one";"1,2"_{h}  
    }
  \end{equation*}
  To solve this, we first form the composite $g\circ r$ which
  completes the data for a lifting problem from the inclusion
  $\inc i:E^-_1->E'_2.$ to $\arrow p:A->B.$,
  because $r$ is a retract of the inclusion from $E^-_3$ to $C$ and
  the (skewed) square of inclusions to the left of our diagram
  commutes. Solving this problem, which we may do by assumption on
  $p$, we obtain the map $h'$ which in turn furnishes us with a
  lifting problem $(h',g\circ r)$ from $\overinc\subseteq_r:E'_2->C.$
  to $p$. However the proof of the last paragraph tells us that these
  enjoy the lifting property with respect to each other, so we may
  solve this latter problem to obtain the map $h$ and thus solve our
  original problem with the composite of that map and the inclusion
  $\overinc\subseteq_r:E^-_3->C.$ as required.
\end{proof}

As a corollary, we find that a weak inner complicial set is actually a
weak complicial set if and only if each of its thin $1$-simplices has
an equivalence inverse:

\begin{cor}\label{equiv.coher.cor1}
  If $A$ is a weak inner complicial set then it is a weak complicial
  set iff it has the RLP with respect to the inclusion
  $\overinc\subseteq_r:E^-_1->E_2.$.
\end{cor}

\begin{proof}
  For the ``only if'' implication, if $A$ is a weak complicial set
  then we may lift any stratified map $\arrow f:E^-_1->A.$ along the
  anodyne extension $\overinc\subseteq_r:E^-_1->E^-_3.$ and compose
  the resulting map with the inclusion
  $\overinc\subseteq_r:E_2->E^-_3.$ to construct a stratified map
  which provides the required lift of $f$ along
  $\overinc\subseteq_r:E^-_1->E_2.$.

  For the reverse implication, we know that $E_2$ is a quotient of
  $E'_2$ and that we may decompose the inclusion
  $\overinc\subseteq_r:E^-_1->E_2.$ as the composite of the inclusion
  $\inc i:E^-_1->E'_2.$ and the quotient map $\epi q:E'_2->E_2.$. So
  if $A$ has the lifting property of the statement then we may show
  that it has the RLP with respect the inclusion $i$ by lifting along
  $\overinc\subseteq_r:E^-_1->E_2.$ and composing with $q$. If follows
  that we may apply lemma~\ref{equiv.coher.lem} to show that $A$ also
  has the RLP with respect to the inclusion
  $\overinc\subseteq_r:E^-_1->E^-_3.$ which, in turn, allows us to
  apply lemma~\ref{inner.equiv.lem} and demonstrate that it is a weak
  complicial set as required.
\end{proof}

\begin{cor}\label{equiv.coher.cor2}
  If $A$ and $B$ are weak complicial sets and $\arrow p:A->B.$ is an
  inner complicial fibration then it is a complicial fibration iff it
  has the RLP with respect to the inclusion
  $\overinc\subseteq_r:E^-_0->E^-_1.$.
\end{cor}

\begin{proof}
  The ``only if'' direction is immediate, since the cited inclusion is
  isomorphic to the elementary anodyne extension
  $\overinc\subseteq_r:\Lambda^0[1]->\Delta^0[1].$. For the reverse
  implication we start by demonstrating that $p$ has the RLP with
  respect to the inclusion
  $\overinc\subseteq_r:\Lambda^0[2]->\Delta^0[2].$. To that end define
  three stratified sets
  \begin{displaymath}
    W_1\defeq E^-_1\join\Delta[0] \mkern30mu
    W_2\defeq (E^-_0\join\Delta[0])\cup(E^-_2\join\boundary\Delta[0]) 
    \mkern30mu W_3\defeq E^-_2\join\Delta[0]
  \end{displaymath}
  for which $W_1,W_2\subseteq_r W_3$ and
  \begin{align*}
    W_1\cap W_2 & {}=(E^-_0\join\Delta[0])\cup(E^-_1\join
    \boundary\Delta[0]) \\
    W_1\cup W_2 & {}=(E^-_1\join\Delta[0])\cup(E^-_2\join
    \boundary\Delta[0])
  \end{align*}
  consequently, arguing as in observation~\ref{join.anodyne}, we find
  that the horn inclusion of the last sentence is isomorphic to
  $\overinc\subseteq_r:W_1\cap W_2-> W_1.$. To show that $p$ has the
  RLP with respect to this latter inclusion consider the lifting
  problem $(f,g)$ shown in the outer square of the following diagram:
  \begin{displaymath}
    \xymatrix@R=4em@C=12em{
      {W_1\cap W_2} 
      \save []+<5.5em,-1.4em>*+{W_2}="w2"\restore
      \save []+<11em,-3em>*+{W_3}="w3"\restore
      \ar[r]^f\ar@{u(->}[d]_{\subseteq_r} &
      {A}\ar[d]^{p} \\
      {W_1}\save [] +<5.5em,1.4em>*+{W_1\cup W_2}="w1uw2"\restore
      \ar[r]_g & {B}
      \ar@{d(->}"1,1";"w2"_-{\subseteq_r} 
      \ar@{u(->}"2,1";"w1uw2"^-{\subseteq_r}
      \ar@{u(->}"w2";"w1uw2"_-{\subseteq_r}
      \ar@{u(->}"w2";"w3"^-{\subseteq_r}
      \ar@{u(->}"w1uw2";"w3"\ar@{..>}"w2";"1,2"_h
      \ar@{..>}"w3";"2,2"^k\ar@{..>}"w1uw2";"2,2"^l
      \ar@{..>}"w3";"1,2"_m
    }
  \end{displaymath}
  Here, we construct the various dotted maps in the following sequence:
  \begin{itemize}[itemsep=1ex, labelindent=2ex, leftmargin=*]
  \item The inclusion $\overinc\subseteq_r:W_1\cap W_2->W_2.$ may be
    constructed by taking the pushout along the inclusion
    $\overinc\subseteq_r:E^-_1\join\boundary\Delta[0]-> W_2.$ of the
    corner join of the anodyne extension
    $\overinc\subseteq_r:E^-_1->E^-_2.$ and the inclusion
    $\overinc\subseteq_r:\emptyset->\boundary\Delta[0].$. So, applying
    lemma~\ref{wcs.decalage}, we find that this inclusion is an
    anodyne extension and thus that we may factor the map $f$ through
    $W_2$, since $A$ is a weak complicial set by assumption, to give
    the map labelled $h$.
  \item Now the map $l$ may be constructed using the pushout property
    of the pasting square determined by $W_1$ and $W_2$ to the left of
    the diagram.
  \item The inclusion $\overinc\subseteq_r:W_1\cup W_2->W_3.$ is
    simply the corner join of the anodyne extension
    $\overinc\subseteq_r: E^-_1->E^-_2.$ (cf.\
    observation~\ref{equiv.decomp.strat}) and the inclusion
    $\overinc\subseteq_r:\boundary\Delta[0]->\Delta[0].$. Applying
    lemma~\ref{wcs.decalage}, we find that this inclusion is an
    anodyne extension and thus that we may factor the map $l$ through
    $W_3$, since $B$ is a weak complicial set by assumption, to give
    the map labelled $k$.
  \item The inclusion $\overinc\subseteq_r:W_2->W_3.$ is the
    corner join of the right anodyne extension
    $\overinc\subseteq_r:E^-_0->E^-_2.$ (cf.\
    observation~\ref{equiv.decomp.strat}) and the inclusion
    $\overinc\subseteq_r:\boundary\Delta[0]->\Delta[0].$. Applying
    observation~\ref{inner.joins}, we find that this inclusion is an
    {\it inner\/} anodyne extension, and thus that we may solve the
    lifting problem $(h,k)$ into $p$, since this latter map is an
    inner complicial fibration by assumption, to give the map $m$.
  \end{itemize}
  So we obtain the desired solution to our original lifting problem
  $(f,g)$ by composing $m$ with the inclusion
  $\overinc\subseteq_r:W_1->W_2.$. However we know, by the comment in
  definition~\ref{thinness.ext} and the fact that $A$ is a weak
  complicial set, that $p$ has the RLP with respect to the elementary
  thinness extension $\overinc\subseteq_e:\Delta^0[2]'->
  \Delta^0[2]''.$. Combining this with the lifting property
  established above, it follows that $p$ also has the RLP with respect
  to the thin horn inclusion $\overinc\subseteq_r:\Lambda^0[2]'->
  \Delta^0[2]''.$.

  Now we know that the inclusion $\inc i:E^-_1->E'_2.$, of
  definition~\ref{equiv.coher.obs} is constructed by composing
  $\overinc\subseteq_r:E^-_1->E^-_2.$ with a pushout of
  $\overinc\subseteq_r:E^+_1->E^-_2.$ and that each of these latter
  inclusions may be constructed as a pushout of
  $\overinc\subseteq_r:\Lambda^0[2]'->\Delta^0[2]''.$ as in
  observation~\ref{equiv.decomp.strat}. So it follows immediately,
  from the result of the last paragraph, that $p$ has the RLP with
  respect to $\inc i:E^-_1->E'_2.$ and that we may thus apply
  lemma~\ref{equiv.coher.lem} to show that it also has the RLP with
  respect to $\overinc\subseteq_r:E^-_1->E^-_3.$. Finally the
  assumption that $p$ also has the RLP with respect to
  $\overinc\subseteq_r:E^-_0->E^-_1.$ allows us to
  apply lemma~\ref{inner.equiv.lem} and show that $p$ is a complicial
  fibration as required.
\end{proof}

\begin{thm}\label{almost.weakinner.cor}
  Suppose that the stratified set $A$ is {\bf almost} a weak inner
  complicial set, in the sense that we insist that it has the RLP with
  respect to all inner elementary anodyne extensions except
  $\overinc\subseteq_e:\Delta[2]'-> \Delta[2]''.$. Furthermore,
  suppose that its set of thin $1$-simplices is the subset
  \begin{equation}\label{equivs.thin}
    \left\{ 
      v\in A_1\mid \exists \arrow f:\triv_1(\usimp{E}_2)->A.\text{
        with } f(e^-_1)=v
    \right\}
  \end{equation}
  of those $1$-simplices with equivalence inverses then $A$ is a weak
  complicial set.
\end{thm}

\begin{proof}
  The involution $\arrow\neg:E->E.$ of observation~\ref{iso.symm}
  restricts to an involution $\arrow\neg:E_2->E_2.$ which carries
  $e^-_1$ to $e^+_1$. So it follows, from the description of the thin
  $1$-simplices of $A$ given in the statement, that if $\arrow
  f:\triv_1(\usimp{E}_2)->A.$ is a stratified map then both of
  $f(e^-_1)$ and $f(e^+_1)=(f\circ\neg)(e^-_1)$ are thin in $A$ and
  thus that $A$ has the RLP with respect to the inclusion
  $\overinc\subseteq_r:E^-_1->E_2.$. Consequently, we may apply
  corollary~\ref{equiv.coher.cor1} to show that $A$ is a weak
  complicial set so long as we have verified that it is a weak inner
  complicial set, that is we need to show that it also has the RLP
  with respect to the elementary thinness extension
  $\overinc\subseteq_e:\Delta^1[2]'->\Delta^1[2]''.$. 

  In other words, we must demonstrate that if $c\in A$ is a thin
  $2$-simplex and its $1$-dimensional faces $v_0\defeq c\cdot\face_0$
  and $v_2\defeq c\cdot\face_2$ are both thin then so is $v_1\defeq
  c\cdot\face_1$. However, we know that the $1$-simplex $v_0$ (resp.\
  $v_2$) is thin in $A$ iff we have have corresponding thin
  $2$-simplices $a_0$, $a'_0$ (resp.\ $a_2$, $a'_2$) as depicted in
  display~\eqref{lr.equiv.inv}, so our task will be to construct
  similar thin simplices $a_1$ and $a'_1$ for $v_1$.

  \begin{figure}
    \begin{equation*}
      \def\a{0}\def\b{1}\def\c{2}\def\d{3}\def\e{4}
      \def\ab{v_2}\def\bc{v_0}\def\cd{w_0}\def\de{w_2}\def\ae{}
      \def\ac{v_1}\def\ad{}\def\ae{=}\def\bd{=}\def\be{w_2}
      \def\ce{?w_1}
      \def\abc{\stackrel{c}\simeq}\def\abd{=}\def\abe{\stackrel{a_2}\simeq}
      \def\acd{}\def\ace{\stackrel{?a_1}\simeq}\def\ade{\stackrel{a_2}\simeq}
      \def\bcd{\stackrel{a_0}\simeq}\def\bce{\stackrel{?e}\simeq}
      \def\bde{=}
      \def\cde{\stackrel{?d}\simeq}
      \def\abcd{}\def\abce{\stackrel{?t}\simeq}
      \def\abde{=}\def\acde{}
      \def\bcde{\stackrel{?s}\simeq}
      \def\abcde{}
      \def\baselen{0.9mm}
      \begin{xy}
        0;<\baselen,0mm>:
        \POS(0,30)*[o]{\pent 1{\baselen}}="one"
        \POS(45,0)*[o]{\pent 4{\baselen}}="four"
        \POS(90,30)*[o]{\pent 2{\baselen}}="two"
        \ar "one";"four"_{\objectstyle\bcde}
        \ar "two";"four"^{\objectstyle\abce}
      \end{xy}
    \end{equation*}
    \caption{Composing Thin $1$-Simplices in $A$\label{comp.1-simp}}
  \end{figure}
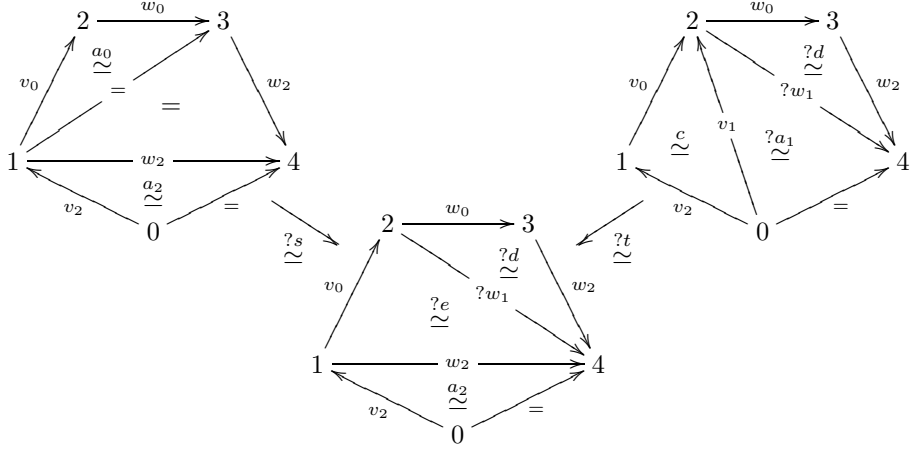

  We illustrate the construction of these in figure~\ref{comp.1-simp}
  which depicts a stratified map with domain
  $X\subseteq_r\triv_1(\Delta[4])$ the regular subset generated by the
  $3$-simplices $\face_0$ and $\face_3$ and with codomain $A$. To aid
  our discussion $0$-simplices in this diagram have been labelled to
  identify them in the domain $X$ whereas all other simplices take
  names intended to represent simplices of the codomain
  $A$. Furthermore, a question mark appended to the front of a simplex
  label indicates that the corresponding simplex will be constructed
  in $A$ by filling some complicial horn. We identify other simplices
  in $X$ by listing their vertices so, for instance, $014$ denotes the
  (unique) $2$-simplex whose $0$-dimensional faces (vertices) are $0$,
  $1$ and $4$.

  We commence the building of this map by initialising the left hand
  pentagon with the data we are given and then working rightward,
  filling complicial horns as we go.  So we map $123$ to $a_0$, $014$
  to $a_2$ and $134$ to the degenerate simplex $w_2\cdot\degen_0$,
  whose faces are mutually compatible as shown in the diagram.  Now in
  the middle pentagon we may fill the $1$-complicial horn on vertices
  $2,3,4$ thereby constructing a mapping of $234$ to the thin
  $2$-simplex $d\in A$ and obtaining a new $1$-face $w_1\in A$. This
  completes the data for a $2$-complicial horn on the vertices
  $1,2,3,4$, which we fill to map $1234$ to a thin $3$-simplex $s\in
  A$ and obtain a new $2$-face $e\in A$ (which is thin since all the
  other $2$-faces of $s$ are thin). Finally, we may map $012$ to the
  thin $2$-simplex $c\in A$ that we started with in the second
  paragraph of this proof.  In doing so we complete the data for a
  $1$-complicial horn on the vertices $0,1,2,4$, which we fill to map
  $0124$ to a thin $3$-simplex $t\in A$ and obtain a new $2$-face
  $a_1\in A$ (which is thin since all the other $2$-faces of $t$ are
  thin). This is the thin $2$-simplex we seek, witnessing that $w_1$
  is a left equivalence inverse of $v_1$. Dually we may replay the
  construction above in the alternating dual $A\dual$ to derive a thin
  $2$-simplex $a'_1$ which demonstrates that the $1$-simplex $v'_1$
  obtained by ``composing'' $v'_0$ and $v'_2$ is a left equivalence
  inverse of $w_1$ as required.
\end{proof}

\begin{eg}[quasi-categories as weak complicial sets]\label{quasi.wcs}
  We are now in a position to validate example~\ref{kan.joyal.eg} by
  showing that any quasi-category $A$ may be given a stratification
  that makes it into a $1$-trivial weak complicial set.  

  First note that the quasi-categorical inner horn filler conditions
  simply translate to postulating that the $1$-trivialisation
  $\triv_1(A)$ is almost a weak inner complicial set (in the sense of
  the last corollary). This property places no restriction on thin
  $1$-simplices, so we may extend the stratification of $\triv_1(A)$
  without disrupting it by making thin all $1$-simplices in the subset
  shown in display~\eqref{equivs.thin}, thereby giving a stratified
  set we denote by $A^e$. It follows that we may apply
  theorem~\ref{almost.weakinner.cor} to this latter stratification to
  show that $A^e$ is a $1$-trivial weak complicial set as suggested.

  Consequently, whenever we speak of quasi-categories in future we
  will implicitly assume that they carry the stratification of the
  last paragraph. The construction is clearly functorial, thereby
  demonstrating that we may identify the category of quasi-categories
  with a certain full subcategory of the category of $1$-trivial weak
  complicial sets.  Indeed theorem~\ref{almost.weakinner.cor} tells
  us, amongst other things, that we may characterise the objects of
  this full subcategory as being those $1$-trivial weak complicial
  sets $A$ which have the RLP with respect to the inclusion
  $\overinc\subseteq_e:\triv_1( \usimp{E}_2)->E_2.$.
\end{eg}

\section{Gray Tensor Products}
\label{gray.tensor.sec}

In this section we generalise the complicial theory of Gray tensor
products and their closures, as presented in
\cite{Verity:2006:Complicial}, to the weak complicial context. While
many of the proofs given there generalise directly we still feel that
an independent presentation is warranted here, since  it simplifies some
aspects of the strict theory and recasts it more clearly as a piece of
homotopy theory.

\subsection{Some Tensor Products of Stratified Sets}

\begin{obs}[a motivating analogy with bicategory theory]
  In the theory of bicategories, as explicated by Street
  in~\cite{Street:1980:FibBicat}, the (strict, algebraic) cartesian
  product of bicategories makes the category of bicategories and
  homomorphisms (pseudo-functors) into a symmetric monoidal category.
  Given a pair of bicategories $\bicat{B}$ and $\bicat{C}$ we may form
  the bicategory $\Hom(\bicat{B},\bicat{C})$ of homomorphisms, strong
  transformations (pseudo-naturals) and modifications between them,
  which provides this monoidal structure with a {\em weak closure\/}
  in the sense that there is a canonical biequivalence
  $\Hom(\bicat{B}\times\bicat{C},\bicat{D})\simeq_b\Hom(\bicat{B},
  \Hom(\bicat{C},\bicat{D}))$. In other words, in bicategory theory
  the cartesian product takes the role of the Gray tensor product in
  $2$-category theory. This insight motivates the next three
  definitions:
\end{obs}

\begin{defn}[Gray tensor product of stratified sets]
  \label{gray.tensor.def}
  The {\em Gray tensor product\/} of stratified sets $X$ and $Y$ is
  simply defined to be their cartesian product $X\gray Y$ in the
  category $\Strat$. Explicitly, $X\gray Y$ is the stratified set
  whose $n$-simplices are pairs of $n$-simplices $(x,y)$ with $x\in
  X_n$ and $y\in Y_n$, whose simplicial action is given pointwise
  $(x,y)\cdot\alpha=(x\cdot\alpha,y\cdot\alpha)$ and whose thin
  simplices are those $(x,y)$ with $x$ thin in $X$ {\bf and} $y$ thin
  in $Y$. 

  We have two reasons for not adopting the usual cartesian product
  notation $\times$ for the Gray tensor product. Firstly we would like
  to stress that we are primarily interested in regarding this as the
  appropriate generalisation of the 2-categorical Gray tensor product,
  the fact that it actually coincides with the categorical product in
  this context is an important but secondary fact. Secondly, it helps
  us to avoid certainly notational difficulties which might arise when
  manipulating simplicial sets $X$ and $Y$ under the minimal
  stratification convention, since then it is not the case that the
  minimal stratification of their simplicial cartesian product
  $X\times Y$ coincides with the stratification of their Gray tensor
  product $X\gray Y$ as stratified sets. In other words, the minimal
  stratification operation does not preserve cartesian products.

  Observation~\ref{LFP.quasitopos} reminds us that $\Strat$ is a
  quasi-topos, so in particular it is cartesian closed with closure
  (function space construction) between stratified sets $X$ and $Y$
  denoted by $\hom(X,Y)$. This is often referred to as the {\em
    stratified set of strong transformations}, since it is the true
  weak complicial analogue of the bicategory theorist's bicategory of
  homomorphisms, strong transformations and modifications. We also
  denote the corresponding corner product and closure by $\cgray$
  and $\chom$ respectively (cf.\ recollection~\ref{corner.tensor}).
\end{defn}

\begin{defn}[partition operators]\label{part.op}
  We say that a pair $p,q\in\mathbb{N}$ is a {\em partition\/} of
  $n\in\mathbb{N}$ if $p+q=n$ and associate with it four {\em
    partition operators}:
  \begin{itemize}
  \item face operators $\arrow\partinj^{p,q}_1:[p]->[n].$ given by
    $\partinj^{p,q}_1(i)=i$ and $\arrow\partinj^{p,q}_2:[q]->[n].$
    given by $\partinj^{p,q}_2(j)=j+p$, and
  \item degeneracy operators $\arrow\partproj^{p,q}_1:[n]->[p].$ and
    $\arrow\partproj^{p,q}_2:[n]->[q].$ given by
    \begin{equation*}
      \partproj^{p,q}_1(i) =
      \begin{cases}
        i & \text{when $i\leq p$} \\
        p & \text{when $i>p$}
      \end{cases}\mkern20mu\text{and}\mkern20mu
      \partproj^{p,q}_2(i) =
      \begin{cases}
        0 & \text{when $i< p$} \\
        i - p & \text{when $i\geq p$}
      \end{cases}
    \end{equation*}
    respectively.
  \end{itemize}
\end{defn}

\begin{defn}[associative lax Gray tensor product of stratified sets]
  \label{laxgray.tensor.def}
  The {\em (associative) lax Gray tensor product\/} $X\laxgray Y$ of
  stratified sets $X$ and $Y$ (definition~128
  of~\cite{Verity:2006:Complicial}) is formed by taking the product of
  underlying simplicial sets and endowing it with the stratification
  under which the $n$-simplex $(x,y)$ is thin in $X\laxgray Y$ iff for
  each partition $p,q$ of its dimension we either have that
  $x\cdot\partinj^{p,q}_1$ is a thin $p$-simplex in $X$ or that
  $y\cdot\partinj^{p,q}_2$ is a thin $q$-simplex in $Y$.

  Notice, in particular, that by definition stratifications can have
  no thin $0$-simplices so this condition applied to the extremal
  partitions $n,0$ and $0,n$ imply that if $(x,y)$ is thin in
  $X\laxgray Y$ then $x$ is thin in $X$ and $y$ is thin in $Y$. In
  other words, $X\laxgray Y$ is an entire subset of $X\gray Y$.
\end{defn}

\begin{obs}[lax Gray tensors in strict complicial set theory]
  In the theory of (strict) complicial sets presented
  in~\cite{Verity:2006:Complicial}, the relationship between the lax
  Gray tensors of stratified sets and (strict) \infc-categories extends
  well beyond mere analogy. Indeed, section~11.4 of that work
  demonstrates that the lax Gray tensor product of (strict)
  \infc-categories, as defined by Steiner~\cite{Steiner:1991:Tensor} or
  Crans~\cite{Crans:1995:PhD}, may be obtained by reflecting the lax
  Gray tensor of stratified sets to the equivalent subcategory of
  (strict) complicial sets.
\end{obs}

\begin{obs}\label{laxgray.obs}
  The primary properties of $\laxgray$ as a tensor product on $\Strat$
  are established in lemmas~129 and~131
  of~\cite{Verity:2006:Complicial} and may be derived directly from
  the {\em partition identities\/} between partition operators given
  in notation~5 of that work. In summary, $\laxgray$ may be extended
  to stratified maps and equipped with canonical associativity and
  identity isomorphisms which make it into the tensor of a
  (non-symmetric) monoidal structure on $\Strat$. This structure is
  completely characterised by the fact that the forgetful underlying
  simplicial set functor becomes a strict monoidal functor from the
  monoidal category $(\Strat,\laxgray,\Delta[0])$ to the cartesian
  closed category $(\Simp,\times,\Delta[0])$. Furthermore, $\laxgray$
  is well behaved with respect to alternating duals, with the ``swap''
  map on underlying simplicial products providing us with canonical
  isomorphisms $(X\laxgray Y)\dual\cong Y\dual\laxgray X\dual$.

  However, as discussed in observation~136
  of {\em loc.\ cit.}, while $\laxgray$ provides $\Strat$
  with a genuine monoidal structure it fails to be well behaved with
  respect to colimits of stratified sets. This leads us to define the
  following, closely related, tensor product for which left and right
  tensoring does preserve colimit but which fails to be coherently
  associative.  To simplify our presentation here a little the next
  definition doesn't quite follow that of the corresponding
  construction of that work, but nevertheless shares all of its
  important properties.
\end{obs}

\begin{defn}[lax Gray pre-tensor product of stratified sets]
  \label{pretens.defn}
  The {\em lax Gray pre-tensor product\/} $X\pretens Y$ of stratified
  sets $X$ and $Y$ (definition~135 of~\cite{Verity:2006:Complicial})
  is formed by taking the product of underlying simplicial sets and
  endowing it with a stratification under which an $r$-simplex
  $(x,y)\in X\times Y$ is thin if either
  \begin{itemize}
  \item there exists $0<k<r$ such that $x=x'\cdot\degen_{k-1}$ and
    $y=y'\cdot\degen_k$ for some pair of simplices $x'\in X$ and
    $y'\in Y$, or
  \item there exists a partition $p,q$ of its dimension and simplices
    $x'\in X_p$ and $y'\in Y_q$ such that $x=x'\cdot\partproj^{i,j}_1$ and
    $y=y'\cdot\partproj^{i,j}_2$ and either $x'$ is thin in $X$ or
    $y'$ is thin in $Y$.
  \end{itemize}
  It is easily demonstrated that this is a stratification which makes
  $X\pretens Y$ into an entire subset of $X\laxgray Y$.
\end{defn}

\begin{obs}\label{pretens.gen.obs}
  Here again, it is easily shown that the action of cartesian product
  as a bifunctor of underlying simplicial sets may be lifted to make
  $\pretens$ into a bifunctor on $\Strat$. This time, however, it is
  {\bf not} the case that canonical associativity isomorphisms also
  lift in this way, but it is still true that identity isomorphism
  lift to give $X\pretens\Delta[0]\cong X\cong \Delta[0]\pretens X$
  and that the ``swap'' map provides a canonical isomorphism
  $(X\pretens Y)\dual\cong Y\dual\pretens X\dual$.

  Most importantly, lemma~142 of~\cite{Verity:2006:Complicial}
  demonstrates that the pre-tensor $\pretens$ preserves colimits in
  each variable. Consequently, since $\Strat$ is locally finitely
  presentable, it follows that it possesses closures $\lax_l(X,Z)$ and
  $\lax_r(Y,Z)$ which are right adjoint to the endo-functors
  $X\pretens -$ and $-\pretens Y$ respectively. We often call these
  the {\em stratified set of left (resp.\ right) lax
    transformations\/} since they generalise the bicategory theorist's
  bicategories of homomorphisms, left (resp.\ right) lax
  transformations and modifications. Again we adopt the notations
  $\cpretens$, $\clax_l$ and $\clax_r$ to denote the corresponding
  corner tensor and its closures (cf.\
  recollection~\ref{corner.tensor}).
\end{obs}

\begin{obs}\label{pretens.tens.comp}
  Finally, it is also worth pointing out that lemma~139
  of~\cite{Verity:2006:Complicial} demonstrates that for each pair of
  stratified sets $X$ and $Y$ the entire inclusion
  $\overinc\subseteq_e:X\pretens Y->X\laxgray Y.$ is an inner anodyne
  extension. In other words, we might say that weak complicial sets
  ``do not see'' the difference between the tensors $\pretens$ and
  $\laxgray$. 

  We may also define a bifunctor
  $\arrow\spretens:\Strat\times\Strat->\Strat.$ for which $X\spretens Y$
  is the entire superset of $X\pretens Y$ constructed by making thin
  all simplices of the form $(x\cdot\partproj^{r,s}_1,
  y\cdot\partproj^{r,s}_2)$ with $x\in X$, $y\in Y$ and
  $r,s>0$. Notice that each of these simplices is thin in $X\gray Y$,
  so we have an entire inclusion $\overinc\subseteq_e:X\spretens Y->
  X\gray Y.$. Furthermore, as the reader may readily verify, we may
  modify the proof given in lemma~139 of loc.\ sit.\ to show that this
  inclusion is also an inner anodyne extension.
\end{obs}

\subsection{Gray Tensors and Anodyne Extensions}

Our primary interest in the remainder of this section will be to
demonstrate that these tensors are well behaved with respect to
certain anodyne extensions. In the process we show that the functors
$\hom(X,{-})$, $\lax_l(X,{-})$ and $\lax_r(Y,{-})$ all preserve weak
(inner) compliciality.

\begin{defn}\label{shuffle.depth}
  As ever, the non-degenerate $(n+m)$-simplices of the simplicial set
  $\Delta[n]\times\Delta[m]$ are called {\em shuffles}. An easy and
  useful characterisation of these is that they are precisely the
  $(n+m)$-simplices $(\alpha,\beta)$ which satisfy the {\em ordinate
    summation property\/} which states that $\alpha(i)+\beta(i)=i$ for
  all $i\in[n+m]$.  We define the {\em depth\/} (cf.\ Porter and
  Ehlers~\cite{Porter:2005:SpecialPasting}) of such a shuffle to be
  the integer:
  \begin{equation*}
    \depth(\alpha,\beta)\defeq\sum_{i=0}^{n+m} \min(\alpha(i), m-\beta(i))
  \end{equation*}
\end{defn}

\begin{obs}[more about shuffles]\label{shuffle.facts}
  We may depict a shuffle in $\Delta[n]\times\Delta[m]$ as a path of
  strictly horizontal (rightward) and vertical (upward) moves on an
  $[n]\times[m]$ grid, which starts at its bottom left corner and ends
  at its top right one. Then, as observed in
  \cite{Porter:2005:SpecialPasting}, the depth of that shuffle is
  simply the number of squares of that grid which occur to the left of
  and above that path. For instance the depth of the example (solid
  line) in figure~\ref{shuffle.pic} is equal to the number of squares
  that have been outlined with dotted boundaries, which in this case
  is $8$.
  
  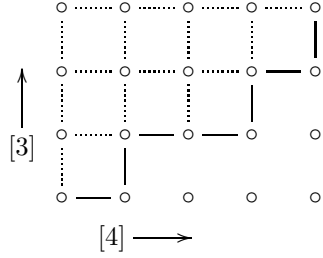
\begin{figure}
    \begin{equation*}
      \xymatrix@!=0em{
        {\circ} \ar@{..} '[r] '[rr] '[rrr] [rrrr] 
        \ar@{..} '[d] '[dd] [ddd]
        & {\circ} \ar@{..} '[d] [dd]
        & {\circ} \ar@{..} '[d] [dd]
        & {\circ} \ar@{..} [d]
        & {\circ} \\
        {\circ} \ar@{..} '[r] '[rr] [rrr] 
        & {\circ} & {\circ} & {\circ} & {\circ} \\
        {\circ}  \ar@{..} [r] & {\circ} & {\circ} & {\circ} & {\circ} \\
        {\circ} \ar@{-} '[r] '[ru] '[ru] '[rur] '[rurr] '[rurru]
        '[rurrur] '[rurrur] [rurruru]
        & {\circ} & {\circ} & {\circ} & {\circ}
        \save "4,2"+<-0.5em,-1.5em>*+{[4]}="one"
        \ar "one";"one"+<3em,0em> \restore
        \save "3,1"+<-1.5em,-0.5em>*+{[3]}="two"
        \ar "two";"two"+<0em,3em> \restore
      }
    \end{equation*}
    \caption{A shuffle in
      $\Delta[4]\times\Delta[3]$\label{shuffle.pic}}
  \end{figure}
  
  In line with this depiction, we make the following simple
  observations and definitions:
  \begin{enumerate}[label=(\alph*)]
  \item\label{shuffle.facts.a} The only depth $0$ shuffle is the one
    which would be depicted as a sequence of $m$ upward moves followed
    by $n$ rightward ones, in other words the simplex
    $(\partproj^{m,n}_2,\partproj^{m,n}_1)$.
  \item\label{shuffle.facts.b} All shuffles have depth less than or
    equal to $nm$ and the only depth $nm$ shuffle is the one which
    would be depicted as a sequence of $n$ rightward moves followed by
    $m$ upward ones, in other words the simplex
    $(\partproj^{n,m}_1,\partproj^{n,m}_2)$.
  \item\label{shuffle.facts.c} If $(\alpha,\beta)$ is a shuffle we say
    that its $i\oth$ vertex (for $0<i<n+m$) is a {\em left-upper
      corner\/} if $\alpha(i-1)=\alpha(i)$ and $\beta(i)=\beta(i+1)$
    and we say it is a {\em right-lower corner\/} if
    $\beta(i-1)=\beta(i)$ and $\alpha(i)=\alpha(i+1)$. These are
    simply the right and left handed right angle turning points in its
    depiction.
  \item\label{shuffle.facts.d} The only shuffle with no left-upper
    corners is the maximal depth shuffle
    $(\partproj^{n,m}_1,\partproj^{n,m}_2)$. Dually the only shuffle
    with no right-lower corners is the minimal depth shuffle
    $(\partproj^{m,n}_2,\partproj^{m,n}_1)$.
  \item\label{shuffle.facts.e} No two left-upper (resp.\ right-lower)
    corners of $(\alpha,\beta)$ can be immediately adjacent, that is
    to say if $i<j$ are indices of two such left upper corners then we
    actually have $i+1<j$.
  \item\label{shuffle.facts.f} The $i\oth$ vertex of our shuffle is
    neither a left-upper nor a right-lower corner iff the face
    $(\alpha,\beta)\cdot\face_i$ is a simplex of
    $\boundary(\Delta[n]\times\Delta[m])$.
  \end{enumerate}
\end{obs}

\newcommand{\nc}{_\bullet}

\begin{notation}\label{P.notation}
  In the next lemma we will assume that $P$ is a stratified set with
  underlying simplicial set $\Delta[n]\times\Delta[m]$ and which
  satisfies the condition that whenever $(\phi,\psi)$ is a
  non-degenerate $r$-simplex of $P$ and $l$ is some integer with
  $0<l<r$ such that $\phi(l-1)=\phi(l)$ and $\psi(l)=\psi(l+1)$
  (upper-left corner) then:
  \begin{enumerate}[label=(\alph*)]
  \item\label{P.notation.a} $(\phi,\psi)$ is thin in $P$, and
  \item\label{P.notation.b} if the face $(\phi,\psi)\cdot\face_l$ is
    thin in $P$ then so are $(\phi,\psi)\cdot\face_{l-1}$ and
    $(\phi,\psi)\cdot\face_{l+1}$.
  \end{enumerate}
  The reader might like to compare these conditions to the
  corresponding clauses of definition~\ref{gen.horn}, as we do in
  detail in the proof of lemma~\ref{main.prod.lem} later on.

  Most importantly, if we are given entire supersets $N$ and $M$ of
  $\Delta[n]$ and $\Delta[m]$ respectively then each of the stratified
  sets $N\laxgray M$ and $N\gray M$ satisfies the condition required
  of $P$ in the last paragraph. Indeed, it is also the case that
  $\Delta[n]\laxgray\Delta[m]$ is minimal for those conditions, in the
  sense that it is an entire subset of any $P$ which satisfies
  them. The proofs of these facts are a matter of routine combinatorial
  verification, directly from the definitions of $\gray$ and
  $\laxgray$, which we leave to the reader. As a guide, sections~7
  and~8 of~\cite{Verity:2006:Complicial} contain numerous examples of
  detailed calculations involving the stratification of the lax Gray
  tensor product.

  We will also have reason to consider the following stratified subsets of
  $P$:
  \begin{itemize}
  \item $\boundary P$, the {\em boundary\/} of $P$, which is
    the regular subset whose underlying simplicial set is
    $\boundary(\Delta[n]\times\Delta[m])
    \defeq(\boundary\Delta[n]\times\Delta[m])\cup(\Delta[n]\times
    \boundary\Delta[m])\subseteq_s\Delta[m]\times\Delta[n]$.
  \item $P_d$ which is the regular subset generated by the set of
    shuffles in $P$ of depth less than or equal to $d\in[n+m]$.
  \item $\boundary P_d$, the boundary of $P_d$, which is the
    intersection of $P_d$ and $\boundary P$.
  \item $\tilde{P}_d$ and $\boundary\tilde{P}_d$ which are the entire
    subsets of $P$ defined by
    $\tilde{P}_d\defeq(\Delta[n]\laxgray \Delta[m])\cup P_d$ and
    $\boundary\tilde{P}_d\defeq(\Delta[n]\laxgray \Delta[m])\cup \boundary
    P_d$.
  \end{itemize}

  Before moving on, it is worth noting that by convention we take
  $P_{-1}$ to be empty and that $P_{nm}$ is equal to $P$
  itself. Furthermore, it is easily seen that $P_{nm-1}$ is precisely
  the regular subset of $P$ containing those simplices which do not
  have $(n,0)$ as a vertex. This latter set, its boundary and their
  associated unions with $\Delta[n]\laxgray\Delta[m]$ are our real
  objects of interest in the following lemma (and its corollary) and
  so we adopt the denotations $P\nc$, $\boundary P\nc$,
  $\tilde{P}\nc$ and $\boundary \tilde{P}\nc$ for these in
  order to avoid tedious repetition of the index $nm-1$.
\end{notation}


\begin{lemma}\label{main.prod.lem}
  For each integer $d\in\mathbb{N}$ with $0\leq d < nm$ the regular
  inclusion $\overinc\subseteq_r:P_{d-1}\cup\boundary P_d->P_d.$ (cf.\ 
  notation~\ref{P.notation}) is an inner anodyne extension. It follows
  that the regular inclusion $\overinc\subseteq_r:\boundary P\nc
  -> P\nc.$ is also an inner anodyne extension.
\end{lemma}

\begin{proof}
  This result depends on a few simple combinatorial observations:

  \begin{enumerate}[label=\textbf{observation (\roman*)},
                    ref={(\roman*)}, itemsep=1.5ex, fullwidth]
  \item\label{main.prod.lem.1} {\em If $(\alpha,\beta)$ is a shuffle
      in $P$ and $0<t<n+m$ is an integer with $\beta(t-1)=\beta(t)$
      and $\alpha(t)=\alpha(t+1)$ (right-lower corner) then
      $(\alpha,\beta)\cdot\face_t$ is a face of some shuffle of lower
      depth.}\vspace{1ex}

    Observe that the ordinate summation property of shuffles given in
    definition~\ref{shuffle.depth} may be applied to the conditions on
    $\alpha$ and $\beta$ at $t$ in the statement to establish that
    $\alpha(t)=\alpha(t-1)+1$ and $\beta(t+1)=\beta(t)+1$. It follows
    easily that we may construct a well defined $(n+m)$-simplex
    $(\alpha',\beta')$ by letting
    \begin{equation*}
      \alpha'(i)=
      \begin{cases}
        \alpha(i) & \text{if $i\neq t$} \\
        \alpha(t)-1 & \text{if $i=t$}
      \end{cases}
      \mkern20mu\text{and}\mkern20mu
      \beta'(i)=
      \begin{cases}
        \beta(i) & \text{if $i\neq t$} \\
        \beta(t)+1 & \text{if $i=t$}
      \end{cases}
    \end{equation*}
    which simply turns the original right-lower corner in
    $(\alpha,\beta)$ into a left-upper corner in
    $(\alpha',\beta')$. Now it is clear that this again satisfies the
    ordinate summation property, making it a shuffle which only
    differs from our original one at $t$ and thus has
    $(\alpha',\beta')\cdot\face_t=(\alpha,\beta)\cdot\face_t$. Furthermore,
    the expressions for the depths of these only differ at $t$ where
    $\min(\alpha'(t),m-\beta'(t)) = \min(\alpha(t),m-\beta(t)) - 1$ so
    it follows that $\depth(\alpha',\beta')=\depth(\alpha,\beta) - 1$
    as required.

  \item\label{main.prod.lem.2} {\em Suppose that $(\alpha,\beta)$ and
      $(\alpha',\beta')$ are distinct shuffles of the same depth $d$
      say then any face common to both of them is an element of
      $P_{d-1}$.}\vspace{1ex}

    Observe that the integer $$s\defeq\max\{i\in[n+m] \mid (\forall
    j\leq i) \alpha(j)=\alpha'(j) \}$$ is well defined, since all
    shuffles have $(0,0)$ as their $0\oth$ vertex, and that it is less
    than $n+m$, because our shuffles are distinct. Notice also that
    the ordinal summation property combined with the definition of $s$
    ensures that our shuffles agree at all vertices up to and
    including their $s\oth$ ones and that we may assume w.l.o.g. that
    $\alpha(s+1)=\alpha(s)+1$, $\beta(s+1)=\beta(s)$ and
    $\alpha'(s+1)=\alpha'(s)=\alpha(s)$,
    $\beta'(s+1)=\beta'(s)+1=\beta(s)+1$, by swapping the identities
    of our shuffles if necessary. Now we also know that the
    integer $$t\defeq\max\{i\in[n+m]\mid i> s \wedge
    \beta(i)=\beta(s)\}$$ is well defined, since the set we are taking
    this maximum over certainly contains $s+1$, and it must be less
    than $n+m$, because $\beta(s+1)<\beta'(s+1)\leq m$ whereas the
    ordinate summation property implies that $\beta(n+m)=m$. By
    construction our shuffles disagree at $t$, since
    $\beta(t)=\beta(s)=\beta'(s)<\beta'(s+1)\leq\beta'(t)$, so any
    face common to both of them must be a face of
    $(\alpha,\beta)\cdot\face_t$. Furthermore the maximality of $t$
    combined with the ordinal summation property implies that we also
    have $\beta(t-1)=\beta(t)$ and $\alpha(t)=\alpha(t+1)$
    (right-lower corner), so we may apply the previous observation to
    show that the simplex $(\alpha,\beta)\cdot\face_t$ is a face of a
    shuffle of depth less than $d=\depth(\alpha,\beta)$ and that it is
    thus an element of $P_{d-1}$. However, since the shuffles
    $(\alpha,\beta)$ and $(\alpha',\beta')$ disagree at their $t\oth$
    vertex it follows that any simplex $(\phi,\psi)$ common to both of
    them must of necessity be a face of $(\alpha,\beta)\cdot\face_t$
    and therefore must also be an element of $P_{d-1}$ as required.

  \item\label{main.prod.lem.3} {\em If $(\alpha,\beta)$ is a shuffle
      of depth $d$ then we may factor the corresponding Yoneda map
      $\arrow\yoneda{(\alpha,\beta)}:\Delta[n+m]-> P_d.$ as a
      composite of an entire inclusion
      $\overinc\subseteq_e:\Delta[n+m]->N.$ and a regular inclusion
      $\inc \yoneda{(\alpha,\beta)}:N->P_d.$. The entire superset $N$
      of $\Delta[n+m]$ is a $\vec{k}$-complicial $(n+m)$-simplex 
      with respect to the family 
      \begin{equation*}
        \vec{k} =\{k\in\mathbb{N}\mid 0<k<n+m \wedge
        \alpha(k-1)=\alpha(k) \wedge \beta(k)=\beta(k+1)\}
      \end{equation*}
      of (indices of) the left-upper corners of
      $(\alpha,\beta)$.}\vspace{1ex}

    Firstly it is clear that the simplicial set
    $\Delta[n]\times\Delta[m]$ enjoys the property that the Yoneda
    maps associated with its non-degenerate simplices are all
    simplicial inclusions. So applying the entire coimage
    factorisation of definition~\ref{reg.ent.def} to the Yoneda map
    associated with the shuffle $(\alpha,\beta)$ of $P_d$ we obtain
    the entire superset $N$ of $\Delta[n+m]$ and the regular inclusion
    $\overinc:N->P_d.$ of the statement. In other words, a simplex
    $\gamma$ is thin in $N$ if and only if its image
    $(\alpha,\beta)\cdot\gamma=(\alpha\circ\gamma,\beta\circ\gamma)$
    under the Yoneda map $\yoneda{(\alpha,\beta)}$ is thin in
    $P_d$. Now index the elements of $\vec{k}$ in increasing order
    $\{k_1<k_2<...<k_t\}$ and notice that
    observation~\ref{shuffle.facts} ensures that this satisfies the
    conditions of definition~\ref{gen.horn}, in particular its
    point~\ref{shuffle.facts.d} tells us that $\vec{k}$ is non-empty,
    because by assumption the depth $d$ of our shuffle is less than
    $mn$, whereas its point~\ref{shuffle.facts.e} implies that for
    each index $1\leq i<t$ we have $k_i+1<k_{i+1}$. Notice also that
    all of the elements of $\vec{k}$ are greater than $0$ and less
    than $n$, so all of our arguments here will involve inner
    (generalised) horns.
    
    Next let us examine the stratification of $N$ in more
    detail. First, suppose $\gamma$ is a $k$-admissible $r$-simplex of
    $\Delta[n+m]$ for some $k\in\vec{k}$ and let $l\in[r]$ be the
    unique integer with $\gamma(l)=k$, $\gamma(l-1)=k-1$ and
    $\gamma(l+1)=k+1$ (cf. observation~\ref{kadmis}). To see if
    $\gamma$ is thin in $N$ we consider the associated simplex
    $(\alpha,\beta)\cdot\gamma=(\alpha\circ\gamma,\beta\circ\gamma)$
    of $P_d$, for which the defining property of the elements
    $k\in\vec{k}$ and the definition of $l$ provides us with the
    equalities $\alpha\circ\gamma(l-1)=
    \alpha(k-1)=\alpha(k)=\alpha\circ\gamma(l)$ and
    $\beta\circ\gamma(l)=\beta(k)=\beta(k+1)=\beta\circ\gamma(l+1)$.
    These equalities demonstrate that $(\alpha,\beta)\cdot\gamma$ is
    thin in $P$, by condition~\ref{P.notation.a} of
    notation~\ref{P.notation}, and thus it is thin in the regular
    subset $P_d\subseteq_r P$.  Consequently $\gamma$ is thin in $N$
    and thus, quantifying over $k\in\vec{k}$ and all $k$-admissible
    simplices $\gamma\in\Delta[n+m]$, we have demonstrated that $N$
    satisfies condition~\ref{gen.horn.a} of
    definition~\ref{gen.horn}. However, we can take this argument a
    bit further and observe that if $\gamma\circ\face_l$ is thin in
    $N$ then $(\alpha,\beta)\cdot
    (\gamma\circ\face_l)=(\alpha\circ\gamma, \beta\circ\gamma)
    \cdot\face_l$ is thin in $P_d\subseteq_r P$, so we may apply the
    condition~\ref{P.notation.b} of notation~\ref{P.notation} to show
    that $(\alpha,\beta)\cdot
    (\gamma\circ\face_{l-1})=(\alpha\circ\gamma,
    \beta\circ\gamma)\cdot\face_{l-1}$ and $(\alpha,\beta)\cdot
    (\gamma\circ\face_{l+1})=(\alpha\circ\gamma,
    \beta\circ\gamma)\cdot\face_{l+1}$ are also both thin in $P_d$ and
    thus that $\gamma\circ\face_{l-1}$ and $\gamma\circ\face_{l+1}$
    are thin in $N$ as required by condition~\ref{gen.horn.b} of
    definition~\ref{gen.horn}. In other words, we have shown that $N$
    is a $\vec{k}$-complicial $(n+m)$-simplex as required.

  \item\label{main.prod.lem.4} {\em The generalised horn
      $\Lambda^{\vec{k}}N\subseteq_r N$ is the inverse image of the
      regular subset $P_{d-1}\cup\boundary P_d\subseteq_r P_d$ along
      the inclusion $\inc\yoneda{(\alpha,\beta)}:N->P_d.$.}\vspace{1ex}

    The inverse image of a regular subset along any stratified map is
    always a regular subset, so all we need do is check that the
    inverse image $L\subseteq_r N$ of $P_{d-1}\cup\boundary
    P_d\subseteq_r P_d$ along $\inc\yoneda{(\alpha,\beta)}:H->P_d.$
    contains the same simplices as $\Lambda^{\vec{k}}N$. So to show
    that $\Lambda^{\vec{k}}N\subseteq_r L$ we recall that
    $\Lambda^{\vec{k}}N$ is the regular subset generated by the set of
    faces $\{\face_i\mid i\in[n]\setminus\vec{k}\}$ and infer that it
    is enough to show that each of these is a simplex of $L$. However,
    consulting point~\ref{shuffle.facts.f} of
    observation~\ref{shuffle.facts} we see that if an integer
    $i\in[n]$ is not in $\vec{k}$, that is to say it is not the index
    of a left-upper corner, then it must either be the index of a
    right-lower corner, in which case we may apply
    observation~\ref{main.prod.lem.1} of this proof to show that
    $(\alpha,\beta)\cdot\face_i$ is in $P_{d-1}$, or it must be an
    element of the boundary $\boundary P_d$. In other words, for each
    $i\in[n]\setminus\vec{k}$ the face $(\alpha,\beta)\cdot\face_i$ is
    in $P_{d-1}\cup\boundary P_d$ and so $\face_i$ is in its inverse
    image $L$ as required.

    Conversely, we know that a simplex $\gamma\in N$ is not an element
    of $\Lambda^{\vec{k}} N$ if and only if each element of
    $[n]\setminus\vec{k}$ is also an element of $\im(\gamma)$, or
    equivalently iff any element of $[n]$ which is not in
    $\im(\gamma)$ is in $\vec{k}$. It then follows, from the
    definition of $\vec{k}$ with respect to $(\alpha,\beta)$, that
    each of the operators $\alpha\circ\gamma$ and $\beta\circ\gamma$
    are surjective and in particular we know that
    $\yoneda{(\alpha,\beta)}(\gamma)=
    (\alpha\circ\gamma,\beta\circ\gamma)$ is not an element of the
    boundary $\boundary P_d$. Furthermore, suppose that
    $(\alpha',\beta')$ is any other shuffle that has
    $(\alpha\circ\gamma,\beta\circ\gamma)$ as a face and observe that
    we then have $\alpha\circ\gamma=\alpha'\circ\gamma$ and
    $\beta\circ\gamma=\beta'\circ\gamma$, so if we consider an index
    $k$ at which $(\alpha,\beta)$ and $(\alpha',\beta')$ differ, it
    follows that it cannot be an element of $\im(\gamma)$ and thus
    that it must be an element of $\vec{k}$. Conversely, we know that
    if $k$ is in $\vec{k}$ then neither $k-1$ nor $k+1$ can be in
    there so it follows that they must both be in $\im(\gamma)$ and
    thus that $(\alpha,\beta)$ and $(\alpha',\beta')$ must agree at
    those indices. Now since $k$ is in $\vec{k}$ we know that
    $(\alpha,\beta)$ has a left-upper corner there and since our two
    shuffles agree at $k-1$ and $k+1$ but disagree at $k$ it is clear
    that $(\alpha',\beta')$ must have a right-lower corner at that
    index. It follows, therefore, that $\alpha'(k)=\alpha(k)+1$ and
    $\beta'(k)=\beta(k)-1$ and thus that $\min(\alpha'(k),m-\beta'(k))
    =\min(\alpha(k),m-\beta(k))+1$, so since this is true at any index
    where these shuffles differ it follows that
    $d=\depth(\alpha,\beta)<\depth(\alpha',\beta')$. This demonstrates
    that $\yoneda{(\alpha,\beta)}(\gamma)=(\alpha\circ\gamma,
    \beta\circ\gamma)$ is not an element of $P_{d-1}$, since we've
    shown that any shuffle of which it is a face must have depth at
    least $d$, so combining this with the corresponding fact with
    respect to the boundary $\boundary P_d$ we find that $\gamma$
    cannot be an element of the inverse image $L$ of
    $P_{d-1}\cup\boundary P_d$ along $\yoneda{(\alpha,\beta)}$ as
    required.
  \end{enumerate}\vspace{1ex}

  Now we may apply these observations to proving the result described
  in the first sentence of the statement. To do so enumerate the
  shuffles of depth $d$ over a suitable index set $I=\{1,2,...,s\}$
  and for each $i\in I$ let $N_i$ denote the entire superset of
  $\Delta[n+m]$ and $\vec{k}_i$ denote the family of integers
  associated with the $i\oth$ shuffle $(\alpha_i,\beta_i)$ as in
  observation~\ref{main.prod.lem.3}. Now define an increasing sequence
  of regular subsets $X_i$ of $P_d$ by starting at $X_0=P_{d-1}\cup
  \boundary P_d$ and letting each successive $X_i$ be the smallest
  regular subset of $P_d$ containing its predecessor $X_{i-1}$ and the
  shuffle $(\alpha_i,\beta_i)$, thus ensuring that the last member of
  this sequence $X_s$ is actually equal to $P_s$ itself. Now suppose
  that $\gamma$ is an $r$-simplex in $N_i$ and consider the face
  $(\alpha_i,\beta_i)\cdot\gamma\in P_d$ which is its image under the
  regular inclusion $\inc\yoneda{(\alpha_i,\beta_i)}:N_i->P_d.$ of
  observation~\ref{main.prod.lem.3}. If this is an element of
  $X_{i-1}$ then, by definition, it must either be an element of
  $P_{d-1}\cup\boundary P_d$ or it must also be a face of some other
  shuffle $(\alpha_j,\beta_j)$ with $j\leq i-1$, and in the latter
  case we may apply observation~\ref{main.prod.lem.2} to show that it
  is again an element of $P_{d-1}$. Consequently, the inverse image of
  $X_{i-1}\subseteq_r P_d$ along
  $\inc\yoneda{(\alpha_i,\beta_i)}:N_i->P_d.$ coincides with the
  inverse image its regular subset $P_{d-1}\cup\boundary P_d$ along
  the same map which we know, by observation~\ref{main.prod.lem.4}, is
  the generalised horn $\Lambda^{\vec{k}_i}N_i\subseteq_r
  N_i$. Furthermore, the definition of $X_i$ may be trivially recast
  to say that it is the union of $X_{i-1}$ and the direct image of
  $N_i$ under the regular inclusion
  $\inc\yoneda{(\alpha_i,\beta_i)}:N_i->P_d.$. Summarising these facts
  by applying recollection~\ref{glueing} we obtain a glueing square
  \begin{equation*}
    \xymatrix@R=2em@C=4em{
      *+!!<0em,0.4em>{\Lambda^{\vec{k}_i}N_i}\pbexcursion
      \ar@{u(->}[r]^-{\subseteq_r}
      \ar@{u(->}[d] & *+!!<0em,0.4em>{N_i}
      \ar@{u(->}[d]^{\yoneda{(\alpha_i,\beta_i)}}\\
      {X_{i-1}}\ar@{u(->}[r]_-{\subseteq_r} & 
      {X_i}\poexcursion}
  \end{equation*}
  which demonstrates that its lower horizontal inclusion is an inner
  anodyne extension since its upper horizontal is the inner anodyne
  extension of lemma~\ref{gen.horn.lem}. It follows, therefore, that
  $\overinc\subseteq_r:P_{d-1}\cup \boundary P_d-> P_d.$ may be
  decomposed as a composite of inner anodyne extensions
  $\overinc\subseteq_r:X_{i-1}->X_i.$ and is thus itself an inner
  anodyne extension as required.

  Finally, to prove the last sentence of the statement observe that we
  have equalities $P_d\cup(P_{d-1}\cup\boundary P\nc)=
  P_d\cup\boundary P\nc$ and $P_d\cap(P_{d-1}\cup\boundary
  P\nc)= P_{d-1}\cup\boundary P_d$ so we may apply
  recollection~\ref{glueing} to obtain a glueing square which displays
  the regular inclusion $\overinc\subseteq_r: P_{d-1}\cup\boundary
  P\nc-> P_d\cup\boundary P\nc.$ as a pushout of the anodyne
  extension $\overinc\subseteq_r:P_{d-1}\cup\boundary P_d->P_d.$. It
  follows therefore that we have a sequence of regular subsets
  $P_d\cup\boundary P\nc$ of $P\nc$ ($d=-1,0,1,...,nm-1$)
  whose first member is $\boundary P\nc$, whose last is $P\nc$
  and in which each inclusion of a sequence member into its successor
  was shown to be an anodyne extension (as a pushout of such) in the
  last sentence. It follows therefore that their composite
  $\overinc\subseteq_r: \boundary P\nc->P\nc.$ is also an
  anodyne extension as required.
\end{proof}

\begin{cor}\label{main.prod.cor}
  For each integer $d\in\mathbb{N}$ with $0\leq d< nm$ the entire
  inclusion $\overinc \subseteq_e: \tilde{P}_{d-1} \cup\boundary
  P_d->\tilde{P}_d.$ (cf.\ notation~\ref{P.notation}) is an inner
  anodyne extension. It follows that the entire inclusion
  $\overinc\subseteq_e:\boundary \tilde{P}\nc ->
  \tilde{P}\nc.$ is also an inner anodyne extension.
\end{cor}

\begin{proof}
  A routine reprise of the method used in the proof of the last lemma,
  replacing pushouts of generalised horn extensions by pushouts of the
  generalised thinness extensions of corollary~\ref{gen.horn.cor}
  wherever necessary. We leave the details to the reader.
\end{proof}

In the next lemma we use the notation $\anytens$ and $\canytens$ to
represent either one of the tensors $\gray$ or $\laxgray$ on $\Strat$
and its associated corner tensor.

\begin{lemma}\label{left.anod.prod.lem}
  If $k$ is an integer with $0\leq k<n$ then each of the corner
  tensors
  \begin{gather*}
    (\overinc\subseteq_r:\Lambda^k[n]->\Delta^k[n].)\canytens
    (\overinc\subseteq_r:\boundary\Delta[m]->\Delta[m].) \\
    (\overinc\subseteq_r:\Lambda^k[n]->\Delta^k[n].)\canytens
    (\overinc\subseteq_e:\Delta[m]->\Delta[m]_t.) \\
    (\overinc\subseteq_r:\Delta^k[n]'->\Delta^k[n]''.)\canytens
    (\overinc\subseteq_r:\boundary\Delta[m]->\Delta[m].) \\
    (\overinc\subseteq_r:\Delta^k[n]'->\Delta^k[n]''.)\canytens
    (\overinc\subseteq_e:\Delta[m]->\Delta[m]_t.)
  \end{gather*}
  is a left anodyne extension and is an inner anodyne extension if
  $0<k$.
\end{lemma}

\begin{proof}
  We prove the stated result for the first corner tensor in the list
  above in detail. Firstly arguing just as we did in
  observation~\ref{corner.join} we see that the corner tensor of the
  two maps in the statement is (isomorphic to) the regular subset
  inclusion $\overinc\subseteq_r:(\Lambda^k[n]
  \anytens\Delta[m])\cup(\Delta^k[n]\anytens\boundary\Delta[m])->
  \Delta^k[n]\anytens\Delta[m].$ and we adopt the letter $R$ to denote
  the codomain of this inclusion. To prove that this is a left anodyne
  extension we will apply lemma~\ref{main.prod.lem} {\em twice\/} to
  the stratified sets $Q\defeq\Delta[n-1]\anytens \Delta[m]$ and
  $P\defeq\Delta^k[n]\anytens \Delta[m]$ respectively. So consider the
  increasing sequence $R \subseteq_r R\cup \boundary P\nc \subseteq_r
  R\cup P\nc \subseteq_r P$ of regular subset inclusions, which are
  subject to the following observations:
  \begin{enumerate}[
    label=\textbf{observation (\roman*)},
    fullwidth, leftmargin=1em, itemsep=1ex]
  \item A simplex $(\alpha,\beta)$ is in $R$ iff either
    $\im(\alpha)\cup\{k\}\neq[n]$ or $\im(\beta)\neq[m]$ whereas it is
    in $\boundary P\nc$ iff it doesn't have $(n,0)$ as a simplex and
    either $\im(\alpha)\neq[n]$ or $\im(\beta)\neq[m]$. So under the
    assumption that $k<n$ we define $W$ to be the regular subset of
    $P$ of those simplices $(\alpha,\beta)$ which don't have $(n,0)$
    as a simplex and for which $k\notin\im(\alpha)$ and may then
    easily demonstrate that $R\cup\boundary P\nc=R\cup W$.

  \item The stratified map corresponding by Yoneda's lemma to the
    $(n-1)$-simplex $\face_k$ in $\Delta^k[n]$ is a regular inclusion
    $\inc\yoneda{\face_k}:\Delta[n-1]->\Delta^k[n].$. Furthermore,
    each of the tensors $\gray$ and $\laxgray$ preserves inclusions
    and regularity, as the reader may readily verify, so it follows
    that we obtain a regular inclusion
    $\yoneda{\face_k}\anytens\Delta^k[n]$ from
    $Q=\Delta[n-1]\anytens\Delta[m]$ to $P$. Under the assumption that
    $k<n$ it is then easily seen that the regular subset $W$ of the
    last observation is simply the direct image of $Q\nc\subseteq_r Q$
    under this inclusion and that $\boundary Q\nc$ is the inverse
    image of $R\subseteq_r P$ along the same inclusion. It follows, by
    applying recollection~\ref{glueing}, that we have a glueing square
    \begin{equation*}
      \xymatrix@=2em{
        {\boundary Q\nc}\pbexcursion
        \ar@{u(->}[r]^-{\subseteq_r}\ar@{u(->}[d] &
        {Q\nc}\ar@{u(->}[d] \\
        {R}\ar@{u(->}[r]_-{\subseteq_r} &
        {R\cup W}\save []+R*!L{{}=R\cup\boundary P\nc}\restore
        \poexcursion
      }
    \end{equation*}
    whose upper horizontal is an anodyne extension by
    lemma~\ref{main.prod.lem}. Consequently its pushout, the inclusion
    of $R$ into $R\cup\boundary P$, is also an inner anodyne
    extension.

  \item Clearly we have $(R\cup\boundary P\nc)\cup P\nc=R\cup P\nc$
    and $(R\cup\boundary P\nc)\cap P\nc = R\cup\boundary P\nc$, where
    the latter equality holds because $R$ is a subset of the boundary
    $\boundary P$, so we may apply recollection~\ref{glueing} to
    obtain a glueing square:
    \begin{equation*}
      \xymatrix@=2em{
        {\boundary P\nc}\pbexcursion
        \ar@{u(->}[r]^-{\subseteq_r}\ar@{u(->}[d]_{\subseteq_r} &
        {P\nc}\ar@{u(->}[d]^{\subseteq_r} \\
        {R\cup\boundary P\nc}\ar@{u(->}[r]_-{\subseteq_r} &
        {R\cup P\nc}\poexcursion
      }
    \end{equation*}
    Again we may apply lemma~\ref{main.prod.lem} to show that the
    upper horizontal here and its pushout, the inclusion
    $\overinc\subseteq_r: R\cup\boundary P\nc->R\cup P\nc.$, are both
    inner anodyne extensions.

  \item Only two simplices of $P$ are not elements of $R\cup P\nc$,
    those being the maximal depth shuffle
    $(\partproj^{n,m}_1,\partproj^{n,m}_2)$ and its $k\oth$
    $(n+m-1)$-dimensional face
    $(\partproj^{n,m}_1,\partproj^{n,m}_2)\cdot\face_k$. Consequently
    we are led to considering the Yoneda map corresponding to this
    shuffle, which is an inclusion
    $\inc\yoneda{(\partproj^{n,m}_1, \partproj^{n,m}_2)}: \Delta[n+m]
    -> P.$. In fact this may be lifted to a stratified map whose
    domain is $\Delta^k[n+m]$, although the combinatorial details of
    the argument demonstrating that fact (which we leave to the
    reader) depend upon precisely which of the tensors $\gray$ or
    $\laxgray$ we are studying. Indeed these cases diverge a little
    further at this point, since it turns out that the face
    $(\partproj^{n,m}_1,\partproj^{n,m}_2)\cdot\face_k$ is thin in $P$
    when it is defined using $\gray$ but that this simplex is not thin
    in there when we consider $\laxgray$.  In that first case, it
    turns out that the flanking faces
    $(\partproj^{n,m}_1,\partproj^{n,m}_2)\cdot\face_i$
    ($i\in\{k-1,k+1\}\cap[n+m]$) are also thin in $P$ so we may lift
    our map further to one with domain $\Delta^k[n+m]''$. Ultimately
    however, regardless of tensor, we may apply
    recollection~\ref{glueing} and obtain one of the following glueing
    squares:
    \begin{equation*}
      \xymatrix@=2em{
        {\Lambda^k[n+m]'}\pbexcursion\ar@{u(->}[r]^{\subseteq_r}\ar@{u(->}[d] & 
        {\Delta^k[n+m]''}
        \ar@{u(->}[d]^{\yoneda{(\partproj^{n,m}_1,\partproj^{n,m}_2)}} \\
        {R\cup P\nc}\ar@{u(->}[r]_{\subseteq_r} &
        {\Delta^k[n]\gray\Delta[m]}\poexcursion }
      \mkern20mu
      \xymatrix@=2em{
        {\Lambda^k[n+m]}\pbexcursion\ar@{u(->}[r]^{\subseteq_r}\ar@{u(->}[d] & 
        {\Delta^k[n+m]}
        \ar@{u(->}[d]^{\yoneda{(\partproj^{n,m}_1,\partproj^{n,m}_2)}} \\
        {R\cup P\nc}\ar@{u(->}[r]_{\subseteq_r} &
        {\Delta^k[n]\laxgray\Delta[m]}\poexcursion }
    \end{equation*}
    Consequently the inclusion $\overinc\subseteq_r:R\cup P\nc-> P.$
    is a pushout of a left outer (and possibly thin) horn if $k=0$ and
    of an inner (and possibly thin) horn otherwise and is thus a left or
    inner anodyne extension.
  \end{enumerate}
  Summarising these observations, we see that each of the first two
  inclusions in our sequence above are inner anodyne extensions
  whereas the last one is a left anodyne extension if $k=0$ and an
  inner anodyne extension if $0<k<n$. It follows therefore that their
  composite, the corner join under study, is an inner or left
  anodyne extension as described in the statement.

  The remaining three corner tensors of the statement, each of which
  is an entire inclusion, may all be shown to be left or inner anodyne
  extensions using a routine reprise of the argument above. The
  primary modification required is that we replace pushouts of
  the inner anodyne extension of lemma~\ref{main.prod.lem} by pushouts
  of the corresponding one of corollary~\ref{main.prod.cor}. We leave
  the details to the reader.
\end{proof}

\begin{obs}
  Notice that the assumption $k<n$ was vital to all of the
  observations made in the proof above. A completely
  different combinatorial argument would be required to directly prove
  the corresponding result for right outer horns. This however need
  not bother us here, since everything we need with regard to right
  compliciality may be derived from lemma~\ref{left+inner=outer.cor}
  as we do in theorem~\ref{comp.hom.stab.thm} below.
\end{obs}


\begin{obs}\label{anod.ctens.obs}
  The Gray tensor product $\gray$ preserves colimits in each variable,
  so we may apply observation~\ref{corner.cofibration} to the result
  of the last lemma and show that if $\inc e:U->V.$ is a left (resp.\
  inner) anodyne extension and $\inc i:X->Y.$ is any inclusion then
  their corner join $e\cgray i$ is also a left (resp.\ inner) anodyne
  extension.

  Things are, however, somewhat less straightforward for the lax Gray
  tensor product $\laxgray$ which is not well behaved with respect to
  colimits. Unfortunately it is also not possible to replace this by
  the related ``colimit friendly'' pre-tensor $\pretens$ because this
  does not satisfy the conditions required to make the arguments of
  the last few lemmas work.  We will return to resolve this issue in
  section~\ref{quillen.model.sec}, for now however we have the
  following useful theorem for the closed structure associated with
  $\gray$:
\end{obs}



\begin{thm}\label{comp.hom.stab.thm}
  If $A$ is a weak complicial set and $X$ is any stratified set then
  the stratified set of strong transformations $\hom(X,A)$ is also a
  weak complicial set. Furthermore if $\arrow p:A->B.$ is a complicial
  fibration between complicial sets and $\overinc i:X->Y.$ is any
  inclusion then the corner closure $\chom(i,p)$ is also a complicial
  fibration.
\end{thm}

\begin{proof}
  Applying observation~\ref{corner.fibration} to the result of the
  last observation we find that $\hom(X,A)$ is a weak left complicial
  set to which we may apply corollary~\ref{left+inner=outer.cor} to
  demonstrate that it is actually a weak complicial set.  A similar
  argument shows that the corner closure $\chom(i,p)$ is a left
  complicial fibration, whose codomain is the pullback in the
  following diagram:
  \begin{equation*}
    \xymatrix@R=2em@C=4em{
      {\hom(Y,A)}\ar[r]^-{\chom(i,p)} &
      {\hom(Y,B)\times_{\hom(X,B)}\hom(X,A)}\pbexcursion
      \ar[r]\ar[d] & {\hom(X,A)}\ar[d]^{\hom(X,p)} \\
      & {\hom(Y,B)}\ar[r]_{\hom(i,B)} & {\hom(X,B)}
    }
  \end{equation*}
  Of course, the right hand vertical of this square is a left
  complicial fibration by the result of the last sentence, since it is
  the right corner closure of $p$ and the inclusion
  $\overinc:\emptyset-> X.$. It follows that its left hand vertical is
  also a left complicial fibration, since these are stable under
  pullback, whose codomain $\hom(Y,B)$ is a weak left complicial set
  as already discussed. Consequently its domain, our pullback, is also
  a weak left complicial set so we may apply
  corollary~\ref{left+inner=outer.cor} to show that it is actually a
  weak complicial set which then enables us to apply the same
  corollary again to show that $\chom(i,p)$ is a complicial fibration
  as required.
\end{proof}

\begin{cor}\label{anod.tensor.stab.cor}
  If $\inc e:U->V.$ is an anodyne extension and $\inc i:X->Y.$ is any
  inclusion then their corner tensor $e\cgray i$ has the LLP with
  respect to all complicial sets and all complicial fibrations between
  complicial sets.
\end{cor}

\begin{proof}
  Apply observation~\ref{corner.fibration} to show that this result is
  simply dual to that of the last theorem under the adjunction
  $-\cgray i\dashv\chom(i,*)$.
\end{proof}

\begin{obs}[the Gray-category of weak complicial sets]
  \label{wcs.gray.obs}
  We may canonically enrich $\Strat$ with respect to its Gray tensor
  $\gray$, to obtain an enriched category by taking $\hom(X,Y)$ as its
  stratified homset between the stratified sets $X$ and $Y$ (cf.\
  Kelly~\cite{Kelly:1982:ECT} for the
  details). Theorem~\ref{comp.hom.stab.thm} now tells us that the
  homsets of its enriched full subcategory $\Wcs$ of weak complicial
  sets are all themselves weak complicial sets. We call such gadgets
  {\em Gray-categories\/} and the reader may find out much more about
  these structures by consulting the companion
  paper~\cite{Verity:2005:GrayNerves}.
\end{obs}

\subsection{A Characterisation of Strict Complicial
  Sets}\label{str.comp.char}

Before moving on, these results allow us to establish another important
characterisation of strict complicial sets amongst the weak
ones. 

\begin{defn}[reprise of definition~117 of~\cite{Verity:2006:Complicial}]
  If $X$ is a stratified set, we say that an $n$-simplex $x\in X$ is
  {\em pre-degenerate at $k$} if its face $x\cdot\alpha$ is thin
  whenever $\arrow \alpha:[m]->[n].$ is a simplicial operator whose
  image contains the vertices $k,k+1\in[n]$. Most importantly, the
  degeneracy condition on stratifications ensures that if $x$ is
  degenerate at $k$ then it is pre-degenerate at $k$.

  Conversely, we say that $X$ is {\em well-tempered\/} if whenever
  $x\in X$ is pre-degenerate at $k$ then it is actually degenerate at
  $k$. The slogan here is that in a well-tempered stratified set,
  thinness is a sufficient property for the detection of degeneracy.
\end{defn}

\begin{lemma}\label{welltemp=>trivhty}
  If $Y$ is a well-tempered stratified set then every stratified map
  $\arrow h:X\gray\Delta[1]_t->Y.$ factors through the projection map
  $\arrow\pi_X:X\gray\Delta[1]_t->X.$.
\end{lemma}

\begin{proof}(essentially that of corollary~164
  of~\cite{Verity:2006:Complicial}) It is clear that if $h$ may be
  factored through $\pi_X$ to give a stratified map $\arrow
  \hat{h}:X->Y.$ then this must be given by $\hat{h}(x)\defeq h(x,0)$,
  in which expression the symbol $0$ represents the instance of the
  constant operator of definition~\ref{simp.of.1} whose dimension
  matches that of $x$.  To show that this does indeed provide us with
  an appropriate factor we need to demonstrate that
  $h(x,\rho_{k+1})=\bar{h}(x)=h(x,0)$ for each one of the operators
  $\rho_{k+1}$ of definition~\ref{simp.of.1}.

  For definiteness let $r$ be the dimension of $x\in X$ and we'll
  decorate the operators of definition~\ref{simp.of.1} with their
  superscripted dimension. Consider the $(r+1)$-simplex
  $(x\cdot\degen_k,\rho^{r+1}_{k+1})$ of $X\gray\Delta[1]_t$ and
  observe that it is pre-degenerate at $k$ if and only if each one of
  its ordinates is pre-degenerate at $k$. However $x\cdot\degen_k$ is
  certainly pre-degenerate at $k$, since it is degenerate there, and
  every simplex of dimension greater than $0$ is thin in $\Delta[1]_t$
  which clearly implies that every one of its $(r+1)$-simplices is
  pre-degenerate at $k$. Now stratified maps clearly preserve
  pre-degeneracy so it follows that
  $h(x\cdot\degen_k,\rho^{r+1}_{k+1})$ is pre-degenerate at $k$ in $Y$
  and so there exists an $r$-simplex $y\in Y$ with
  $h(x\cdot\degen_k,\rho^{r+1}_{k+1})=y\cdot\degen_k$ since $Y$ is
  well-tempered. Observe now that we have the following calculations
  \begin{align*}
    y = (y\cdot\degen_k)\cdot\face_k & {}=
    h(x\cdot\degen_k,\rho^{r+1}_{k+1})\cdot\face_k =
    h((x\cdot\degen_k)\cdot\face_k,\rho^{r+1}_{k+1}\circ\face_k)=
    h(x,\rho^r_k) \\
    y = (y\cdot\degen_k)\cdot\face_{k+1} & {}=
    h(x\cdot\degen_k,\rho^{r+1}_{k+1})\cdot\face_{k+1} \\
    & {}=
    h((x\cdot\degen_k)\cdot\face_{k+1},\rho^{r+1}_{k+1}\circ\face_{k+1})=
    h(x,\rho^r_{k+1})
  \end{align*}
  wherein we rely repeatedly on the simplicial identities
  $\degen_k\circ\face_k=\id=\degen_k\circ\face_{k+1}$ and the
  easy observations that $\rho^{r+1}_{k+1}\circ\face_k=\rho^r_k$ and
  $\rho^{r+1}_{k+1}\circ\face_{k+1}=\rho^r_{k+1}$. In other words we
  have shown that for each $k=0,1,...,r$ we have
  $h(x,\rho^r_k)=h(x,\rho^r_{k+1})$ and composing these equalities we
  find that $h(x,\rho^r_k)=h(x,\rho^r_{r+1})=h(x,0)$ as required 
  where the last equality simply expresses the fact that the operators
  $\rho^r_{r+1}$ and $0$ are identical.
\end{proof}

\begin{thm}\label{wcs+welltemp=>scs}
  A stratified set $A$ is a (strict) complicial set if and only if it
  is a weak complicial set and it is well-tempered.
\end{thm}

\begin{proof}
  The ``only if'' part follows from the argument of
  example~\ref{comp.sets} and lemma~163
  of~\cite{Verity:2006:Complicial}.

  To prove the converse, first observe that it is enough to show that
  if $A$ is well-tempered and a weak complicial set then it has unique
  fillers for inner complicial horns. So suppose that we have a
  stratified map $\arrow f:\Lambda^k[n]->A.$ and with a pair of
  extensions $\arrow k_0,k_1:\Delta^k[n]->A.$ along the inclusion
  $\overinc\subseteq_r:\Lambda^k[n]->\Delta^k[n].$. From this
  information build a stratified map $\arrow
  h:(\Lambda^k[n]\gray\Delta[1]_t)\cup(\Delta^k[n]\gray
  \boundary\Delta[1])->A.$ by letting $h(\alpha,\beta)\defeq
  f(\alpha)$ on $\Lambda^k[n]\gray\Delta[1]_t$ and letting
  $h(\alpha,0) \defeq k_0(\alpha)$ and $h(\alpha,1) \defeq
  k_1(\alpha)$ on $\Delta^k[n]\gray\boundary\Delta[1]$, where $0$ and
  $1$ denote appropriate instances of the constant operators given in
  definition~\ref{simp.of.1}. Of course each of these pieces of $h$ is
  stratified and they match where mutually defined, because $k_0$ and
  $k_1$ both extend $f$, thus demonstrating that it is a well defined
  stratified map. Furthermore, it may be extended to a stratified map
  $\arrow \bar{h}:\Delta^k[n]\gray\Delta[1]_t->A.$, because $A$ is a
  weak complicial set and the inclusion of the domain of $h$ into
  $\Delta^k[n]\gray\Delta[1]_t$ is the corner tensor of the inner horn
  inclusion $\overinc\subseteq_r:\Lambda^k[n]-> \Delta^k[n].$ and the
  inclusion $\overinc\subseteq_r:\boundary\Delta[1]->\Delta[1]_t.$
  which is an (inner) anodyne extension by
  observation~\ref{anod.ctens.obs}. Since $A$ is well-tempered we may
  now apply lemma~\ref{welltemp=>trivhty} and factor $\bar{h}$ through
  the projection $\arrow\pi_{\Delta^k[n]}:\Delta^k[n]\gray\Delta[1]_t
  -> \Delta^k[n].$ to give a map $\arrow \hat{h}:\Delta^k[n]->A.$ and
  now we find that
  $k_0(\alpha)=f(\alpha,0)=\bar{f}(\alpha,0)=\hat{f}(\alpha)$ and
  $k_1(\alpha)=f(\alpha,1)=\bar{f}(\alpha,1)=\hat{f}(\alpha)$, which
  demonstrates that $k_0=k_1$ and thus that $f$ has exactly one
  extension as required.
\end{proof}



\section{Quillen Model Structures on Stratified Sets}
\label{quillen.model.sec}

In this section we muster the machinery developed in the last few
sections to demonstrate that the category of stratified sets $\Strat$
supports a natural Quillen model structure whose fibrant objects are
precisely the weak complicial sets. We do so using Jeffery Smith's
theorem for locally presentable categories, the conditions of which
we've recounted as theorem~\ref{jsmith.thm} in the appendix. 
As discussed in observation~\ref{LFP.quasitopos}, the category
$\Strat$ is locally finitely presentable and thus provides a context
within which to apply this theorem.

\begin{defn}\label{def.triv.fib}
  We define $I$ to be the set of boundary and thin simplex inclusions
  \begin{equation*}
    \{\overinc\subseteq_r:\boundary\Delta[n]->\Delta[n].\mid
    n=0,1,...\}\cup\{\overinc\subseteq_e:\Delta[n]->\Delta[n]_t.\mid 
    n=1,2,...\}
  \end{equation*}
  whose cellular completion $\cell(I)$ is the class of all inclusions
  of stratified sets (cf.\
  observation~\ref{strat.inc.obs}). Consequently, the members of the
  corresponding class of fibrations $\fib(I)$, called {\em trivial
    fibrations}, all enjoy the RLP with respect to arbitrary
  inclusions of stratified sets. Notice that it is immediate that all
  trivial fibrations are also complicial fibrations.
\end{defn}

\subsection{Homotopy Equivalences of Weak Complicial Sets}

The next few definitions and results are appropriated, with
appropriate modifications, from classical simplicial homotopy
theory.

\begin{defn}
  If $\arrow f,g:X->Y.$ are stratified maps then a {\em simple
    homotopy\/} from $f$ to $g$ is a stratified map $\arrow h:
  X\gray\Delta[1]_t->Y.$ for which $h(x,0)=f(x)$ and $h(x,1)=g(x)$ for
  all $x\in X$. Notice that in order to make sense of these
  expressions we assume that $0$ and $1$ denote suitable instances of
  the constant operators introduced in notation~\ref{simp.of.1}.  

  We write $f\sim_1 g$ if there exists a simple homotopy from $f$ to
  $g$ and let $f\sim g$, the {\em homotopy\/} relation, denote the
  transitive closure of that relation. 
\end{defn}

\begin{obs}
  Taking duals under the adjunction $X\gray-\dashv \hom(X,*)$ and
  appealing to Yoneda's lemma we see that a simple homotopy
  corresponds to a thin 1-simplex $\hat{h}$ in the stratified set
  $\hom(X,Y)$ whose vertices are $f=\hat{h}\cdot\vertex_0$ and
  $g=\hat{h}\cdot\vertex_1$. This presentation immediately tells us
  that the simple homotopy relation $\sim_1$ is already transitive
  (and is thus identical to $\sim$) whenever the codomain of our maps
  is a weak complicial set $A$.  To verify this fact simply observe
  that, by theorem~\ref{comp.hom.stab.thm}, $\hom(X,A)$ is a weak
  complicial set whenever $A$ is and demonstrate transitivity of
  simple homotopy using fillers for suitable (thin) $1$-dimensional
  horns in $\hom(X,A)$ to compose the witnessing simple homotopies.
\end{obs}

\begin{defn}\label{hty.equiv.def}
  If $X$ and $Y$ are stratified sets then a {\em homotopy
    equivalence\/} between them is a stratified map $\arrow e:X->Y.$
  which has an {\em equivalence inverse\/} $\arrow e':Y->X.$ for which we
  have $e'\circ e\sim\id_X$ and $e\circ e'\sim\id_Y$.
\end{defn}

\begin{obs}\label{equiv.basic.obs}
  The homotopy relation is preserved by pre-composition and
  post-composition in $\Strat$, so we may form a {\em homotopy
    category\/} $\Pi(\Strat)$ by taking the quotient of each of the
  homsets of $\Strat$ under the homotopy relation. Then a stratified
  map $\arrow e:X->Y.$ is a homotopy equivalence if and only if its
  corresponding homotopy class $\arrow[f]_{\sim}:X->Y.$ is an
  isomorphism in $\Pi(\Strat)$. 

  Using this observation we may immediately derive many useful
  properties of homotopy equivalences directly from the corresponding
  facts about isomorphisms in any category. In particular, in the
  sequel we will make use of the following very simple observations:
  \begin{itemize}
  \item {\bf homotopy stability} If $e$ is a homotopy equivalence then
    so is any stratified map homotopic to $e$.
  \item {\bf 2-of-3 property} If two of the stratified maps $e$, $f$
    and their composite $f\circ e$ are homotopy equivalences then so
    is the third.
  \item {\bf stability under retract} Retracts of homotopy
    equivalences are again homotopy equivalences.
  \item {\bf left inverse property} If $\arrow e:X->Y.$ is a homotopy
    equivalence and $\arrow\bar{e}:Y->X.$ is such that $\bar{e}\circ
    e\sim\id_X$ (left equivalence inverse) then we also have
    $e\circ\bar{e}\sim\id_Y$ (right equivalence inverse).
  \end{itemize}
\end{obs}

\begin{obs}\label{hty.eqv.coh}
  The Gray-category $\Wcs$ of weak complicial sets, which we
  introduced in observation~\ref{wcs.gray.obs}, gives rise to a Kan
  complex enriched category by applying the Gray tensor product
  preserving $0$-superstructure functor $\sst_0$ to its homsets. We
  may apply Cordier and Porter's homotopy coherent nerve
  functor~\cite{Cordier:1986:HtyCoh} to this structure to obtain a
  quasi-category $\nerv_{hc}(\Wcs)$. A presentation of this nerve
  construction suited to our needs here is provided in the companion
  paper~\cite{Verity:2005:GrayNerves}, which generalises the classical
  homotopy coherent nerve construction to provide a faithful embedding
  of the category of Gray-categories into the category of weak
  complicial sets.

  Now definition~\ref{hty.equiv.def} above may simply be regarded as
  saying that the stratified map $\arrow e:A->B.$ has an equivalence
  inverse in $\nerv_{hc}(\Wcs)$ in the sense of
  theorem~\ref{almost.weakinner.cor}. So applying that theorem, it
  follows that any homotopy equivalence gives rise to a simplicial map
  $\arrow:E^-_3->\nerv_{hc}(\Wcs).$ which may be pictured as:
  \begin{equation*}
      \begin{xy}
      \POS(0,0) *[o]{\xybox{\xymatrix@R=2.5em@C=1.5em{
          & {B}\ar[rr]^{e'} && {A}\ar[rd]^{e} & \\
          {A}\ar[ur]^{e}\ar[urrr]|(0.4){=}^{}="two" \ar[rrrr]_{e}^{}="one"
          &&&& {B} \ar@{}"one";"1,4"|(0.4){=}
          \ar@{}"1,2";"two"|{\objectstyle\stackrel{h}{\simeq}} }}}="a"
      \POS(55,0) *[o]{\xybox{\xymatrix@R=2.5em@C=1.5em{ &
          {B}\ar[rr]^{e'}\ar[drrr]|(0.6){=}^{}="one"
          && {A}\ar[rd]^{e} & \\
          {A}\ar[ur]^{e}="two"\ar[rrrr]_{e}^{}="two" &&&& {B}
          \ar@{}"two";"1,2"|(0.4){=}
          \ar@{}"one";"1,4"|{\objectstyle\stackrel{k}{\simeq}} }}}="b"
     \ar@{->}"a";"b"^t_{\simeq}
    \end{xy}
  \end{equation*}
  Unwinding the definition of $\nerv_{hc}(\Wcs)$ given
  in~\cite{Verity:2005:GrayNerves} it is easily seen that this data
  amounts to a choice of simple homotopies
  \begin{align*}
    &\arrow h:A\gray\Delta[1]_t->A.\text{ with $h(a,0)=e'(
      e(a))$ and $h(a,1)=a$} \\
    &\arrow k:B\gray\Delta[1]_t->B.\text{ with $k(b,0)=e(
      e'(b))$ and $k(b,1)=b$} 
  \end{align*}
  and a ``double'' homotopy
  \begin{align*}
    \arrow t:A\gray\Delta[1]_t\gray\Delta[1]_t->B.&
    \text{ with $t(a,0,\beta)=e(h(a,\beta))$ and
      $t(a,\alpha,0)=k(e(a),\alpha)$}\\
    & \text{ and $t(a,\alpha,1)=t(a,1,\beta)=e(a)$}
  \end{align*}
  connecting them.
\end{obs}

\begin{defn}\label{def.retr.defn}
  We say that a stratified map $\arrow e:X->Y.$ is a {\em (simple)
    deformation retraction\/} if there is a stratified map $\arrow
  m:X->Y.$ with $e\circ m=\id_Y$ and a simple homotopy $\arrow
  d:A\gray\Delta[1]_t->A.$ from $m\circ e$ to $\id_A$ with
  $e(d(a,\alpha))=e(a)$ (for all $a\in A$ and $\alpha\in\Delta[1]_t$).
\end{defn}

\begin{lemma}\label{cfib.equiv.lem}
  If $B$ is a weak complicial set and $\arrow e:A->B.$ is a
  complicial fibration then the following are equivalent:
  \begin{enumerate}[label=(\roman*)]
  \item\label{he.1} $e$ is a homotopy equivalence,
  \item\label{he.2} $e$ is a deformation retraction, and
  \item\label{he.3} $e$ is a trivial fibration.
  \end{enumerate}
\end{lemma}

\begin{obs}\label{comp.fib.comp}
  Notice that we need not assume explicitly that $A$ is a weak
  complicial set because we may immediately infer that this is the
  case from the compliciality assumptions on $B$ and $e$ and the fact
  that the class of complicial fibrations is closed under
  composition. In future in these cases we will simply say that such a
  map is a {\em complicial fibration of weak complicial sets}.
\end{obs}

\begin{proof} (of lemma~\ref{cfib.equiv.lem})
  This is fundamentally a classical result. Clearly a deformation
  retraction is a special sort of homotopy equivalence, so the
  implication $\text{\ref{he.2}}\Rightarrow\text{\ref{he.1}}$ is
  trivial.  The implication
  $\text{\ref{he.3}}\Rightarrow\text{\ref{he.2}}$ is also routine, we
  simply use the trivial fibration assumption on $e$ to make
  successive lifts
  \begin{equation*}
    \xymatrix@=2em{
      {\emptyset}\ar[r]\ar@{u(->}[d] &
      {A}\ar[d]^e \\
      {B}\ar[r]_{\id_B}\ar@{..>}[ur]|{\exists m} & {B} 
    }\mkern50mu
    \xymatrix@R=2em@C=4em{
      {A + A}\ar[r]^-{\langle m\circ e,\id_A\rangle}
      \ar@{u(->}[d] & {A}\ar[d]^e \\
      {A\gray\Delta[1]_t}\ar[r]_-{e\circ\pi_A}\ar@{..>}[ur]|{\exists d} &
      {B}}
  \end{equation*}
  to construct $m$ and $d$ satisfying the properties required by
  the last definition.

  To prove the reverse implication $\text{\ref{he.1}}\Rightarrow
  \text{\ref{he.2}}$, assume that $e$ is an equivalence with inverse
  $e'$ and that we are given the simple homotopies $h$ and $k$ and the
  double homotopy $t$ described in observation~\ref{hty.eqv.coh}. Now
  consider the squares
  \begin{equation}\label{def.lift.1}
    \xymatrix@R=2em@C=4em{
      {B\gray\Lambda^0[1]}\ar[r]^-{e'\circ\pi_B}\ar@{u(->}[d]_{\subseteq_r} &
      {A}\ar[d]^e \\
      {B\gray\Delta[1]_t}\ar[r]_-{k}\ar@{..>}[ur]|{\exists \bar{k}} & {B} 
    }\mkern50mu
    \xymatrix@R=2em@C=4em{
      {X}\ar[r]^-{f}
      \ar@{u(->}[d]_{\subseteq_r} & {A}\ar[d]^e \\
      {A\gray\Delta[1]_t\gray\Delta[1]_t}
      \ar[r]_-{t}\ar@{..>}[ur]|{\exists \bar{t}} &
      {B}}
  \end{equation}
  wherein the left hand vertical of the left hand square is the corner
  tensor
  \begin{equation*}
    (\overinc\subseteq_r:\emptyset->B.)\cgray
    (\overinc\subseteq_r:\Lambda^0[1]->\Delta[1]_t.)
  \end{equation*}
  and the stratified set $X$ is defined to be the regular subset which
  makes the left hand vertical of the right hand square into
  the corner tensor:
  \begin{equation*}
    (\overinc\subseteq_r:\emptyset->A.)\cgray
    (\overinc\subseteq_r:\Lambda^0[1]->\Delta[1]_t.)\cgray
    (\overinc\subseteq_r:\boundary\Delta[1]->\Delta[1]_t.)
  \end{equation*}
  These are both corner tensors of an inclusion with the elementary
  anodyne extension $\overinc\subseteq_r:\Lambda^0[1]->\Delta[1]_t.$
  so we may apply corollary~\ref{anod.tensor.stab.cor} to show that
  they both have the LLP with respect to the complicial fibration
  $\arrow e:A->B.$. This fact explains the existence of the lift
  $\bar{k}$ in the left hand square of display~(\ref{def.lift.1}),
  which we may use to define a stratified map $\arrow m:B->A.$ by
  $m(b)\defeq \bar{k}(b,1)$. Furthermore, the commutativity of its upper
  triangle tells us that $\bar{k}$ is a simple homotopy from $e'$ to $m$
  whereas its lower triangle tells us that
  $e(m(b))=e(\bar{k}(b,1))=k(b,1)=b$ (amongst other things).

  Now we may turn to the right hand square of
  display~(\ref{def.lift.1}), and proceed to define the stratified map
  $f$, whose domain $X$ is the stratified set
  $A\gray((\Delta[1]_t\gray\boundary\Delta[1])\cup
  (\Lambda^0[1]\gray\Delta[1]_t))$. This splits naturally into three
  components each of which is isomorphic to $A\gray\Delta[1]_t$ and
  upon which we define $f$ in a piecewise manner:
  \begin{equation*}
    f(a,\alpha,0) \defeq \bar{k}(e(a),\alpha) \mkern40mu
    f(a,\alpha,1) \defeq  a \mkern40mu
    f(a,0,\beta) \defeq h(a,\beta)
  \end{equation*}
  In other words, the first two clauses specify how $f$ acts on the
  disjoint components of $A\gray\Delta[1]_t\gray\boundary\Delta[1]$
  and the last one specifies how it acts on
  $A\gray\Lambda^0[1]\gray\Delta[1]_t$. To check that $f$ is well
  defined it is enough to observe that the pieces of its definition
  match at the ``corners'' where they meet, since
  $\bar{k}(e(a),0)=e'(e(a))= h(a,0)$ and $h(a,1)=a$, and that it
  respects the stratification on each component (since $\bar{k}$, $e$
  and $h$ are all stratified maps). Furthermore, comparing the
  definition of $f$ with the properties of the boundary of $t$ laid
  out at the bottom of observation~\ref{hty.eqv.coh} and applying the
  defining property $e\circ\bar{k}=k$ of $\bar{k}$ it is now easily
  seen that the right hand square in display~(\ref{def.lift.1})
  commutes and thus that the lift $\bar{t}$ exists as advertised
  there.

  Finally, it remains to define the simple homotopy $\arrow
  d:A\gray\Delta[1]_t->A.$ by letting
  $d(a,\beta)\defeq\bar{t}(a,1,\beta)$. Applying the commutativities of
  the triangles in the right hand square of display~(\ref{def.lift.1})
  and the properties of $t$ given in observation~\ref{hty.eqv.coh} we
  discover that
  \begin{align*}
    d(a,0)&{}=\bar{t}(a,1,0)=f(a,1,0)=\bar{k}(e(a),1)=m(e(a)) \\
    d(a,1)&{}=\bar{t}(a,1,1)=f(a,1,1) = a \\
    e(d(a,\beta))&{}=e(\bar{t}(a,1,\beta))=t(a,1,\beta)=e(a)
  \end{align*}
  thus verifying that $d$ completes the data required to demonstrate
  that $e$ is a deformation retraction as required.

  All that remains for us is to prove
  $\text{\ref{he.2}}\Rightarrow\text{\ref{he.3}}$, so suppose that $e$
  is a deformation retraction witnessed by the stratified map $m$ and
  the simple homotopy $d$ of definition~\ref{def.retr.defn} and consider the
  lifting problem depicted in the left hand square below:
  \begin{equation}\label{more.liftings}
    \xymatrix@=2em{
      {U}\ar@{u(->}[d]_{\subseteq_s}\ar[r]^{f} & {A}\ar[d]^{e} \\
      {V}\ar[r]_{g} & {B} 
    }\mkern50mu
    \xymatrix@R=2em@C=4em{
      {(U\gray\Delta[1]_t)\cup(V\gray\Lambda^0[1])}
      \ar@{u(->}[d]_{\subseteq_s}\ar[r]^-{\bar{f}} & {A}\ar[d]^{e} \\
      {V\gray\Delta[1]_t}\ar[r]_{g\circ\pi_V}\ar@{..>}[ur]|{\exists
        \bar{h}} & {B}
    }
  \end{equation}
  We provide a solution to this problem by constructing the right hand
  square, in which $\bar{f}$ is defined in a piecewise manner:
  \begin{equation*}
    \bar{f}(u,\alpha) \defeq d(f(u),\alpha)\text{ on
      $U\gray\Delta[1]_t$}\mkern20mu\text{and}\mkern20mu
    \bar{f}(v,0) \defeq m(g(v)) \text{ on $V\gray\Lambda^0[1]$}
  \end{equation*}
  This is well defined because the actions on these components respect
  stratifications, (since $d$, $m$, $f$ and $g$ are all stratified
  maps) and they match at the intersection of their domains where
  $d(f(u),0)=m(e(f(u))=m(g(u))$, in which the former equality holds
  because $d$ is a simple homotopy from $m\circ e$ to $\id_A$ and the latter
  one simply follows from the commutativity of the original lifting
  problem. Notice now that using the properties of $e$, $m$ and $d$ as
  the components of a deformation retraction and the commutativity of
  our original lifting problem again we have
  \begin{align*}
    e(\bar{f}(u,\alpha)) &{}= e(d(f(u),\alpha))=e(f(u))=g(u)  \\
    e(\bar{f}(v,0)) & {}= e(m(g(v))) = g(v)
  \end{align*}
  or, in other words, the square displayed does indeed commute.
  Furthermore its left hand vertical is the corner tensor
  \begin{equation*}
    (\overinc\subseteq_s:U->V.)\cgray
    (\overinc\subseteq_r:\Lambda^0[1]->\Delta[1]_t.)
  \end{equation*}
  which has the LLP with respect to the complicial fibration $e$ by
  corollary~\ref{anod.tensor.stab.cor} and so the lift $\bar{h}$
  (doted arrow) exists as depicted in the right hand square of
  display~(\ref{more.liftings}). Now it is trivially verified,
  directly from the properties of $\bar{h}$ as the stated lift, that
  the stratified map $\arrow \bar{g}:V->A.$ defined by
  $\bar{g}(v)=\bar{h}(v,1)$ is the required solution to the left hand
  lifting problem. Ultimately, it follows that $e$ is a trivial fibration
  since, in particular, we have shown that it has the RLP with respect
  to all boundary and thin simplex inclusions.
\end{proof}

\begin{obs}\label{tfib.are.hty}
  The proofs of the implications
  $\text{\ref{he.3}}\Rightarrow\text{\ref{he.2}}$ and
  $\text{\ref{he.2}}\Rightarrow\text{\ref{he.1}}$ in in the last lemma
  were independent of any weak compliciality assumptions, so we may
  infer that {\bf any} trivial fibration between arbitrary stratified
  sets is a deformation retraction and that these in turn are homotopy
  equivalences.
\end{obs}
\subsection{Weak Equivalences of Stratified Sets}

\begin{defn}\label{J.defn}\label{weak.equiv.def}
  For this subsection and the next we will do everything relative to a
  fixed (small) set of inclusions $J$ in $\Strat$, which we assume
  satisfies the condition
  \begin{enumerate}[label=(\roman*)]
  \item\label{J.cond.1} each elementary anodyne extension is an
    element of $J$
  \end{enumerate}
  thereby ensuring that every $J$-fibrant object is a weak complicial
  set. 

  We say that a stratified map $\arrow w:X->Y.$ is a {\em $J$-weak
    equivalence\/} if and only if the associated stratified map
  $\arrow\hom(w,A):\hom(Y,A)->\hom(X,A).$ is a homotopy equivalence
  for each $J$-fibrant stratified set $A$ and we let $\mathcal{W}_J$
  denote the class of all $J$-weak equivalences in $\Strat$.
  Unless otherwise stated, we will generally also assume that our set
  $J$ satisfies the condition
  \begin{enumerate}[label=(\roman*), resume]
  \item\label{J.cond.2} each element of $J$ is a $J$-weak equivalence
  \end{enumerate}
  which postulates a stability property closely related to the result
  established in corollary~\ref{anod.tensor.stab.cor} for anodyne
  extensions.

  The construction to follow provides a Quillen model structure whose
  fibrant objects are the $J$-fibrant stratified sets and whose
  fibrations between fibrant objects are precisely the $J$-fibrations
  between those stratified sets. This will allow us to construct model
  structures whose fibrant objects are weak complicial sets,
  $n$-trivial weak complicial sets, quasi-categories under their
  standard stratification and so forth.

  To make our nomenclature match with that of previous sections we
  shall call the $J$-fibrant objects {\em $J$-weak complicial sets},
  the $J$-fibrations {\em $J$-complicial fibrations}, the $J$-cell
  complexes {\em $J$-anodyne extensions\/} and so on. Also let
  $\Wcs_J$ denote the full subcategory of $\Strat$ whose objects are
  the $J$-weak complicial sets.
\end{defn}

\begin{lemma}\label{ccof.char}
  Suppose that $J$ is a small set of stratified inclusions that
  satisfies condition~\ref{J.cond.1} of definition~\ref{J.defn} and
  suppose that $\inc e:U->V.$ is an inclusion of stratified sets then
  the following are equivalent:
  \begin{enumerate}[label=(\roman*)]
  \item\label{ccof.char.1} $e$ is a $J$-weak equivalence,
  \item\label{ccof.char.2} $\arrow\hom(e,A):\hom(V,A)->\hom(U,A).$ is
    a trivial fibration for all $J$-weak complicial sets $A$, and
  \item\label{ccof.char.3} for all inclusions $\inc i:X->Y.$ the
    corner tensor $e\cgray i$ has the LLP with respect to each
    $J$-weak complicial set $A$.
  \end{enumerate}
\end{lemma}

\begin{proof}
  Every $J$-weak complicial set is, in particular, a weak complicial
  set so we may apply theorem~\ref{comp.hom.stab.thm} to the inclusion
  $e$ to show that $\hom(e,A)$ is a complicial fibration of weak
  complicial sets whenever $A$ is a $J$-weak complicial
  set. Consequently, applying lemma~\ref{cfib.equiv.lem} we see that
  $\hom(e,A)$ is a homotopy equivalence if and only if it is a trivial
  fibration. So, quantifying over all $J$-weak complicial sets and
  applying the $J$-weak equivalence definition, we have established
  the equivalence $\text{\ref{ccof.char.1}}\Leftrightarrow
  \text{\ref{ccof.char.2}}$. The remaining equivalence follows
  routinely by applying observation~\ref{corner.fibration} to the
  adjunction $e\cgray - \dashv \chom(e,*)$.
\end{proof}

\begin{eg}\label{comp.mod.eg}
  Recasting the result of corollary~\ref{anod.tensor.stab.cor} using
  observation~\ref{corner.fibration} and applying the last lemma we
  may verify that the countable set
  \begin{align*}
    J_{c} \defeq {} & \{ \overinc\subseteq_r:\Lambda^k[n]->\Delta^k[n].\mid
    n=1,2,... \text{ and } 0\leq k\leq n \} \cup {} \\
    & \{ \overinc\subseteq_e:\Delta^k[n]'->\Delta^k[n]''.\mid
    n=2,... \text{ and } 0\leq k\leq n \}
  \end{align*}
  of all elementary anodyne extensions provides a {\em minimal\/} set
  satisfying the conditions given in definition~\ref{J.defn}. 
\end{eg}

\begin{obs}[homotopy equivalence implies $J$-weak equivalence]
  \label{hty=>weak}
  The contravariant functor $\hom(*,A)$ has a canonical enrichment
  whose action on homsets $\overarr:\hom(X,Y)->\hom(\hom(Y,A),
  \hom(X,A)).$ is constructed by taking the dual of the composition
  $\arrow\circ:\hom(Y,A)\gray\hom(X,Y)->\hom(X,A).$ of $\Strat\sgray$
  under the appropriate closure adjunction. In this way, for each weak
  complicial set $A$ we obtain a $\Strat\sgray$-enriched functor
  $\arrow\hom(*,A):\Strat\sgray\op->\Wcs.$ which carries the thin
  1-simplices in the homsets of $\Strat\sgray$, that is to say simple
  homotopies, to thin 1-simplices in the homsets of $\Wcs$ and thus
  preserves the homotopy relation $\sim$ between stratified
  maps. Consequently it maps left (resp.\ right) homotopy inverses to
  right (resp.\ left) homotopy inverses and therefore preserves
  homotopy equivalences, thus demonstrating that any homotopy
  equivalence of stratified sets is a $J$-weak equivalence.
\end{obs}

\begin{obs}[a partial converse]\label{weak+cmpl=>hty}
  Suppose that $\arrow w:A->B.$ is a $J$-weak equivalence between
  $J$-weak complicial sets then, since $A$ is a $J$-weak complicial
  set, we know that the associated stratified map
  $\arrow\hom(w,A):\hom(B,A)->\hom(A,A).$ has a homotopy inverse
  $\arrow\bar{w}:\hom(A,A)->\hom(B,A).$ for which the right inverse
  homotopy $\hom(w,A)\circ\bar{w}\sim\id_{\hom(A,A)}$ may be witnessed
  by a simple homotopy $\arrow
  \bar{h}:\hom(A,A)\gray\Delta[1]_t->\hom(A,A).$. So if we define maps
  $\arrow w':B->A.$ by $w'=\bar{w}(\id_A)$ and $\arrow h:A\gray
  \Delta[1]_t->A.$ by $h(a,\alpha)=\bar{h}(\id_A\cdot\eta,\alpha)(a)$
  then we have
  \begin{align*}
    h(a,0)&{}=\bar{h}(\id_A,0)(a)=\hom(w,A)(\bar{w}(\id_A))(a)=w'(w(a))\\
    h(a,1)&{}=\bar{h}(\id_A,1)(a)=\id_{\hom(A,A)}(\id_A)(a)=a
  \end{align*}
  or, in other words, $h$ is a simple homotopy from $w'\circ w$ to
  $\id_A$. 

  Applying $\hom(-,B)$ to this, and consulting the last observation,
  we obtain a simple homotopy from $\hom(w,B)\circ\hom(w',B)$ to
  $\id_{\hom(B,B)}$ thus demonstrating that $\hom(w',B)$ is a right
  equivalence inverse of $\hom(w,B)$. This latter map is, however, a
  homotopy equivalence, since $B$ is a $J$-weak complicial set and $w$
  is a $J$-weak equivalence, so it follows that $\hom(w',B)$ is also a
  left equivalence inverse of $\hom(w,B)$ by
  observation~\ref{equiv.basic.obs}. Finally, applying the argument
  used above to obtain $h$ from $\bar{h}$ to the resulting simple
  homotopy $\bar{k}$ from $\hom(w',B)\circ\hom(w,B)$ to
  $\id_{\hom(B,B)}$ we obtain a simple homotopy $k$ from $w\circ w'$
  to $\id_B$ thus completing the demonstration that $w$ is a homotopy
  equivalence.
\end{obs}

\begin{obs}\label{inc.corner.obs}
  Fix an inclusion $\inc i:X->Y.$ and observe that we may apply
  clause~\ref{ccof.char.3} of lemma~\ref{ccof.char} and the fact that
  every element of $J$ is both a $J$-weak equivalence and an inclusion
  (under the assumptions of definition~\ref{J.defn}) to show that for
  each $e\in J$ the corner tensor $e\cgray i$ has the LLP with respect
  to each $J$-weak complicial set $A$. Applying
  observation~\ref{corner.fibration} to the adjunction $-\cgray
  i\dashv \chom(i,*)$, we find that this is equivalent to saying that
  $\arrow \hom(i,A):\hom(Y,A)->\hom(X,A).$ has the RLP with respect to
  each inclusion in $J$ and that it is thus a $J$-complicial
  fibration.

  Applying this result to the (unique) inclusion $\inc
  !:\emptyset->X.$ whose domain is the empty stratified set, we find
  that whenever $A$ is a $J$-weak complicial set the (also unique)
  stratified map $\arrow\hom(!,A):\hom(X,A)->\hom(\emptyset,A)\cong
  1.$ is a $J$-complicial fibration. In other words, in that
  circumstance $\hom(X,A)$ is also a $J$-weak complicial set.
\end{obs}

\begin{obs}[verifying the conditions of Jeffery Smith's theorem]
  \label{j.smith.verif.1}\label{anod=>weak}
  Our intention is to show that the set of inclusions $I$ and the
  class of $J$-weak equivalences satisfy the conditions of
  theorem~\ref{jsmith.thm}. Indeed, we are already in a position to
  verify the first three clauses of its statement.
  \begin{enumerate}[label=(\arabic*)]
  \item Observation~\ref{equiv.basic.obs} tells us that the class of
    homotopy equivalences is closed under retracts and enjoys the
    2-of-3 property. It is thus clear, directly from
    definition~\ref{weak.equiv.def} and the functoriality of
    $\hom(*,A)$, that the class of $J$-weak equivalences $\mathcal{W}_J$
    also possesses these properties.
  \item From observation~\ref{tfib.are.hty} we know that all
    $I$-fibrations (trivial fibrations) are homotopy equivalences and
    thus that they are all $J$-weak equivalences by
    observation~\ref{hty=>weak}.
  \item Observation~\ref{strat.inc.obs} tells us that the class of all
    inclusions of stratified sets is closed under pushout, transfinite
    composition and retraction and also that it is equal to
    $\cof(I)$. Lemma~\ref{ccof.char} reveals that an inclusion
    $\inc e:U->V.$ is in $\mathcal{W}_J$ if and only if
    $\arrow\hom(e,A):\hom(V,A)->\hom(U,A).$ is a trivial fibration for
    each $J$-weak complicial set $A$. However $\hom(*,A)$ carries
    transfinite composites and pushouts to the corresponding limits in
    $\Strat$ and any class of the form $\fib(I)$, such as the class of
    trivial fibrations, is closed under those limits. Combining these
    facts, we see that if $\inc e:U->V.$ is a transfinite composite of
    pushouts of elements of $\cof(I)\cap\mathcal{W}_J$ then
    $\arrow\hom(e,A):\hom(V,A)->\hom(U,A).$ is a corresponding
    transfinite limit of pullbacks of trivial fibrations which is thus
    itself a trivial fibration. That however implies that $\hom(e,A)$
    is a trivial fibration for each $J$-weak complicial set $A$ and
    thus that the inclusion $e$ is a $J$-weak equivalence as required.
  \end{enumerate}

  Under our assumption that the elements of $J$ are all inclusions,
  condition~\ref{J.cond.2} of definition~\ref{J.defn} simply states
  that the set $J$ is a subset of the class
  $\cof(I)\cap\mathcal{W}_J$.  However, the last of the properties
  above verifies that this latter class is closed under the operations
  used to derive the class of $J$-anodyne extensions from $J$
  itself. It follows, therefore, that all $J$-anodyne extensions are
  also $J$-weak equivalences.

  We demonstrate the fourth condition of theorem~\ref{jsmith.thm} by
  showing that the class of $J$-weak equivalences is actually an
  {\em accessible class\/} of maps (cf.\ observation~\ref{inj.acc} and
  Beke~\cite{Beke:2000:SheafHomotopy}), which we do in two steps:
\end{obs}

\begin{obs}[the class of trivial fibrations is accessible]
  By observation~\ref{LFP.quasitopos} we know that the category
  $\Strat$ is locally finitely presentable. Furthermore we may argue,
  as in observation~\ref{small.obj.arg}, that the full subcategory
  $\TFib$ of $\Strat^2$ whose objects are the trivial fibrations is
  simply the injectivity class associated with the set of squares of
  the form:
  \begin{equation*}
    \xymatrix@R=2em@C=3em{
      {\boundary\Delta[n]}\ar@{u(->}[r]^-{\subseteq_r}
      \ar@{u(->}[d]_{\subseteq_r} & {\Delta[n]} 
      \ar[d]^{\id} \\
      {\Delta[n]}\ar@{->}[r]_{\id} & {\Delta[n]}
    }\mkern40mu
    \xymatrix@R=2em@C=3em{
      {\Delta[n]}\ar@{u(->}[r]^-{\subseteq_e}
      \ar@{u(->}[d]_{\subseteq_e} & {\Delta[n]_t} 
      \ar[d]^{\id} \\
      {\Delta[n]_t}\ar@{->}[r]_{\id} & {\Delta[n]_t}
    }
  \end{equation*}
  It follows that we may apply observation~\ref{inj.acc} to show that
  the class of trivial fibrations is an accessible class of maps in
  $\Strat$. Indeed, with a little more work we may show that $\TFib$
  is $\aleph_1$-accessible and $\aleph_0$-accessibly embedded in
  $\Strat^2$, although we will not need that result here.
\end{obs}

\begin{obs}[the class of $J$-weak equivalences is accessible]
  \label{j.smith.verif.2}
  In a similar fashion we may describe the class of $J$-complicial
  fibrations between $J$-weak complicial sets as an injectivity class
  $\CFib_J$ in $\Strat^2$. To be precise the objects of this subcategory
  are the morphisms which are injective with respect to the
  squares of the form
  \begin{equation}\label{inj.squares.2}
    \xymatrix@=2em{
      {\emptyset}\ar@{->}[r]\ar[d] & {\emptyset} 
      \ar[d] \\
      {U_j}\ar@{u(->}[r]_{j} & {V_j}}\mkern40mu
    \xymatrix@=2em{
      {U_j}\ar@{u(->}[r]^-{j}
      \ar@{u(->}[d]_{j} & {V_j} 
      \ar[d]^{\id} \\
      {V_j}\ar@{->}[r]_{\id} & {V_j}
    }
  \end{equation}
  for each $\inc j:U_j->V_j.$ in $J$.  Here injectivity with respect
  to squares of the left hand form ensures that the codomains of
  morphisms in our class are $J$-weak complicial sets and injectivity
  with respect to squares of the right hand form ensures that these
  morphisms are themselves $J$-complicial fibrations. We will use $K$
  to denote the set of those morphisms of $\Strat^2$ depicted in the
  above display.

  Now we may apply observation~\ref{small.obj.arg} to the locally
  finitely presentable category $\Strat^2$ and the set $K$ of its
  squares to obtain an accessible weak reflection of $\Strat^2$ into
  $\CFib_J$.  We denote the weak reflection of an object $\arrow
  f:X->Y.$ of $\Strat^2$ into $\CFib_J$ by $\arrow f^*:X_f^*->Y_f^*.$
  and the use the notation
  \begin{equation}\label{unit.cpt.1}
    \xymatrix@R=2em@C=3em{
      {X}\ar@{u(->}[r]^-{\eta^d_f}\ar[d]_{f} &
      {X_f^*}\ar[d]^{f^*} \\
      {Y}\ar@{u(->}[r]_-{\eta^c_f} & {Y_f^*}
    }
  \end{equation}
  for the associated component of the unit of this weak
  reflection. This is simply the component of the defining colimiting
  cone of $(-)^*$ from the first element $\id_{\Strat^2}$ of the chain
  constructed in observation~\ref{small.obj.arg}. Since all colimits
  in $\Strat^2$ and its endo-functor category are constructed
  pointwise in $\Strat$, we know that the maps $\eta^d_f$ and
  $\eta^c_f$ are $J$-anodyne extensions because they are constructed
  as transfinite composites of pushouts of coproducts of the
  horizontal maps in the squares in display~(\ref{inj.squares.2}),
  each of which is an element of $J$ or an identity. Now $J$-anodyne
  extensions are $J$-weak equivalences, by
  observation~\ref{anod=>weak}, so we may apply the 2-of-3 property
  for these to the square in display~\ref{unit.cpt.1} to show that $f$
  is a $J$-weak equivalence if and only if $f^*$ is a $J$-weak
  equivalence. However, since $f^*$ is a $J$-complicial fibration of
  $J$-weak complicial sets we may apply
  observation~\ref{weak+cmpl=>hty} and lemma~\ref{cfib.equiv.lem} to
  show that $f^*$ is a $J$-weak equivalence if and only if it is a
  trivial fibration.

  In summary, we have constructed an accessible endo-functor $(-)^*$
  of $\Strat^2$ with the property that a stratified map $\arrow
  f:X->Y.$ is an object of $\WEqv_J$, the full subcategory of
  $\Strat^2$ whose objects are $J$-weak equivalences, if and only if
  its weak reflection $\arrow f^*:X_f^*->Y_f^*.$ is an object of the
  accessible and accessibly embedded full subcategory $\TFib$ of
  trivial fibrations. In other words, $\WEqv_J$ is a pseudo-pullback
  of $\TFib$ along the endo functor $(-)^*$ and consequently we may
  apply theorem~5.1.6 of~\cite{Makkai:1989:Access} to show that it too
  is an accessibly embedded, accessible subcategory of $\Strat^2$. In
  this way we have shown that the $J$-weak equivalences form an
  accessible class, which is finally all we need in order to apply
  Jeffery Smith's theorem and establish the next theorem.
\end{obs}

\subsection{The $J$-Complicial Model Structure}

\begin{thm}\label{model.str.thm}
  Each set of stratified inclusions $J$ satisfying the conditions
  given in definition~\ref{J.defn} gives rise to a cofibrantly
  generated Quillen model structure on the category $\Strat$ of
  stratified sets, called the {\em $J$-complicial model structure},
  whose:
  \begin{itemize}
  \item {\bf weak equivalences} are the $J$-weak equivalences of
    observation~\ref{weak.equiv.def},
  \item {\bf cofibrations} are simply inclusions of stratified sets,
    and whose
  \item {\bf fibrant objects} are the $J$-weak complicial sets. 
  \end{itemize}
\end{thm}

\begin{proof}
  Apply Jeffery Smith's theorem~\ref{jsmith.thm}, using our
  observations~\ref{j.smith.verif.1} and~\ref{j.smith.verif.2} to
  verify the required properties of $I$ and $\mathcal{W}$. The
  proof that the fibrant objects in this model structure are exactly
  the $J$-weak complicial sets is postponed to
  lemma~\ref{cfib+wcs=>ccfib} below.
\end{proof}

\begin{notation}\label{comp.cof.note}
  We call the model structure derived in the last theorem the {\em
    $J$-complicial model structure} and refer its trivial cofibrations
  as {\em $J$-complicial cofibrations}. Observation~\ref{anod=>weak}
  tells us that all $J$-anodyne extensions are $J$-complicial
  cofibrations, but there is no reason in general to believe that
  these classes coincide.

  Similarly, we call the fibrations of this model structure {\em
    completely $J$-complicial fibrations}. The dual of our last
  observation is that every completely $J$-complicial fibration is a
  $J$-complicial fibration but that these classes may not coincide.

  Of course, it is a general result in the theory of Quillen model
  categories that a stratified map is a $J$-complicial cofibration if
  and only if it has the LLP with respect to all completely
  $J$-complicial fibrations (and indeed vice versa), and it follows
  therefore that the class of such things is closed under pushout,
  transfinite composition and retraction.

  If we omit the prefix ``$J$-'' altogether then we will implicitly
  assume that we are working relative the set $J_c$ discussed in
  example~\ref{comp.mod.eg}. So the fibrant objects of the complicial
  model structure are precisely the weak complicial sets.
\end{notation}

\begin{obs}
  If $J\subseteq J'$ are two sets of stratified inclusions satisfying the
  conditions of definition~\ref{J.defn} then the corresponding
  complicial model structures share the same sets of cofibrations but
  differ in their sets if weak equivalences and completely complicial
  fibrations. We know, however, that every $J$-weak equivalence
  is also a $J'$-weak equivalence, which implies a corresponding
  relationship between respective classes of complicial
  cofibrations. Dually, it immediately follows that every (complete)
  $J'$-complicial fibration is a (complete) $J$-complicial fibration.
\end{obs}





\begin{obs}[localising an existing complicial model structure]
  \label{loc.mod.obs}
  If we start with a set of inclusions $J$ which satisfies the
  conditions of definition~\ref{J.defn} and $K$ is any other set of
  inclusions then a $J\cup K$-fibrant object is also $J$-fibrant so it
  follows that any $J$-weak equivalence is also a $J\cup K$-weak
  equivalence. Consequently every element of $J$ is a $J\cup K$-weak
  equivalence so to show that $J\cup K$ satisfies
  condition~\ref{J.cond.2} of definition~\ref{J.defn} all we need do
  is check that each $k\in K$ is a $J\cup K$-weak
  equivalence. However, since $J\cup K$ satisfies
  condition~\ref{J.cond.1} of definition~\ref{J.defn} we may do this
  using one of the equivalent characterisations of
  lemma~\ref{ccof.char}.

  It is a standard, and easily demonstrated, result of Quillen model
  category theory that a map $\arrow p:A->B.$ between fibrant objects
  is a fibration (resp.\ trivial fibration) if and only if it has the
  RLP with respect to all trivial cofibrations (resp.\ cofibrations)
  $\arrow i:X->Y.$ for which $X$ and $Y$ are fibrant. Now we know, by
  assumption and theorem~\ref{model.str.thm}, that $J$ gives rise to a
  complicial model structure and that any $J\cup K$-weak complicial
  set $A$ is a $J$-weak complicial set, so we may apply
  observation~\ref{inc.corner.obs} to show that any $\hom(X,A)$ is
  also a $J$-weak complicial set. In particular, if $\inc k:U_k->V_k.$
  is an element of $K$ then it follows that the domain and codomain of
  the stratified map $\arrow\hom(k,A):\hom(V_k,A)->\hom(U_k,A).$ of
  clause~\ref{ccof.char.2} in lemma~\ref{ccof.char} are $J$-weak
  complicial sets and thus that we may apply the result recalled in
  the first sentence to show that $\hom(k,A)$ is a trivial fibration
  if and only if it has the RLP with respect to all inclusions of
  $J$-weak complicial sets. Equivalently, it follows that we may show
  that $k$ is a $J\cup K$-weak equivalence by showing that it
  satisfies clause~\ref{ccof.char.3} in lemma~\ref{ccof.char} for
  those inclusions in this restricted class.
\end{obs}

\begin{eg}[the $n$-trivial complicial model structure]\label{triv.mod.eg}
  To obtain a Quillen model structure whose fibrant objects are the
  $n$-trivial weak complicial sets, define the set
  \begin{equation*}
    J_n\defeq J_c\cup\{\overinc\subseteq_e:\Delta[r]->\Delta[r]_t.\mid
    r\in\mathbb{N}\wedge r>n\}
  \end{equation*}
  for which the $J_n$-fibrant objects are those stratified sets with
  the desired properties. Now the $n$-superstructure functor $\sst_n$
  acts as the identity on $n$-trivial stratified sets and we know, by
  lemma~\ref{super.wcs}, that it preserves weak
  compliciality. Combining that observation with lemma~\ref{ccof.char}
  we find that $t_r\defeq \overinc\subseteq_r:\Delta[r]->\Delta[r]_t.$
  ($r>n$) is a $J_n$-weak equivalence if and only if its corner tensor
  $t_r\cgray i$ with each inclusion $\overinc i:X->Y.$ has the LLP
  with respect to $\sst_n(A)$ for each weak complicial set $A$. Taking
  duals under the adjunction $\triv_n\dashv\sst_n$ we find that this
  is equivalent to saying that $\triv_n(t_r\cgray i)$ has the LLP with
  respect to each weak complicial set $A$. Now, it is easily seen that
  $\triv_n$ preserves the Gray tensor $\gray$, in the sense that we
  literally have $\triv_n(X\gray Y)= \triv_n(X)\gray\triv_n(Y)$, so it
  follows that $\triv_n(t_r\cgray i)\cong\triv_n(t_r)\cgray\triv_n(i)$
  (in $\Strat^2$). However, since $r>n$ we find that $\triv_n(t_r)$ is
  actually the identity on the stratified set $\triv_n(\Delta[r])$,
  from which it follows that its corner tensor with any inclusion is
  an isomorphism and thereby demonstrate that $\triv_n(t_r\cgray i)$
  has the LLP with respect to any stratified set. This certainly
  establishes that the condition of the last paragraph holds for $J_n$
  and thus that it satisfies definition~\ref{J.defn}. Consequently it
  gives rise to a Quillen model structure whose fibrant objects we may
  show to be the $n$-trivial weak complicial sets (by applying the
  subsequent lemma).
\end{eg}

\begin{lemma}\label{cfib+wcs=>ccfib}
  If $\arrow p:A->B.$ is a $J$-complicial fibration of $J$-weak
  complicial sets then it is a completely $J$-complicial fibration. In
  particular, it follows that the fibrant objects of the
  $J$-complicial model structure are precisely the $J$-weak complicial
  sets. 
\end{lemma}

\begin{proof} (following an argument due to
  Quillen~\cite{Quillen:1967:Model}) Suppose that the inclusion $\inc
  e:U->V.$ is a $J$-complicial cofibration and consider the following
  defining diagram for the corner tensor $\chom(e,p)$:
  \begin{equation*}
    \xymatrix@R=2em@C=3em{
      \save []+<-2em,+3em>*+{\hom(V,A)}\ar[r]|-{\chom(e,p)}
      \ar@/^/[rr]^{\hom(e,A)}\ar@/_/[rd]_{\hom(V,p)}\restore &
      {\hom(V,B)\times_{\hom(U,B)}\hom(U,A)}\pbexcursion
      \ar[r]_-{q}\ar[d] & {\hom(U,A)}\ar[d]^{\hom(U,p)} \\
      & {\hom(V,B)}\ar[r]_{\hom(e,B)} & {\hom(V,B)}
    }
  \end{equation*}
  in which each object is a $J$-weak complicial set by
  observation~\ref{inc.corner.obs}, which applies here since $A$ and
  $B$ are both $J$-weak complicial sets.  Applying
  lemma~\ref{ccof.char} we may show that the maps $\hom(e,A)$ and
  $\hom(e,B)$ are trivial fibrations, from which it follows that the
  pullback $q$ of the latter is also a trivial fibration and thus that
  these are all homotopy equivalences by
  lemma~\ref{cfib.equiv.lem}. Now we can use the 2-of-3 property to
  demonstrate that $\chom(e,p)$ is also a homotopy equivalence and
  furthermore show that it is a complicial fibration, by applying
  theorem~\ref{comp.hom.stab.thm} to the inclusion $e$, thus allowing
  us to infer that it is a trivial fibration by applying
  lemma~\ref{cfib.equiv.lem}.

  Finally, the fact that trivial fibrations have the right lifting
  property with respect to the (unique) inclusion
  $\inc !:\emptyset->\Delta[0].$ implies that they are all
  surjective on $0$-simplices. Applying this to $\chom(e,p)$ we
  immediately see that $p$ has the RLP with respect to $e$ as required.
\end{proof}

\begin{cor}\label{ccof.char.b}
  An inclusion $\inc e:U->V.$ is a $J$-complicial cofibration if and only
  if it has the LLP with respect to each $J$-complicial fibration $\arrow
  p:A->B.$ of $J$-weak complicial sets.
\end{cor}

\begin{proof}
  The ``only if'' part was established in the last lemma. For the
  ``if'' part observe that if $A$ is a $J$-weak complicial set, and
  $\inc i:X->Y.$ is any inclusion then we may apply
  observation~\ref{inc.corner.obs} to show that
  $\arrow\hom(i,A):\hom(Y,A)->\hom(X,A).$ is a $J$-complicial
  fibration of $J$-weak complicial sets. So by the assumption of the
  statement we know that $e$ has the LLP with respect to $\hom(i,A)$
  and may apply observation~\ref{corner.fibration} under the
  adjunction $-\cgray i\dashv\chom(i,*)$ to show that this is
  equivalent to saying that $e\cgray i$ has the LLP with respect
  $A$. It follows that $e$ satisfies clause~\ref{ccof.char.3} of
  lemma~\ref{ccof.char} by which we may infer that it is a
  $J$-complicial cofibration as postulated.
\end{proof}

\begin{obs}\label{cmod.mon.gray.obs}
  The $J$-complicial model structure is {\em monoidal\/} with respect
  to the Gray tensor product $\gray$. This amounts to showing that if
  $\inc e:U->V.$ is a $J$-complicial cofibration and $\inc i:X->Y.$ is any
  inclusion then their corner tensor $e\cgray i$ is also a $J$-complicial
  cofibration. 

  We use clause~\ref{ccof.char.3} of lemma~\ref{ccof.char} to
  demonstrate this result, so suppose that $\inc j:S->T.$ is any other
  inclusion of stratified sets and consider $(e\cgray i)\cgray j\cong
  e\cgray (i\cgray j)$ (in $\Strat^2$). We know that $i\cgray j$ is an
  inclusion and that $e$ is a $J$-complicial cofibration, so we may
  apply lemma~\ref{ccof.char}\ref{ccof.char.3} to show that
  $e\cgray(i\cgray j)$ has the LLP with respect to each $J$-weak
  complicial set. It follows that the isomorphic map $(e\cgray
  i)\cgray j$ also has this property for each $j$, which fact allows
  us to apply the same characterisation in the reverse direction to
  show that $e\cgray i$ is a $J$-complicial cofibration as required.
\end{obs}

\subsection{Monoidality of the  Complicial Model Structure}

In this subsection we round out the results presented in
lemma~\ref{super.wcs} and corollary~\ref{anod.tensor.stab.cor} by
extending them to encompass complicial cofibrations and to the lax
Gray tensor product. As a result we establish that the complicial
model structure makes $\Strat$ into a {\em monoidal model category\/}
with respect to either one of the Gray tensors $\gray$ or
$\laxgray$. Indeed, it is trivially the case that it is also a
monoidal model category with respect to the join $\join$, a result we
leave to the reader to verify.

\begin{cor}[of lemma~\ref{super.wcs}]\label{ccof.th.cor}
  For each $n\in\mathbb{N}$ the $n$-trivialisation functor $\triv_n$
  preserves complicial cofibrations. It follows that if $\inc e:U->V.$
  is a complicial cofibration then so is the associated inclusion
  $\overinc :\triv_n(U)\vee_U V->\triv_n(V).$.
\end{cor}

\begin{proof}
  Corollary~\ref{ccof.char.b} tells us that
  $\inc\triv_n(e):\triv_n(U)->\triv_n(V).$ is a complicial cofibration
  if and only if it has the LLP with respect to complicial fibrations
  of complicial sets $\arrow p:A->B.$. Taking duals under the
  adjunction $\triv_n\dashv\sst_n$ we find that this is the case iff
  $\inc e:U->V.$ has the LLP with respect to $\arrow
  \sst_n(p):\sst_n(A)->\sst_n(B).$. This latter fact follows directly
  from lemma~\ref{super.wcs} which tells us that $\sst_n$ preserves
  weak complicial sets and complicial fibrations.

  For the second part of the statement, assume w.l.o.g.\ that $e$ is
  actually a stratified subset inclusion $\overinc\subseteq_s:U->V.$
  and that, consequently, the pushout $\triv_n(U)\vee_U V$ is
  actually the union $\triv_n(U)\cup V$ of subsets in
  $\triv_n(V)$. Now consider the following diagram
  \begin{equation*}
    \xymatrix@R=2em@C=3em{
      {U}\ar@{u(->}[r]^-{\subseteq_s}
      \ar@{u(->}[d]_{\subseteq_e} & {V}
      \ar@{u(->}[d]^{\subseteq_e} & \\
      {\triv_n(U)}\ar@{u(->}[r]_-{\subseteq_s} &
      {\triv_n(U)\cup V}\poexcursion
      \ar@{u(->}[r]_-{\subseteq_e} & {\triv_n(V)}
    }
  \end{equation*}
  in which it is easily checked that the square is a pushout as
  marked. By assumption the upper horizontal here is a complicial
  cofibration, so it follows that the lower horizontal in this square
  is also a complicial cofibration. Furthermore, in the last paragraph
  we demonstrated that the composite of the lower horizontals is a
  complicial cofibration. So we may apply the 2-of-3 property to show
  that the right hand lower horizontal is also a complicial
  cofibration as required.
\end{proof}

\begin{thm}\label{ccfib.ctens.thm}
  If $\inc e:U->V.$ is a complicial cofibration and $\inc i:X->Y.$ is
  any inclusion of stratified sets then their corner tensors
  $e\claxgray i$ and $e\cpretens i$ are complicial
  cofibrations. Dually it is the case that $i\claxgray e$ and
  $i\cpretens e$ are also complicial cofibrations. 
\end{thm}

\begin{proof}
  We avoided proving this kind of result in
  section~\ref{gray.tensor.sec} because at the time we had no
  completely satisfactory way of relating the properties of corner
  tensors with respect to the lax Gray tensors $\laxgray$ and
  $\pretens$. However, the complicial cofibrations of our model
  structure provide a solution to this problem and allow us to easily
  prove properties of $\laxgray$ using the colimit preservation
  properties of $\pretens$. To that end consider the following
  commutative square
  \begin{equation*}
    \xymatrix@R=2em@C=5em{
      {(U\pretens Y)\vee_{U\pretens X}(V\pretens X)}
      \ar@{u(->}[r]^-{e\cpretens i}\ar@{u(->}[d]_{\subseteq_e} & {V\pretens Y}
      \ar@{u(->}[d]^{\subseteq_e} \\
      {(U\laxgray Y)\vee_{U\laxgray X}(V\laxgray X)}
      \ar@{u(->}[r]_-{e\claxgray i} & {V\laxgray Y}
    }
  \end{equation*}
  in which the right hand vertical is an anodyne extension by lemma
  139 of~\cite{Verity:2006:Complicial} (cf.\
  observation~\ref{pretens.gen.obs}), as is the left hand vertical
  since it is constructed as a pushout of such comparison
  maps. Applying observation~\ref{anod=>weak} we see that these are
  both weak equivalences and therefore that we can apply the 2-of-3
  property to show that the upper horizontal $e\cpretens i$ is a
  complicial cofibration if and only if the lower horizontal
  $e\claxgray i$ is a complicial cofibration.

  Applying this observation to lemma~\ref{left.anod.prod.lem} we find
  that the corner tensor $\cpretens$ of an elementary left or inner
  anodyne extension with a boundary or thin simplex inclusion is a
  complicial cofibration. However, the class of complicial
  cofibrations is closed under pushout and transfinite composition
  (cf.\ notation~\ref{comp.cof.note}) so we may apply
  lemma~\ref{corner.cofibration} to the tensor $\pretens$, which does
  preserve colimits in each variable, to show that if $e$ is a left
  anodyne extension and $i$ is any inclusion then their corner tensor
  $e\cpretens i$ is a complicial cofibration. In particular, we may
  now apply corollary~\ref{cfib+wcs=>ccfib} to show that $e\cpretens
  i$ has the LLP with respect to each complicial fibration $\arrow
  p:A->B.$ of weak complicial sets. Dually, applying
  observation~\ref{corner.fibration} to the adjunction $-\cpretens
  i\dashv \clax_r(i,*)$ we see that $\clax_r(i,p)$ is a left
  complicial fibration.

  Finally, to extend the work of the last paragraph to all complicial
  cofibrations $e$ observe that the result given in its last sentence
  now allows us to apply the argument of
  theorem~\ref{comp.hom.stab.thm} with respect to $\cpretens$ and
  $\clax_r$.  This demonstrates that $\lax_r(X,A)$ is a weak
  complicial set whenever $A$ is and that for each inclusion $i$ and
  complicial fibration $p$ of weak complicial sets the corner closure
  $\clax_r(i,p)$ is a completely complicial fibration, and thus enjoys
  the RLP with respect to any complicial cofibration $e$. Dualising
  that result using observation~\ref{corner.fibration} and the
  adjunction $-\cpretens i\dashv \clax_r(i,*)$ we find that
  $e\cpretens i$ has the LLP with respect to any such
  $p$. Consequently we may apply corollary~\ref{ccof.char.b} to show
  that $e\cpretens i$ is a complicial cofibration, and thus that
  $e\claxgray i$ is also a complicial cofibration by the observation
  of the first paragraph of this proof.
\end{proof}

\subsection{A Model Structure for Joyal's Quasi-Categories}

Finally, we derive model structures on $\Strat$ and $\Simp$ whose
fibrant objects are Joyal's quasi-categories. This latter structure
was originally constructed by Joyal, although at the time of writing
none of his published papers contain the detail of his
construction. We provide the following construction as an independent
verification of his work and in order to provide us with a
presentation of this model structure which we shall apply in our
forthcoming work on {\em internal
  quasi-categories\/}~\cite{Verity:2006:IntQCat}.

\begin{defn}\label{Jq.defn}
  Define a set of inclusions
  \begin{equation*}
    J_q\defeq J_1\cup\{\overinc\subseteq_e:
    \usimp{E}_2->E_2.\}
  \end{equation*}
  where $J_1$ is the set defined in example~\ref{triv.mod.eg} whose
  fibrant objects are $1$-trivial weak complicial sets. 
\end{defn}

\begin{defn}\label{cestrat.not}
  Consider the following set of stratified inclusions:
  \begin{displaymath}
    Q\defeq \{ \overinc\subseteq_e:\Delta[n]->\Delta[n]_t. \mid
    n=2,3,4,...\} \cup \{ \overinc\subseteq_e:\usimp{E}_2->E_2.,\, 
    \overinc\subseteq_r:E^-_1->E_2.\}
  \end{displaymath}
  Of course we know that a stratified set $X$ has the RLP with respect
  to the given thin simplex inclusions iff it is
  $1$-trivial. Furthermore, $X$ has the RLP with respect to the
  remaining inclusions in $Q$ iff every $1$-simplex with an
  equivalence inverse is thin (RLP w.r.t.\
  $\overinc\subseteq_e:\usimp{E}_2->E_2.$) {\em and\/} every thin
  $1$-simplex has an equivalence inverse (RLP w.r.t.\
  $\overinc\subseteq_r:E^-_1->E_2.$). It follows that the
  stratification of a $Q$-fibrant stratified set is completely
  determined by the structure of its underlying simplicial set.

  Now, if $X$ is a simplicial set then let $X^e$ denote the entire
  superset of $\triv_1(X)$ constructed by making thin those
  $1$-simplices $x\in X$ for which there is a simplicial map $\arrow
  \hat{x}:\usimp{E}_2->X.$ with $\hat{x}(e^-_1)=x$. It is easily seen
  that $X^e$ is $Q$-fibrant for any simplicial set $X$, so it follows
  from the last paragraph that each simplicial set carries a unique
  $Q$-fibrant stratification. Furthermore, it is clear that we may
  construct the entire inclusion $\overinc\subseteq_e:X->X^e.$ as a
  pushout of a coproduct of copies of thin simplex inclusions and the
  inclusion $\overinc\subseteq_e:\usimp{E}_2->E_2.$, from which it
  follows that it is both a relative $Q$-cell complex and a
  $J_q$-anodyne extension.

  Clearly any simplicial map $\arrow f:X->Y.$ is the underlying map of
  a stratified map $\arrow f^e:X^e->Y^e.$, in other words the
  stratification operation $(-)^e$ provides us with a fully-faithful
  functor from $\Simp$ to $\Strat$ which makes the family of entire
  inclusions $\overinc\subseteq_e:X->X^e.$ into a natural
  transformation. So the stratification operation $(-)^e$ provides us
  with an equivalence between $\Simp$ and the full sub-category of
  $Q$-fibrant stratified sets in $\Strat$.
\end{defn}

\begin{obs}[quasi-categories and $J_q$-weak complicial sets]
  \label{qcat=J_q-wcs}
  Returning to example~\ref{quasi.wcs} we see that the functor $(-)^e$
  generalises the canonical stratification discussed
  there. That example tells us that a simplicial set $A$ is a
  quasi-category iff $A^e$ is a weak complicial set. We know however
  that $A^e$ has the RLP with respect to the set $Q$ of the last
  definition and that $J_q\subseteq J_c\cup Q$, so it follows that $A$
  is a quasi-category iff $A^e$ is a $J_q$-weak complicial set.

  We do not yet know that every $J_q$-weak complicial set is of the
  form $A^e$ for some quasi-category $A$, however that result would
  follow from the observations made in the course of the last
  definition as soon as we demonstrate that every $J_q$-weak
  complicial set is $Q$-fibrant.  To that end observe that every
  element of $Q$ {\em except\/} for the inclusion
  $\overinc\subseteq_r:E^-_1->E_2.$ is in $J_q$ and that it is
  possible to show that this latter inclusion is also a
  ($J_q$-)complicial cofibration, as we do in the next two
  paragraphs. Consequently, we know that any complete $J_q$-complicial
  fibration is also a $Q$-fibration, which result immediately
  specialises to the one outlined in the first sentence above.

  So to complete the proof outlined in the last paragraph, we start by
  considering the inclusion $\overinc\subseteq_r:E^-_0->E_2.$ and its
  right inverse $\epi q:E_2->E^-_0.$ (the unique map into
  $E^-_0\cong\Delta[0]$). Consider now the order preserving map
  $\arrow h:\mathbb{I}\times [1]->\mathbb{I}.$ where $\mathbb{I}$ is
  the two-point chaotic category of definition~\ref{fs.equiv.defn} and
  $h$ is defined by $h(p,0) = p$ and $h(p,1)=-$. Applying the
  categorical nerve construction of observation~\ref{cat.nerves} we
  obtain a simplicial map $\arrow h: \usimp{E}\times\Delta[1] ->
  \usimp{E}.$ which we may stratify and restrict to give a simple
  homotopy $\arrow h: E_2\gray\Delta[1]_t->E_2.$ from the identity on
  $E_2$ to the composite $\overepi q:E_2->E^-_0.\overinc\subseteq_r:
  ->E_2.$. This demonstrates that $q$ and
  $\overinc\subseteq_r:E^-_0->E_2.$ are mutual homotopy inverses and
  thus that they are both weak equivalences, by
  observation~\ref{hty=>weak}.

  Now we may factor the inclusion of the last paragraph as the
  composite of the inclusions $\overinc\subseteq_r:E^-_0->E^-_1.$ and
  $\overinc\subseteq_r:E^-_1->E_2.$. However the first of these is
  isomorphic to the left horn inclusion
  $\overinc\subseteq_r:\Lambda^0[1]->\Delta^0[1].$ and is thus a weak
  equivalence. It follows therefore, by the result of the last
  paragraph and an application of the two of three property, that the
  latter of these is also a weak equivalence. In other words, we have
  established that both of the inclusions
  $\overinc\subseteq_r:E^-_0->E_2.$ and
  $\overinc\subseteq_r:E^-_1->E_2.$ are complicial cofibrations as
  required. 
\end{obs}

\begin{lemma}\label{quasi.strat.lem}
  There exists a cofibrantly generated Quillen model structure on
  $\Strat$ whose fibrant objects are precisely the canonically
  stratified quasi-categories (cf.\ example~\ref{quasi.wcs}).
\end{lemma}

\begin{proof}
  Given the work of the last observation, it is clear that if we may
  apply theorem~\ref{model.str.thm} to the set of inclusions $J_q$
  then the resulting Quillen model structure would satisfy the
  condition postulated with respect to quasi-categories. So following
  observation~\ref{loc.mod.obs} we must show that
  $\overinc\subseteq_e:\usimp{E}_2->E_2.$ is a $J_q$-weak equivalence
  and, as discussed there, we may do so by demonstrating that
  clause~\ref{ccof.char.3} of lemma~\ref{ccof.char} holds for each
  inclusion whose domain and codomain are $J_1$-weak complicial sets.
  So assume w.l.o.g.\ that $i$ is a subset inclusion
  $\overinc\subseteq_s:X->Y.$ and that $X$ and $Y$ are $J_1$-weak
  complicial sets and consider the corner tensor:
  \begin{equation}\label{eqv.corner.tens}
    \overinc\subseteq_e:(E_2\gray X)\cup
    (\usimp{E}_2\gray Y)->E_2\gray Y.
  \end{equation}
  Now let $U$ denote the underlying simplicial set common to the
  domain and codomain of this inclusion. Since the entire inclusion
  $\overinc\subseteq_e:U->U^e.$ of definition~\ref{cestrat.not} is a
  $J_q$-anodyne extension, it is clear that we may demonstrate that
  the inclusion above has the LLP with respect to each $J_q$-weak
  complicial set simply by showing that its codomain $E_2\gray Y$ is
  an entire subset of $U^e$.  However this latter set is $1$-trivial
  so it is enough to check that each thin $1$-simplex of $E_2\gray Y$
  is also thin in $U^e$. To that end suppose that $(e,y)$ is a thin
  $1$-simplex in $E_2\gray Y$, and observe that the thin $y\in Y$
  gives rise to a corresponding Yoneda map $\arrow
  \yoneda{y}:E^-_1\cong\Delta[1]_t->Y.$ which we may lift along the
  inclusion $\overinc\subseteq_r:E^-_1->E_2.$, using the assumption
  that $Y$ is a ($J_1$-)weak complicial set and
  corollary~\ref{equiv.coher.cor1}, to give a map $\arrow
  \hat{y}:\usimp{E}_2->Y.$ with $\hat{y}(e^-_1)=y$. Similarly, the
  identity map on $E_2$, the dual map $\arrow \neg:E_2->E_2.$ of
  observation~\ref{iso.symm} and the maps which carry the whole of
  $E_2$ to the $0$-simplex $-$ or $+$ provide maps which carry the
  simplex $e^-_1$ to each one of the $1$-simplices in $E_2$, so we may
  adopt a corresponding notation
  $\arrow\hat{e}:\usimp{E}_2->\usimp{E}_2.$ for the stratified map
  with $\hat{e}(e^-_1)=e$. It follows, therefore, that $(e,y)$ is the
  image of the simplex $e^-_1$ under the induced map $\arrow
  (\hat{e},\hat{v}):\usimp{E}_2->\usimp{E}_2\gray Y\subseteq_e U.$
  thereby demonstrating that it is thin in $U^e$ as required.
\end{proof}

\begin{cor}
  There exists a cofibrantly generated Quillen model structure on
  $\Simp$ whose:
  \begin{itemize}
  \item cofibrations are the simplicial inclusions, 
  \item weak equivalences are those maps in $\Simp$ which are
    $J_q$-weak equivalences in $\Strat$ (under the minimal
    stratification), and
   \item fibrations are those $\arrow p:A->B.$ in $\Simp$ for which
    $\arrow p^e:A^e->B^e.$ is a completely $J_q$-complicial fibration
    in $\Strat$
  \end{itemize}
  In particular, the fibrant objects in this model category are
  the quasi-categories.
\end{cor}

\begin{proof}
  We construct a Quillen model structure on $\Simp$ by restricting the
  $J_q$-model structure of $\Strat$ along the fully-faithful functor
  $\arrow (-)^e:\Simp->\Strat.$ of definition~\ref{cestrat.not}. In
  other words, we define classes of cofibrations, fibrations and weak
  equivalences by saying that a simplicial map $\arrow f:X->Y.$ is a
  cofibration (resp.\ fibration or weak equivalence) if and only if
  the stratified map $\arrow f^e:X^e->Y^e.$ is a $J_q$-complicial
  cofibration (resp.\ complete $J_q$-complicial fibration or
  $J_q$-weak equivalence). Now we simply verify Quillen's axioms M1 to
  M5 (see definition 7.1.3 of \cite{Hirschhorn:2003:ModCat} for
  instance) for this choice. Axiom M1 (limits and colimits) is
  immediate for the presheaf category $\Simp$ whereas axioms M2 to M4
  (2-of-3, retract and lifting) are all immediate consequences of the
  corresponding axioms for the $J_q$-complicial model structure and
  the fact that $\arrow (-)^e:\Simp->\Strat.$ is a fully-faithful
  functor. 

  That simply leaves us to verify axiom M5 (factorisation), which
  postulates that we may factor each simplicial map $\arrow f:X->Y.$
  as a composite $f=p\circ i$ wherein $p$ is a fibration (resp.\
  trivial fibration) and $i$ is a trivial cofibration (resp.\
  cofibration). However, we know that we may factor the stratified map
  $\arrow f^e:X^e->Y^e.$ as a fibration, trivial cofibration (resp.\
  trivial fibration, cofibration) composite $\overinc
  i:X^e->W.\overarr p:->Y^e.$ in the $J_q$-complicial model structure
  and that this gives rise to an appropriate factorisation in $\Simp$
  under the proposed model structure if and only if the
  stratified set $W$ is of the form $Z^e$. Now we know, from
  the comments in definition~\ref{cestrat.not}, that this latter
  condition holds if and only if $W$ is $Q$-fibrant. Furthermore the
  same passage tells us that $Y^e$ is $Q$-fibrant and
  observation~\ref{qcat=J_q-wcs} demonstrates that the (trivial)
  fibration $p$ is a $Q$-fibration, from which facts we may infer that
  $W$ is also $Q$-fibrant as suggested.

  In the Quillen model structure we have just constructed it is clear
  that a simplicial map is a cofibration iff it is a simplicial
  inclusion. Furthermore we know, from definition~\ref{cestrat.not},
  that the inclusion $\overinc\subseteq_e:X->X^e.$ is a $J_q$-anodyne
  extension for each simplicial set $X$. So if $\arrow w:X->Y.$ is a
  simplicial map then we know that the horizontal inclusions in the
  square
  \begin{equation*}
    \xymatrix@=2em{
      {X} \ar@{u(->}[r]^-{\subseteq_e}\ar[d]_{w} &
      {X^e}\ar[d]^{w^e} \\
      {Y} \ar@{u(->}[r]_-{\subseteq_e} & {Y^e}
    }
  \end{equation*}
  are $J_q$-weak equivalences and thus that we may apply the 2-of-3
  property to show that $w^e$ if a $J_q$-weak equivalence iff $w$ is
  such in $\Strat$. In other words, we find that $w$ is a weak
  equivalence in the Quillen model structure derived in the last two
  paragraphs iff it is a $J_q$-weak equivalence as a minimally
  stratified map as postulated in the statement.
\end{proof}

\begin{notation}
  We call the Quillen model structure derived in the last corollary the {\em
    quasi-categorical model structure\/} and use the terms {\em
    quasi-cofibration\/} for its trivial cofibrations and {\em complete
    quasi-fibration\/} for its fibrations.
\end{notation}

\begin{defn}\label{qfib.defn}
  We say that a map $\arrow p:A->B.$ in $\Simp$ is an {\em inner
    quasi-fibration\/} if it has the RLP with respect to the
  simplicial inclusions
  \begin{equation}\label{simp.inner.horns}
    \{\overinc\subseteq_s:\Lambda^k[n]->\Delta[n].\mid n =
    2,3,... \wedge 0<k<n\}
  \end{equation}
  that is to say these are what Joyal calls {\em mid-fibrations}. We
  also say that $p$ is a {\em quasi-fibration\/} if it is an inner
  quasi-fibration which also has the RLP with respect to the
  simplicial inclusion
  \begin{equation}\label{simp.eq.incs}
    \overinc\subseteq_r:\usimp{E}^-_0->\usimp{E}_2.
  \end{equation}
\end{defn}

\begin{lemma}
  Each of the inclusions in displays~(\ref{simp.inner.horns})
  and~(\ref{simp.eq.incs}) of the last definition is a
  quasi-cofibration. It follows that every complete quasi-fibration is
  actually a quasi-fibration in the sense introduced there.
\end{lemma}

\begin{proof}
  We know, from the last theorem, that it is enough to show that the
  stratified maps obtained by applying the functor $({-})^e$ to the
  simplicial inclusions in the cited displays are all $J_q$-weak
  equivalences. However only the degenerate simplices of
  $\triv_1(\Delta[n])$ have equivalence inverses, so if we apply the
  functor $({-})^e$ to the inner horn inclusion
  $\overinc\subseteq_s:\Lambda^k[n]->\Delta[n].$ then we obtain a
  stratified inclusion which may otherwise be constructed by applying
  $\arrow\triv_1:\Strat->\Strat.$ to the inner complicial horn
  inclusion $\overinc\subseteq_r:\Lambda^k[n]->
  \Delta^k[n].$. Applying lemma~\ref{super.wcs} it follows that this
  inclusion is an inner anodyne extension and is thus also a
  ($J_q$-)weak equivalence.  It is also clear that the simplicial maps
  $\arrow\id:\usimp{E}_2->\usimp{E}_2.$ and
  $\arrow\neg:\usimp{E}_2->\usimp{E}_2.$ demonstrate that the
  $1$-simplices $e^-_1$ and $e^+_1$ are thin in $(\usimp{E}_2)^e$, so
  when we apply $(-)^e$ to the inclusion
  $\overinc\subseteq_r:\usimp{E}^-_0-> \usimp{E}_2.$ we obtain the
  stratified inclusion $\overinc\subseteq_r:E^-_0->E_2.$ which was
  shown to be a ($J_q$-)weak equivalence in
  observation~\ref{qcat=J_q-wcs}.
\end{proof}

\begin{lemma}
  If the simplicial map $\arrow p:A->B.$ is a quasi-fibration between
  quasi-categories then it is a complete quasi-fibration.
\end{lemma}

\begin{proof}
  Since $A$ and $B$ are quasi-categories we know, by
  observation~\ref{qcat=J_q-wcs}, that the stratified sets $A^e$ and
  $B^e$ are $J_q$-weak complicial sets. It follows immediately, by the
  comment in definition~\ref{thinness.ext}, that $\arrow
  p^e:A^e->B^e.$ has the RLP with respect to each elementary thinness
  extension $\overinc\subseteq_e:\Delta^k[n]'->\Delta^k[n]''.$, each
  thin simplex inclusion $\overinc\subseteq_e:\Delta[n]->\Delta[n]_t.$
  ($n>1$) and the equivalence inclusion $\overinc\subseteq_e:
  \usimp{E}_2->E_2.$. Furthermore, it is easily seen that $p^e$ has
  the RLP with respect to the inner horn inclusions
  $\overinc\subseteq_r:\Lambda^k[n]->\Delta^k[n].$ ($0<k<n$) and the
  equivalence inclusion $\overinc\subseteq_r:E^-_0->E_2.$ iff $p$ has
  the RLP with respect to their underlying simplicial maps. These are,
  however, the simplicial inclusions used to describe quasi-fibrations
  in definition~\ref{qfib.defn}, so it follows from the postulated
  properties of $p$ that $p^e$ does indeed have the RLP with respect
  to the stratified inclusions of the last sentence.

  To summarise we have shown that, under the conditions of the
  statement, the stratified map $\arrow p^e:X^e->Y^e.$ is an inner
  complicial fibration between $J_q$-weak complicial sets and that it
  also has the RLP with respect to the inclusion
  $\overinc\subseteq_r:E^-_0->E_2.$. Consequently
  lemma~\ref{equiv.coher.cor2} tells us that $p^e$ is a
  ($J_q$-)complicial fibration iff it has the RLP with respect to the
  inclusion $\overinc\subseteq_r:E^-_0->E^-_1.$, a result we establish
  by constructing a solution to the arbitrary lifting problem $(u,v)$
  depicted in the following square:
  \begin{equation*}
    \xymatrix@R=0.5em@C=2em{
      {E^-_0}\ar[rr]^-{u}\ar@{u(->}[dd]_{\subseteq_r} &&
      {X^e}\ar[dd]^{p^e} \\
      & {E_2}\ar@{->}[dr]^{v'}\ar@{->}[ur]_{w} & \\
      {E^-_1}\ar[rr]_-{v}\ar@{u(->}[ru]^-{\subseteq_r}
      && {Y^e}}
  \end{equation*}
  Start by factoring the map $v$ through the inclusion
  $\overinc\subseteq_r:E^-_1->E_2.$, which we may do since we know
  that this inclusion is a complicial cofibration
  (observation~\ref{qcat=J_q-wcs}) and that $Y^e$ is a weak complicial
  set, to give the map $v'$. Now we have a lifting problem $(u,v')$
  from the inclusion $\overinc\subseteq_r:E^-_0->E_2.$ to $p^e$ for
  which we may find a solution $w$ since we know that the latter map
  has the RLP with respect to the former. Finally we may take the
  composite of $w$ and the inclusion $\overinc\subseteq_r:E^-_1->E_2.$
  as the required solution to our original problem. 

  All that remains now is to apply lemma~\ref{cfib+wcs=>ccfib} to show
  that $p^e$ is actually a complete $J_q$-complicial fibration and
  thus that, by definition, $p$ is a complete quasi-fibration as
  postulated.
\end{proof}



\section{Appendix A - Some Categorical Homotopy Theory}
\label{mor.append}

We recollect here a few basic results of categorical homotopy theory
upon which we rely in the body.

\begin{defn}[categories of morphisms]
  If $\lcat{C}$ is a category then its {\em category of morphisms\/}
  $\lcat{C}^2$ is defined to be the category of functors from the
  ordinal $2 = \{0<1\}$ to $\lcat{C}$. Its objects are simply
  morphisms $\arrow f:A->B.$ of $\lcat{C}$ and its arrows from $f$ to
  another morphism $\arrow g:C->D.$ are pairs of arrows $(u,v)$ of
  $\lcat{C}$ making the obvious naturality square
  \begin{equation*}
    \xymatrix@=2em{
      {A}\ar[r]^u\ar[d]_f & {C}\ar[d]^g \\
      {B}\ar[r]_v & {D}
    }
  \end{equation*}
  commute and which are thus called {\em squares}. In the context of
  Quillen model categories, we often think of the arrows of
  $\lcat{C}^2$ as being {\em lifting problems\/} in $\lcat{C}$.
\end{defn}

\begin{recall}[the corner tensor and its closures]\label{corner.tensor}
  It is common in the theory of Quillen model categories to consider
  {\em pushout corner maps\/} so we recall the basic concepts and
  notation here in a suitably generalised setting. 

  Let $\lcat{C}$, $\lcat{D}$ and $\lcat{E}$ be categories which are
  cocomplete and let $\arrow\anytens:\lcat{C}\times
  \lcat{D}->\lcat{E}.$ be a bifunctor (tensor) which preserves these
  colimits in each variable. Now suppose that $\arrow f:C->C'.$ is an
  arrow of $\lcat{C}$ and $\arrow g:D->D'.$ is an arrow of $\lcat{D}$
  then we may consider the commutative square
  \begin{equation*}
    \xymatrix@=2em{
      {C\anytens D} \ar[r]^{f\anytens D} \ar[d]_{C\anytens g} &
      {C'\anytens D} \ar[d]^{C'\anytens g} \\
      {C\anytens D'} \ar[r]_{f\anytens D'} & {C'\anytens D'} }
  \end{equation*}
  which induces a unique map usually denoted $f\canytens g$ from the
  pushout $(C'\anytens D)\vee_{C\anytens D}(C\anytens D')$ of the
  upper horizontal and left hand vertical maps in this square to its
  lower right vertex $C'\anytens D'$ making the usual triangles
  commute. This map is often called the {\em corner tensor\/} (or
  sometimes the {\em Liebnitz tensor}) of $f$ and $g$ and, for
  instance, it plays a central role in Quillen's theory of simplicial
  model categories \cite{Quillen:1967:Model} (for a suitable
  $\anytens$). This construction provides us with a naturally defined
  bifunctor $\arrow\canytens:\lcat{C}^2 \times\lcat{D}^2->\lcat{E}^2.$
  which again preserves colimits in each variable.

  Generally $\anytens$ will be closed in each variable, meaning that for
  each $C\in\lcat{C}$ the functor $\arrow
  C\anytens{-}:\lcat{D}->\lcat{E}.$ has a right adjoint
  $\arrow\cls_l(C,{*}):\lcat{E}->\lcat{D}.$ (and
  dually for objects $D\in\lcat{D}$). In this case, the corner
  tensor $\canytens$ is also closed in each variable with the (left)
  {\em corner closure\/} $\ccls_l(f,h)$ of morphisms $\arrow
  f:C->C'.\in\lcat{C}$ and $\arrow h:E->E'.\in\lcat{E}$ being
  the unique map induced by the commutative square
  \begin{equation*}
    \xymatrix@R=2em@C=5em{
      {\cls_l(C',E)} \ar[r]^{\cls_l(f,E)}
      \ar[d]_{\cls_l(C',h)} &
      {\cls_l(C,E)} \ar[d]^{\cls_l(C,h)} \\
      {\cls_l(C',E')} \ar[r]_{\cls_l(f,E')} &
      {\cls_l(C,E')}}
  \end{equation*}
  from its upper left vertex to the pullback
  $\cls_l(C,E)\times_{\cls_l(C,E')}\cls_l(C',E')$
  of it right vertical and lower horizontal maps.

  Most importantly, the corner tensor is well behaved with respect
  to cellular completions of sets of morphisms:
\end{recall}

\begin{lemma}\label{corner.cofibration}
  Let $I$ and $J$ be sets of morphisms of $\lcat{C}$ and $\lcat{D}$
  and let $\mathcal{K}$ be a class of morphisms of $\lcat{E}$ which is
  closed under pushout and transfinite composition. In particular, we
  may also take $\mathcal{K}$ to be $\cell(K)$ for some set of
  morphisms $K$ in $\lcat{E}$.

  Suppose also that we know that whenever $i\in I$ and $j\in J$ then
  their corner tensor $i\canytens j$ is in $\mathcal{K}$. Then
  whenever $f$ is a morphism in $\cell(I)$ and $g$ is a morphism in
  $\cell(J)$ we may infer that their corner tensor $f\canytens g$ is in
  $\mathcal{K}$.
\end{lemma}

\begin{proof}
  The proof here is entirely standard and is left to the reader.
\end{proof}

\begin{obs}\label{corner.fibration}
  On interpreting the arrows of $\lcat{E}^2$ and $\lcat{D}^2$ as
  lifting problems in $\lcat{E}$ and $\lcat{D}$ (respectively) it is
  worth observing that if we are given $\arrow f:C->C'.\in\lcat{C}$,
  $\arrow g:D->D'.\in\lcat{D}$ and $\arrow h:E->E'.\in\lcat{E}$ then
  the adjunction $f\canytens{-}\dashv\ccls_l(f,*)$ sets up a bijection
  between lifting problems
  \begin{equation*}
    \vcenter{\xymatrix@R=2em@C=2em{
        {(C'\anytens D)\vee_{C\anytens D}(C\anytens D')}
        \ar[r]^-{u}\ar[d]_{f\canytens g} &
        {E}\ar[d]^{h} \\
        {C'\anytens D'}\ar[r]_-{v}\ar@{..>}[ur]_-{l} & {E'}
      }}\mkern30mu
    \vcenter{\xymatrix@R=2em@C=2em{
        {D}\ar[r]^-{u'}\ar[d]_{g} &
        {\cls_l(C,E)}\ar[d]^{\ccls_l(f,h)} \\
        {D'}\ar[r]_-{v'}\ar@{..>}[ur]^-{\hat{l}} & 
        {\cls_l(C,E)\times_{\cls_l(C,E')}\cls_l(C',E')}
      }}
  \end{equation*}
  in $\lcat{E}$ and $\lcat{D}$ respectively. Furthermore, as indicated
  a map $\arrow l:C'\anytens D'->E.$ in $\lcat{E}$ is a solution of
  the lifting problem on the left iff the dual map $\arrow
  \hat{l}:D'->\cls_l(C',E).$ in $\lcat{D}$ under the adjunction
  $C'\anytens{-}\dashv\cls_l(C',{*})$ is a solution of the dual
  lifting problem on the right. It follows, therefore, that under the
  conditions of the last lemma if $h$ is a $K$-fibration and $f$ is an
  $I$-cofibration then $\ccls_l(f,p)$ is a $J$-fibration.
\end{obs}

\begin{obs}[the small object argument]\label{small.obj.arg}
  Almost all constructions of Quillen model structures rely upon some
  version of Quillen's small object argument. For instance the proof
  of Jeff Smith's construction theorem rests upon a variant of this
  construction presented in subsection~III.6
  of~\cite{Adamek:1993:Inj}. Explicitly, if $J$ is a (small) set of
  morphisms of our locally presentable category $\lcat{C}$ then
  proposition~III.8 of that work allows us to construct a {\em weak
    reflection\/} of $\lcat{C}$ into its full subcategory $\lcat{C}_J$
  of $J$-injective ($J$-fibrant) objects (called the {\em injectivity
    class associated with $J$}).  We recall, and slightly recast, their
  construction here in order to extract a few of the properties of the
  resulting weak reflection which are not discussed explicitly in
  loc.\ cit.

  We will assume that we are given a (small) set of morphisms $J$ in a
  locally presentable category $\lcat{C}$ and adopt the notation $U_j$
  and $V_j$ for the domain and codomain of a morphism $j\in J$
  (respectively). We also assume, by an appeal to corollary~2.3.12
  of~\cite{Makkai:1989:Access}, that we have chosen a regular cardinal
  $\kappa$ for which $\lcat{C}$ is locally $\kappa$-presentable and
  for which the domains $U_j$ and codomains $V_j$ of the elements of
  $J$ are all $\kappa$-presentable. Now we start by constructing a
  pointed endo-functor $(F,\phi)$ on $\lcat{C}$ by forming a
  (pointwise) pushout
  \begin{equation*}
    \xymatrix@R=2em@C=10em{
      {\coprod_{j\in J}\lcat{C}(U_j,-)\tens U_j}
      \ar[r]^{\coprod_{j\in J}\lcat{C}(U_j,-)\tens j}\ar[d] &
      {\coprod_{j\in J}\lcat{C}(U_j,-)\tens V_j}\ar[d] \\
      {\id_{\lcat{C}}}\ar[r]_{\phi} & {F}\poexcursion}
  \end{equation*}
  in the endo-functor category $[\lcat{C},\lcat{C}]$. Here we use
  $X\tens W$ to denote the $X$-fold coproduct of $W$ with
  itself. So it is clear that the component $\overarr:\coprod_{j\in J}
  \lcat{C}(U_j,X)\bullet U_j->X.$ of the left hand vertical in this
  square is naturally defined to be the map induced by the family
  whose component from the copy of $U_j$ corresponding to some
  $f\in\lcat{C}(U_j,X)$ is simple the morphism $f$ itself.  Notice
  that each representable $\arrow \lcat{C}(U_j,-):\lcat{C}->Set.$ is
  {\em $\kappa$-accessible\/} (preserves $\kappa$-filtered colimits),
  since $U_j$ is $\kappa$-presentable, and that the tensor $-\tens W$
  preserves all colimits, so it follows that each functor
  $\lcat{C}(U_j,-)\tens W$ is $\kappa$-accessible. Now the full
  subcategory of $\kappa$-accessible endo-functors is closed under
  colimits, since colimits commute with colimits, so consequently $F$
  is also $\kappa$-accessible since it is a pushout of coproducts of
  $\kappa$-accessible functors. To complete their construction, we
  iterate $F$ to obtain a transfinite chain of powers of $F$ with
  \begin{align*}
    F^0 & {} \defeq \id_{\Strat} && \\
    F^{\alpha^+} & {}\defeq F\circ F^\alpha && \text{at successor
      ordinals $\alpha^+$}\\
    F^{\gamma} & {}\defeq \colim_{\alpha<\gamma}(F^\alpha) && \text{at
      limit ordinals $\gamma$}
  \end{align*}
  and chain maps $\arrow\phi_{\alpha,\beta}:F^\alpha->F^\beta.$ (for
  $\alpha\leq\beta$) determined by:
  \begin{align*}
    \phi_{\alpha,\alpha^+} & \defeq \arrow \phi\circ
    F^\alpha:F^{\alpha}->F\circ F^\alpha=F^{\alpha^+}. &&
    \text{between an ordinal and its successor}\\
    \phi_{\alpha,\gamma} & \defeq \arrow
    \iota^\gamma_\alpha:F^{\alpha}->\colim_{\alpha<\gamma}
    (F^{\alpha})=F^\gamma. && \text{the canonical colimit inclusion.}
  \end{align*}
  Now using the fact that the all of the objects $U_j$ and $V_j$ are
  $\kappa$-presentable in $\lcat{C}$ we may apply proposition~III.8 of
  loc.\ cit.\ to show that $F^\kappa$ is the desired weak reflection,
  which is $\kappa$-accessible since the class of $\kappa$-accessible
  functors is closed under composition and colimits.
  
  Of course the category of morphisms $\lcat{C}^2$ is also locally
  $\kappa$-presentable; its limits and colimits are formed
  pointwise in $\lcat{C}$ and its $\kappa$-presentable
  objects are those morphisms $\arrow f:A->B.$ for which $A$ and $B$
  are both $\kappa$-presentable in $\lcat{C}$. So we may apply the weak
  reflection result above to the set $J_s$ of squares in $\lcat{C}^2$
  of the form
  \begin{equation*}
    \xymatrix@=2em{
      {U_j}\ar[r]^j\ar[d]_j & {V_j}\ar[d]^{\id_{V_j}} \\
      {V_j}\ar[r]_{\id_{V_j}} & {V_j}
    }
  \end{equation*}
  for each arrow $\arrow j:U_j->V_j.$ in $J$, to thereby construct a
  functorial and $\kappa$-accessible factorisation of any arrow
  $\arrow f:A->B.$ of $\lcat{C}$ into a composite $p\circ k$ in which
  $p\in\fib(J)$ and $k\in\cell(J)$.
\end{obs}

\begin{obs}[injectivity classes and accessibility]\label{inj.acc}
  Using the fact that the domains of the maps in $J$ are all
  $\kappa$-presentable, it is easily shown that the injectivity class
  $\lcat{C}_J$ associated with $J$ is closed in $\lcat{C}$ under
  $\kappa$-filtered colimits. Furthermore, applying corollary~III.9
  of~\cite{Adamek:1993:Inj} and the work of subsection~2.3
  of~\cite{Makkai:1989:Access}, we may construct a regular cardinal
  $\nu>\kappa$ for which $\lcat{C}$ is locally $\nu$-presentable and
  for which the weak reflection of the last observation carries each
  $\nu$-presentable object $A\in\lcat{C}$ to an object $F^\kappa(A)$
  which is also $\nu$-presentable in $\lcat{C}$.

  Using this property, it is easily demonstrated that if
  $\lcat{C}_\nu$ is the essentially small, full subcategory of
  $\nu$-presentable objects in $\lcat{C}$ and $A$ is an arbitrary
  $J$-injective then the comma category
  $(\lcat{C}_J\cap\lcat{C}_\nu)\downarrow A$ is $\nu$-filtered and
  cofinal in $\lcat{C}_\nu\downarrow A$. Now since $\lcat{C}$ is
  locally $\nu$-presentable we know that $A$ is the colimit of the
  canonical diagram $\arrow D_A:\lcat{C}_\nu\downarrow A->\lcat{C}.$
  so we may infer, from the last sentence, that $A$ is also the
  colimit in $\lcat{C}_J$ of the restricted diagram $\arrow
  D_A:(\lcat{C}_J\cap \lcat{C}_\nu)\downarrow
  A->\lcat{C}_J.$. Furthermore it is clear that $\lcat{C}_J\cap
  \lcat{C}_\nu$ is essentially small and that each of its objects is
  $\nu$-presentable in $\lcat{C}_J$ so it follows, by definition, that
  $\lcat{C}_J$ is $\nu$-accessible.

  Of course, we may apply the result above to the locally
  $\kappa$-presenable category of morphisms $\lcat{C}^2$ and its set
  $J_s$ of squares derived from $J$ as in the final paragraph of the
  last observation. Doing so we find that the class $\fib(J)$ of
  $J$-fibrations is always an {\em accessible class of maps\/} in
  $\lcat{C}$. In other words, the corresponding full subcategory of
  $\lcat{C}^2$ whose objects are the $J$-fibrations is both
  ($\nu$-)accessible and ($\kappa$-)accessibly embedded in
  $\lcat{C}^2$, simply because it may otherwise be described as the
  injectivity class associated with $J_s$.
\end{obs}

\begin{thm}[Jeffery Smith's theorem]\label{jsmith.thm}
  Let $\lcat{C}$ be a locally presentable category, $\mathcal{W}$ be a
  subclass of its morphisms and $I$ be a small set of its
  morphisms. Suppose further that they satisfy the following
  conditions:
  \begin{enumerate}[label=(\arabic*)]
  \item\label{jsmith.thm.1} $\mathcal{W}$ is closed under retracts and
    has the 2-of-3 property.
  \item\label{jsmith.thm.2} $\fib(I)$ is a subclass of $\mathcal{W}$.
  \item\label{jsmith.thm.3} The class $\cof(I)\cap\mathcal{W}$ is
    closed under transfinite composition and under pushout.
  \item\label{jsmith.thm.4} $\mathcal{W}$ satisfies the {\em solution
      set condition\/} at $I$.
  \end{enumerate}
  Then taking $\mathcal{W}$ as our class of {\bfseries weak
    equivalences}, $\cof(I)$ as our class of {\bfseries
    cofibrations} and $\fib(\cof(I)\cap\mathcal{W})$ as our class of
  {\bfseries fibrations} we obtain a cofibrantly generated Quillen
  model structure on $\lcat{C}$.
\end{thm}

\begin{proof}
  A discussion of the technical details and a full proof of this
  (folkloric) result may be found in Beke's work on simplicial
  sheaves~\cite{Beke:2000:SheafHomotopy}.  When we apply this theorem
  herein we will rely on the fact that our $\mathcal{W}$ is an
  accessible class of maps, which condition then ensures that
  condition~\ref{jsmith.thm.4} holds for any set of maps $I$.
\end{proof}


\bibliographystyle{amsplain}
\bibliography{cattheory}

\end{document}